\documentclass[11pt]{amsart}

\usepackage{amssymb, amsmath, amscd, amsthm,
color, epsfig, graphicx}

\usepackage[all]{xy}          
\xyoption{dvips}              

\newcommand{\R}{{\mathbb R}}
\newcommand{\Z}{{\mathbb Z}}

\newcommand{\C}{{\mathbb C}}

\newcommand{\HH}{\operatorname{H}} 
\newcommand{\Cyl}{\operatorname{Cyl}}
\newcommand{\Emb}{\operatorname{Emb}}
\newcommand{\sEmb}{\operatorname{sEmb}}
\newcommand{\sImm}{\operatorname{sImm}}
\newcommand{\sEbar}{\overline{\sEmb}}
\newcommand{\nlc}{\operatorname{nlc}}
\newcommand{\Ebar}{\overline{\Emb}}

\newcommand{\Id}{\operatorname{Id}}
\newcommand{\calJ}{{\mathcal J}}

\newcommand{\Imm}{\operatorname{Imm}}
\newcommand{\Spectra}{\operatorname{Spectra}}
\newcommand{\Ob}{{\mathcal O}}

\newcommand{\map}{\operatorname{Map}}
\newcommand{\hnat}{\operatorname{hNat}}
\newcommand{\nat}{\operatorname{Nat}}
\newcommand{\holim}{\operatorname{holim}\,}

\newcommand{\hofiber}{\operatorname{hofiber}}

\newcommand{\D}{\operatorname{D}}
\newcommand{\TP}{\operatorname{P}}

\newcommand{\homotopy}{\operatorname{h}}

\newcommand{\comp}{\operatorname{c}}
\newcommand{\supp}{\operatorname{s}}
\newcommand{\excess}{\operatorname{e}}
\newcommand{\Top}{\operatorname{Top}}

\newcommand{\calC}{\mathcal{C}}
\newcommand{\calD}{\mathcal{D}}
\newcommand{\calE}{\mathcal{E}}
\newcommand{\calP}{\mathcal{P}}
\newcommand{\calS}{\mathcal{S}}
\newcommand{\calF}{\mathcal{F}}
\newcommand{\calMP}{\mathcal{MP}}
\newcommand{\calO}{\mathcal{O}}

\newcommand{\GL}{\operatorname{GL}}

\newcommand{\config}{\operatorname{C}}

\newcommand{\OO}{\operatorname{O}}
\newcommand{\Mor}{\mathcal{M}}
\newcommand{\Nerve}{\operatorname{N}}

\newcommand{\iso}{\operatorname{Iso}}
\newcommand{\ad}{\operatorname{Ad}}

\newcommand{\op}{\operatorname{op}}

\theoremstyle{plain}
\newtheorem{theorem}{Theorem}[section]
\newtheorem{proposition}[theorem]{Proposition}
\newtheorem{lemma}[theorem]{Lemma}
\newtheorem{corollary}[theorem]{Corollary}

\theoremstyle{definition}
\newtheorem{definition}[theorem]{Definition}
\newtheorem{example}[theorem]{Example}

\theoremstyle{remark}
\newtheorem{remark}[theorem]{Remark}

\begin{document}


\title{Derivatives of Embedding Functors I: The Stable Case.}


\author{Gregory Arone}
\address{Department of Mathematics, University of Virginia,
Charlottesville, VA} \email{zga2m@virginia.edu}
\thanks{The author gratefully acknowledges that he was supported by the NSF during the preparation of this paper, most recently via  
grant DMS 0605073}



\begin{abstract}
For smooth manifolds $M$ and $N$, let $\Ebar(M, N)$ be the homotopy fiber of the map $\Emb(M, N)\longrightarrow \Imm(M, N)$. Consider the functor from the category of Euclidean spaces to the category of spectra, defined by the formula $V\mapsto \Sigma^\infty\Ebar(M, N\times V)$. In this paper, we describe the Taylor polynomials of this functor, in the sense of M. Weiss' orthogonal calculus, in the case when $N$ is a nice open submanifold of a Euclidean space. This leads to a description of the derivatives of this functor when $N$ is a tame stably parallelizable manifold (we believe that the parallelizability assumption is not essential). Our construction involves a certain space of rooted forests (or, equivalently, a space of partitions) with leaves marked by points in $M$, and a certain ``homotopy bundle of spectra'' over this space of trees. The $n$-th derivative is then described as the ``spectrum of restricted sections'' of this bundle. This is the first in a series of two papers. In the second part, we will give an analogous description of the derivatives of the functor $\Ebar(M, N\times V)$, involving a similar construction with certain spaces of connected graphs (instead of forests) with points marked in $M$.
\end{abstract}

\maketitle


\section{Introduction} Let $M$, $N$ be smooth finite-dimensional manifolds. Let $\Emb(M,N)$ and $\Imm(M,N)$ be the space of smooth embeddings and the space of smooth immersions of $M$ into $N$, respectively. Let $\Ebar(M,N)$ be the homotopy fiber of the inclusion map $\Emb(M,N)\hookrightarrow \Imm(M,N)$. Our goal is to study the space $\Ebar(M,N)$ using Michael Weiss's orthogonal calculus~\cite{WeissOrth}. To be more specific, consider the functor from the category of Euclidean spaces to the category of spectra defined by the formula
$$V\mapsto \Sigma^\infty \Ebar(M, N\times V).$$
This is a functor to which one can apply orthogonal calculus. The Taylor tower of this functor was studied in~\cite{ALV} in the special case when $N=*$ (from a perspective different than the one taken here). The main result there says that the Taylor tower of the functor $V\mapsto \Sigma^\infty \Ebar(M, V)$ rationally splits as a product of its homogeneous layers when the dimension of $V$ is more than twice the embedding dimension of $M$. Furthermore, the derivatives of the functor $\Sigma^\infty\Ebar(M,V)$ were described, without proof, and certain homotopy invariance properties of the derivatives, together with the aforementioned splitting result, were used to prove a statement about the rational homology invariance of $\Ebar(M,V)$. The goal of this paper is to describe, with proofs, the derivatives of the functor $V\mapsto \Sigma^\infty \Ebar(M, N\times V)$. We will do so in the case when $M$ and $N$ are  tame manifolds, and $N$ is also stably parallelizable. Here by ``tame'' we mean that it is diffeomorphic to the interior of a compact manifold with a (possibly empty) boundary. In fact, we will, in some sense,  describe the Taylor polynomials (as opposed to just the layers) of this functor in the case when $N$ is a nice open submanifold of a Euclidean space.

Before we can state our main result (Theorem~\ref{T: MainTheorem} below), we have to review some background material, and to introduce some definitions. We will assume that the reader is familiar with the general idea of calculus of functors, including orthogonal calculus. The basic reference for orthogonal calculus is~\cite{WeissOrth}.

Let $\calJ$ be the category of Euclidean spaces and linear isometric inclusions. Let $\calD$ be a pointed topological model category, e.g., the category of pointed spaces or the category of spectra. Orthogonal calculus is a framework for studying continuous functors from $\calJ$ to $\calD$. Let $F\colon \calJ\longrightarrow \calD$ be such a functor. Orthogonal calculus associates with $F$ a sequence, which we will denote $\partial_1F,  \partial_2F,\ldots, \partial_nF,\ldots $, where $\partial_nF$ is a spectrum with an action of the orthogonal group $\OO( n)$. We will call such a sequence an {\em orthogonal sequence} of spectra. If $F$ is a spectrum-valued functor, the $n$-th homogeneous layer in the Taylor tower of $F$ is determined by $\partial_nF$ by the formula $$\D_nF(V)=\left(\partial_nF\wedge S^{nV}\right)_{\homotopy\OO(n)}.$$
Thus the sequence of derivatives of $F$ potentially contains a lot of information about the homotopy type of $F$.

In this paper, we want to use orthogonal calculus to study covariant functors on manifolds. We will do it in the following way. Let $F$ be a continuous isotopy functor from the category of smooth manifolds and smooth embeddings to the category of pointed spaces (or spectra). We can apply orthogonal calculus to $F$ by considering the functor $V\mapsto F(N\times V)$, where $N$ is a fixed manifold and $V$ ranges over the category of Euclidean spaces. This is a functor to which orthogonal calculus applies. Its Taylor tower, evaluated at $V=\R^0$, is a tower of functors approximating $F(N)$, and under favorable circumstances converging to $F(N)$. In this paper, we call this tower ``the Taylor tower of $F$''. Similarly, we say that a functor on manifolds $F$ is polynomial (or homogeneous), if the associated functor on vector spaces $V\mapsto F(N\times V)$ is polynomial (or homogeneous) in the sense of \cite{WeissOrth}, for all manifolds $N$. 
\begin{example}
Let $N$ be a manifold of dimension $d$. Let $N_+\tilde\wedge S^d$ be the Thom space of the tangent bundle of $N$. The functor
$$N\mapsto \Sigma^\infty N_+\wedge (N_+\tilde\wedge S^d)$$ is homogeneous linear in our sense, because 
$$\Sigma^\infty (N\times V)_+\wedge ((N\times V)_+\tilde\wedge S^{d+V})\simeq \Sigma^\infty N_+\wedge (N_+\tilde\wedge S^{d})\wedge S^V.$$ 
Note that this functor is not a linear functor of $N$ in the sense of Goodwillie's calculus of homotopy functors. If anything, it looks like a quadratic functor from that perspective (we say ``looks like'', because our example is not a homotopy functor, so it is not clear that Goodwillie's calculus applies to it).
\end{example}

We call the derivatives of the functor $V\mapsto F(N\times V)$  simply the derivatives of $F$. Thus in our setting the derivatives of $F$, for a fixed functor $F$, can be thought of as an isotopy functor from the category of manifolds to the category of orthogonal sequences of spectra
$$N\mapsto \{\partial_nF(N\times -)\}_{n\ge 1}.$$

The prime example to which we would like to apply this approach is the functor $$N\mapsto \Sigma^\infty\Ebar(M, N),$$ where $M$ is considered a fixed manifold. Note that the definition of $\Ebar(M,N)$ as a space requires choosing a basepoint 
in $\Imm(M,N)$, and if one wants $\Ebar(M,N)$ to be a pointed space, then one needs to choose a basepoint in $\Emb(M,N)$. Let $h\colon M\longrightarrow N$ be the background map. When convenient, $h$ can be assumed to be an immersion, or even an embedding. However, our constructions are defined if $h$ is any continuous map, and indeed they are determined (up to homotopy) by the  homotopy class of $h$.

Our description of the layers of $\Sigma^\infty\Ebar(M,N)$ utilizes certain spaces of partitions. For us, a {\em partition} $\Lambda$ is an equivalence relation defined on a finite set $s$. We call $s$ the {\em support} of $\Lambda$. We will need to introduce several notions of morphisms between partitions. Suppose $\Lambda$ is a partition of $s$, and let $s'$ be another finite set. Let $f\colon s\longrightarrow s'$ be a map of sets. We let $f(\Lambda)$ denote the equivalence relation on $s'$ that is generated by the image of $\Lambda$ under $f$. Let $\Lambda, \Lambda'$ be partitions of $s$ and $s'$ respectively. A {\em fusion} of $\Lambda$ into $\Lambda'$ is a map of sets $f\colon s\longrightarrow s'$ such that $\Lambda'=f(\Lambda)$. It is easy to see that the composition of fusions is a fusion, and so partitions with fusions form a category.

Let $\Lambda$ be a partition of $s$, and let $c$ be the set of components (equivalence classes) of $\Lambda$. It is often convenient to represent $\Lambda$ by the surjective map of sets $\alpha\colon s\twoheadrightarrow c$ sending each element of $s$ to its component. The surjective map $\alpha$ induces an injective homomorphism $\Z^c\hookrightarrow \Z^s$. Let $\Z^s/\Z^c$ be the quotient group. It is a free abelian group of rank $|s|-|c|$. The number $|s|-|c|$ is an important (to us) invariant of partitions. We call it the {\em excess} of $\Lambda$, and denote it by $\excess(\Lambda)$. 
Let $\Lambda'$ be another partition, represented by a surjection $\alpha'\colon s'\twoheadrightarrow c'$. Let $f\colon \Lambda\longrightarrow \Lambda'$ be a fusion. The underlying map of sets induces a homomorphism $f^\sharp\colon\Z^{s'}\longrightarrow \Z^s$. It is not hard to see that $f^\sharp$ passes to a well-defined homomorphism $f^*\colon\Z^{s'}/\Z^{c'}\longrightarrow \Z^{s}/\Z^{c}$. This homomorphism is, in fact, always injective. We say that $f$ is a {\em strict fusion} if $f^*$ is an isomorphism. In particular, strict fusions preserve excess. An example of a strict fusion is what we call an {\em elementary strict fusion}, which is a fusion of $\Lambda$ into $\Lambda'$ that glues together two points in different components of $\Lambda$, and does nothing else. It can be shown that a fusion is strict if and only if it is a composite of elementary strict fusions (Lemma~\ref{L: StrictFusionsDecompose}).

We say that a partition is {\em irreducible} if none of its components is a singleton. For $n\ge 1$, let $\calE_n$ be the category of irreducible partitions of excess $n$, and strict fusions between them.
The category $\calE_n$ is instrumental in our description of the $n$-th homogeneous layer of $\Sigma^\infty\Ebar(M,N)$. An explicit description of $\calE_n$ for $n=1, 2$ is given Example~\ref{Ex: TheCategoryE_n}, and there is a more general discussion of the structure of $\calE_n$ in Remark~\ref{R: StructureOfE_n}. 

Next, we will define two functors on $\calE_n$, one covariant and one contravariant. The first functor has to do with posets of partitions. Let $\Lambda, \Lambda'$ be two partitions {\em of the same set} $s$. We say that $\Lambda'$ is a {\em refinement} of $\Lambda$ (or, equivalently, that $\Lambda$ is a coarsening of $\Lambda'$) if every component of $\Lambda'$ is a subset of some component of $\Lambda$. In this case we also write that $\Lambda\le\Lambda'$. The relation of refinement is a partial ordering on the set of partitions of $s$. Let $\calP_n$ be the poset of partitions of the standard set with $n$-elements. $\calP_n$ has both an initial and a final object, which we denote $\hat0$ and $\hat1$ (the indiscrete and the discrete partition respectively). Therefore the geometric realization of $\calP_n$, which we denote by $|\calP_n|$, is contractible (for two reasons, as it were). Inside the simplicial nerve of $\calP_n$, consider the simplicial subset consisting of those simplices that do not contain both $\hat0$ and $\hat1$ as a vertex. We denote the geometric realization of this simplicial subset by $\partial|\calP_n|$. It is a sub-complex of $|\calP_n|$, and it does geometrically look like the boundary of this contractible polyhedron. Define $T_n$ to be the quotient space $|\calP_n|/\partial|\calP_n|$. Obviously, $T_n$ has a natural action of $\Sigma_n$. It is well-known that non-equivariantly, 
$T_n\simeq \underset{(n-1)!}{\bigvee}S^{n-1}.$ 
\begin{remark}
This is not the first time that the spaces $T_n$ play a role in calculus of functors. The Spanier-Whitehead dual of $T_n$ is the $n$-th Goodwillie derivative of the identity functor on the category $\Top_*$. Note however the following difference: previously, the spaces $T_n$ occurred in the derivatives of a functor with values in $\Top_*$, while this time they occur in the derivatives of a functor with values in suspension spectra, namely $\Sigma^\infty \Ebar(M, N)$. 
\end{remark}
Now let us generalize the construction $\calP_n$ as follows: for a partition $\Lambda$ of $s$,  let $\calP(\Lambda)$ be the poset of all partitions of $s$ that are refinements of $\Lambda$. Again, $\calP(\Lambda)$ has both an initial object and a final object. Define the subcomplex $\partial|\calP(\Lambda)|\subset|\calP(\Lambda)|$ analogously to $\partial|\calP_n|$, and let 
$$T_\Lambda=|\calP(\Lambda)|/\partial|\calP(\Lambda)|.$$
Suppose that $\Lambda$ has components of sizes $n_1,\ldots, n_i$. It is not hard to see that in this case there is a homeomorphism
$$T_\Lambda\cong T_{n_1}\wedge\ldots\wedge T_{n_i}$$
which describes the homotopy type of $T_\Lambda$ in the general case. Note that $T_\Lambda$ is homotopy equivalent to a wedge of spheres of dimension $\excess(\Lambda)$.

Now let $f\colon \Lambda\longrightarrow \Lambda'$ be a morphism in $\calE_n$. It is easy to see that $f$ induces a map of posets (which we denote by the same letter) $f\colon \calP(\Lambda)\longrightarrow \calP(\Lambda')$ by the formula $\Delta\mapsto f(\Delta)$. It follows that $f$ induces a map of spaces $f\colon |\calP(\Lambda)|\longrightarrow |\calP(\Lambda')|$ (this much would be true for any fusion $f$). It is also true that $f$ takes boundary to boundary, i.e., $f$ restricts to a map $f\colon \partial|\calP(\Lambda)|\longrightarrow \partial|\calP(\Lambda')|$ (this is only true if $f$ is a strict fusion). Thus, we have a covariant functor from the category $\calE_n$ to the category of pairs of spaces
$$\Lambda\mapsto (|\calP(\Lambda)|,\partial|\calP(\Lambda)|).$$
By passing to quotient spaces, we also obtain a functor from $\calE_n$ to pointed spaces
$$\Lambda\mapsto T_\Lambda.$$

Next we need to introduce another, contravariant, functor on $\calE_n$. It is well known that the poset of partitions of $s$ is a lattice, in the sense that any two partitions $\Lambda$ and $\Delta$ of $s$ have a coarsest common refinement, denoted $\Delta\vee\Lambda$ and a finest common coarsening, denoted $\Delta\wedge\Lambda$. Let $\Lambda$  be a partition of $s$. Let $\Delta$ be another partition of $s$. 
For a partition $\Lambda$, let $\comp(\Lambda)$ be the set of components of $\Lambda$. $\Lambda$ defines a natural partition of the set $\comp(\Delta)$: two elements of $\comp(\Delta)$ are equivalent if they belong to the same component of $\Lambda\wedge\Delta$. It is easy to see that the map of sets $s\twoheadrightarrow \comp(\Delta)$, associated with the partition $\Delta$, defines a natural fusion of $\Lambda$ into this partition of $\comp(\Delta)$. Now comes the crucial definition: we say that $\Delta$ is {\em good} relative to $\Lambda$ if the above fusion is a strict fusion. Otherwise, we say that $\Delta$ is {\em bad} relative to $\Lambda$.  


Let $M$ be a smooth manifold. When $\Lambda$ is a partition of $s$, we use the notation $M^\Lambda$ to mean $M^s$. This is the space of maps from $\Lambda$ to $M$. Suppose that $\Lambda$ is represented by the surjection $s\twoheadrightarrow c$. The surjection induces an inclusion map $M^c\hookrightarrow M^s$. We identify $M^c$ with its image in $M^s$. This is the diagonal subspace of $M^s$ associated with $\Lambda$. Now define the space $\Delta^\Lambda M\subset M^\Lambda$ as follows 
$$\Delta^\Lambda M=\underset{\underset{\mbox{relative to } \Lambda}{\Delta \mbox{ is bad }}}{\bigcup}  M^{\comp(\Delta)}.$$
In words, $\Delta^\Lambda M$ is the union of diagonals that are bad relative to $\Lambda$.

Let $f\colon \Lambda\longrightarrow \Lambda'$ be a morphism in $\calE_n$. Clearly, $f$ gives rise to a map $M^{\Lambda'}\longrightarrow M^\Lambda$. It also is true that this map restricts to a map $\Delta^{\Lambda'}M\longrightarrow \Delta^\Lambda M$. Thus we have a contravariant functor from $\calE_n$ to pairs of spaces
$$\Lambda\mapsto (M^\Lambda, \Delta^\Lambda M).$$

Given the category $\calE_n$ and a (contravariant) functor from $\calE_n$ to spaces, we can assemble this data by means of the following standard construction. Define 
$\calE_n\ltimes M^{\Lambda}$ to be the topological category where an object is a pair of the form
$(\Lambda, \alpha)$
where $\Lambda$ is an object of $\calE_n$, and $\alpha\colon\Lambda\longrightarrow M$ is a point of $M^\Lambda$. A morphism
$$(\Lambda, \alpha)\longrightarrow (\Lambda_1, \alpha_1)$$ 
consists of a morphism $\Lambda\longrightarrow \Lambda_1$ in $\calE_n$, that takes $\alpha_1$ to $\alpha$ via the functoriality of $M^\Lambda$. This is a topological category, in which the space of objects is topologized as the disjoint union of spaces of the form $M^\Lambda$, and the space of morphisms also has a topology (topological categories are reviewed in Section~\ref{S: CategoriesAndLimits}). Inside this category we have a subcategory that we will denote by
$\calE_n\ltimes \Delta^\Lambda M$
It is the full subcategory of $\calE_n\ltimes M^{\Lambda}$
consisting of objects $(\Lambda, \alpha)$ where $\alpha\in\Delta^\Lambda M$.

Next, we would like to define a contravariant functor  
$$\Psi\colon \calE_n\ltimes M^\Lambda\longrightarrow \Spectra.$$ Recall that we have a map $h\colon M\longrightarrow N$ at the background, coming from the chosen basepoint in $\Imm(M,N)$ or $\Emb(M,N)$. One can think of $h$ as an immersion or even an embedding if one wishes, but this does not matter for the construction we want to make now. Let $(\Lambda, \alpha)$ be an object of 
$\calE_n\ltimes M^{\Lambda}$. Let $\Cyl_\Lambda$ be the mapping cylinder of the surjection representing $\Lambda$. There is a natural inclusion $i\colon s\hookrightarrow \Cyl_{\Lambda}$ which induces a map $i^*\colon N^{\Cyl_\Lambda}\longrightarrow N^s$. Notice that $h\circ\alpha$ defines a point in $N^s$. 
Let $\Omega^\Lambda_{h\alpha} N$ be the space defined by means of a pullback square
$$\begin{array}{ccc}
\Omega^\Lambda_{h\alpha} N & \longrightarrow & N^{\Cyl_\Lambda}\\
\downarrow & & \,\mbox{}\downarrow i^* \\
* &\stackrel{h\circ\alpha}{\longrightarrow} & N^{s} \end{array}$$
Notice that if, for example, $N$ is connected then there is an equivalence 
$$\Omega^\Lambda_{h\alpha} N\simeq  (\Omega N)^{\excess(\Lambda)}$$ where $\Omega N$ is the ordinary loop space of $N$ (for some choice of basepoint). This is because the quotient $\Cyl_\Lambda/\supp(\Lambda)$ is homotopy equivalent to a wedge of $\excess(\Lambda)$ circles.


Finally, define the functor $\Psi$ on objects of $\calE_n\ltimes M^\Lambda$ by the formula
$$\Psi(\Lambda,\alpha) = \map_*\left(T_\Lambda, \Sigma^\infty \Omega^\Lambda_{h\alpha} N_+ \wedge S^{d\excess(\Lambda)}\right).$$
Since our source category is a topological category, this spectrum really should be thought of as a fiber in a ``homotopy bundle of spectra'' over $M^\Lambda$ (the notion of a homotopy bundle of spectra is reviewed in Section~\ref{S: Fiberwise}). The contravariant functoriality of $\Psi$ is determined by the fact that $T_\Lambda$ is a covariant functor on $\calE_n$, and by the observation that a morphism $f\colon \Lambda\longrightarrow\Lambda'$ in $\calE_n$, together with a choice of $\alpha'\in M^{\Lambda'}$ determines a map (which is, incidentally, a homotopy equivalence)
$$\Omega^{\Lambda'}_{h\alpha'} N_+ \wedge S^{d\excess(\Lambda')}\longrightarrow \Omega^\Lambda_{h\alpha'f} N_+ \wedge S^{d\excess(\Lambda)}.$$

Let $$\widetilde{\nat}_{\calE_n}\left((M^\Lambda,\Delta^\Lambda M); \map_*\left(T_\Lambda,\Sigma^\infty \Omega^\Lambda_{h\alpha} N_+\wedge S^{d\excess(\Lambda)}\right)\right)$$
be the fiber of the map 
$$\underset{\calE_n\ltimes M^\Lambda}{\lim} \Psi \longrightarrow \underset{\calE_n\ltimes \Delta^\Lambda M}{\lim} \Psi.$$
During the introduction, we will some times abbreviate this as $\widetilde{\nat}_{\calE_n}(-,-)$. The notation is meant to suggest a space (more precisely, a spectrum) of ``twisted natural transformations''.  $M^\Lambda$ is a contravariant functor on $\calE_n$. The construction $$\Lambda\mapsto \map_*\left(T_\Lambda,\Sigma^\infty \Omega^\Lambda_{h\alpha} N_+\wedge S^{d\excess(\Lambda)}\right)$$
also tries to be a contravariant from $\calE_n$ to $\Spectra$, but it also is a bundle over $M^\Lambda$. In this situation, one can define a space (or spectrum) of twisted natural transformations, analogous to the usual space of natural transformations, where mapping spaces are replaced with spaces of sections of bundles of spaces (or, as it is happens in our case, bundles of spectra). The relevant definitions are given in sections~\ref{S: Fiberwise}, \ref{S: CategoriesAndLimits}, and especially \ref{S: TwistedNat}. Moreover, the notation is supposed to remind us that we only are looking at natural transformations whose restriction to $\Delta^\Lambda M$ is trivial. 

We are finally ready to state our main result.
\begin{theorem}\label{T: MainTheorem}
Suppose that $M$ is tame, and $N$ is a tame stably parallelizable manifold. There is a natural weak equivalence between $\D_n\Sigma^\infty \Ebar(M,N)$ and 
$$\widetilde{\nat}_{\calE_n}\left((M^\Lambda,\Delta^\Lambda M); \map_*\left(T_\Lambda,\Sigma^\infty \Omega^\Lambda_{h\alpha} N_+\wedge S^{d\excess(\Lambda)}\right)\right).$$
\end{theorem}
\begin{remark}\label{R: Parallelizabily Unnecessary}
We believe that the theorem holds almost as stated without the assumption that $N$ is stably parallelizable, except that in place of the smash product with  $S^{d\excess(\Lambda)}$ one would have a Thom space of a certain evident bundle. However, our proof of the theorem relies on a description of the Taylor polynomial (as opposed to just the homogeneous layer) of the functor $\Sigma^\infty \Ebar(M,N)$, and we do not know how to describe the Taylor polynomials without an assumption of parallelizability.
\end{remark}
\begin{remark}
The spectrum of natural transformations that appears in the statement of the main theorem can be described as the spectrum of restricted sections of a certain homotopy bundle of spectra over the coend
$$\left(M^\Lambda,\Delta^\Lambda M\right)\otimes_{\Lambda\in\calE_n} \left(|\calP(\Lambda)|,\partial|\calP(\Lambda)|\right).$$
Note that this coend is a pair of spaces. ``Restricted'' means that we are only looking at sections that are trivial on the subspace part of the pair. To put it a little differently: we are thinking of the bundle of spectra, whose generic fiber is 
$\Sigma^\infty \Omega^\Lambda_{h\alpha} N_+\wedge S^{d\excess(\Lambda)}$,
as defining a kind of twisted cohomology theory on the category of (pairs of) spaces over 
$M^\Lambda\otimes_{\calE_n} |\calP(\Lambda)|$, and we are taking the relative cohomology of the pair given by the coend above.
\end{remark}
\begin{remark}\label{R: Alternative}
An alternative way to think of the space  $|\calP(\Lambda)|$ is as the space of rooted forests whose underlying partition is $\Lambda$. It follows that $M^\Lambda\otimes_{\calE_n} |\calP(\Lambda)|$ can be thought of as the space of rooted forests, whose leaves are marked by points of $M$, and whose underlying partition has excess $n$. In Figure~\ref{Fig: TwoForests}, we try to illustrate this space of forests in the case $n=2$, together with the fundamental homotopy bundle over it. 
\begin{figure}
\centering
\includegraphics[width=.7\linewidth]{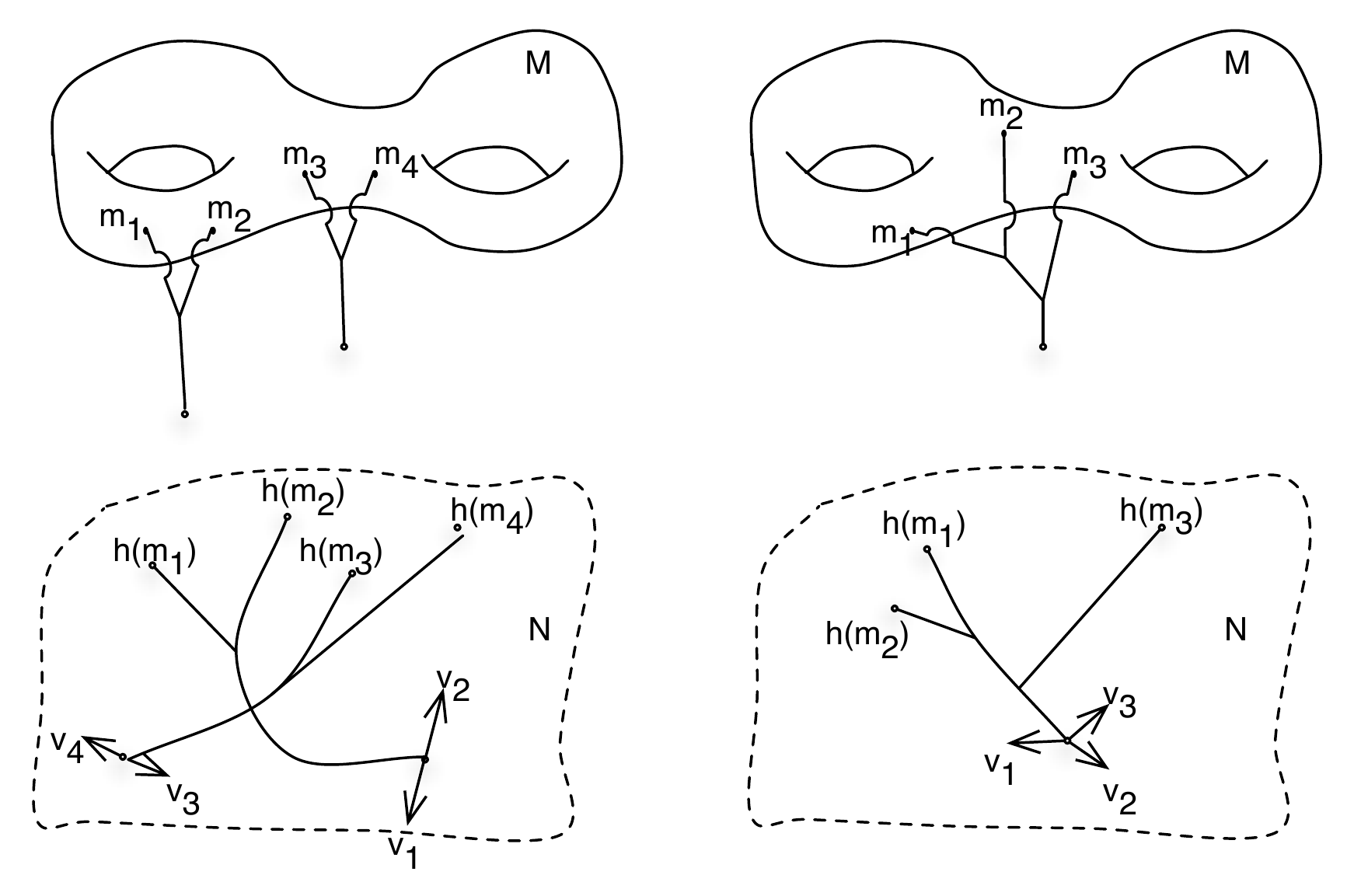}
\caption{The fundamental homotopy bundle in the case $n=2$: Above, we have two points in $M^\Lambda\otimes_{\calE_2} |\calP(\Lambda)|$ - the space of forests of excess $2$ with leaves marked by points in $M$, and underneath each one of them, a point in the corresponding fiber $\Omega^\Lambda_{h\alpha}N_+ \wedge S^{2d}$. The vectors $v_1, v_2,\ldots$ represent the part having to do with $S^{2d}$. The reason this part is there is that our space is really the Thom space of a pullback of the normal bundle of $N^{\comp(\Lambda)}$ in $N^{\supp(\Lambda)}$. In this picture, a normal vector to the diagonal $N$ in $N^i$ is represented by an $i$-tuple of tangent vectors that add up to zero.}
\label{Fig: TwoForests}
\end{figure}
This viewpoint will be developed further in the second part of the series~\cite{Arone}, where we will give an analogous description of the derivatives of $\Ebar(M, N)$ as the twisted cohomology of a certain space of connected {\em graphs} (as opposed to forests) with points marked by elements of $M$.
\end{remark}

It may not be immediately obvious that what we have described in Theorem~\ref{T: MainTheorem} is a homogeneous functor of $N$. The reason is that, first, $\Sigma^\infty \Omega^\Lambda_{h\alpha} N_+\wedge S^{d\excess(\Lambda)}$ is a homogeneous functor of degree $\excess(\Lambda)$, and the spectrum described in Theorem~\ref{T: MainTheorem} is a compact homotopy limit of such functors. To make this more explicit, we will observe (Remark~\ref{R: StructureOfE_n}) that the category $\calE_n$ admits a nice filtration by certain sub-categories $\calE_n^i$, where $1\le i\le n$ and $\calE_n^i$ is the category of partitions of excess $n$ and support of size $\le n+i$. Accordingly, the category $\calE_n\ltimes M^\Lambda$ can be filtered by sub-categories $\calE_n^i\ltimes M^\Lambda$. Obviously, there is a tower of restriction maps
$$\widetilde{\nat}_{\calE_n}(-,-)= \widetilde{\nat}_{\calE_n^n}(-,-) \longrightarrow \widetilde{\nat}_{\calE_n^{n-1}}(-,-) \longrightarrow \cdots \longrightarrow \widetilde{\nat}_{\calE_n^1}(-,-).$$
\begin{proposition} \label{P: EmbeddingOfOrthogonal}
The homotopy fiber of the map $\widetilde{\nat}_{\calE_n^{i}}(-,-) \longrightarrow \widetilde{\nat}_{\calE_n^{i-1}}(-,-) $
is equivalent to the product
$$\prod_{[\Lambda]} \Gamma\left((M^{n+i},\Delta^{n+i} M);\map_*\left(T_\Lambda,\Sigma^\infty \Omega^\Lambda_{h\alpha} N_+\wedge S^{d\excess(\Lambda)}\right)\right)^{\Sigma_{\Lambda}}.$$
Here the product is indexed by isomorphism types of objects of $\calE_n^i\setminus \calE_n^{i-1}$, i.e., isomorphism types of irreducible partitions of excess $n$ and support of size $n+i$. Thus $M^{n+i}=M^\Lambda$, and $\Delta^{n+i} M$ stands for the fat diagonal in $M^{n+i}$. $\Gamma(-,-)$ stands for the spectrum of restricted sections of the natural bundle of spectra over $M^\Lambda$ whose fiber at a point $\alpha$ is $\map_*\left(T_\Lambda,\Sigma^\infty \Omega^\Lambda_{h\alpha} N_+\wedge S^{d\excess(\Lambda)}\right)$ (``restricted'' means that we demand that the section restricts to the trivial section on $\Delta^{n+i}M$). $\Sigma_{\Lambda}$ is the group of automorphisms of $\Lambda$. It acts on the bundle of spectra, and the superscript $\Sigma_{\Lambda}$ indicates that we are taking sections that are invariant under the action.
\end{proposition}

\begin{example}
$$\widetilde{\nat}_{\calE_n^1}(-,-)\simeq \Gamma\left((M^{n+1},\Delta^{n+1}M) ;\map \left(T_{n+1},\Sigma^\infty \Omega^\Lambda_{h\alpha} N_+\wedge S^{dn}\right)\right)^{\Sigma_{n+1}}.$$

On the other extreme, the homotopy fiber of the map $\widetilde{\nat}_{\calE_n^n}(-,-)\longrightarrow \widetilde{\nat}_{\calE_n^{n-1}}(-,-)$ is equivalent to 
$$\Gamma\left((M^{2n},\Delta^{2n} M) ;\map_*\left(T_2^{\wedge n},\Sigma^\infty \Omega^\Lambda_{h\alpha} N_+\wedge S^{dn}\right)\right)^{\Sigma_2\wr\Sigma_n}.$$ In the case when $M$ is one-dimensional, this term corresponds to the ``chord diagrams'' familiar from knot theory. 
\end{example}

\begin{remark}\label{R: Consequences}
Since $\Sigma^\infty \Ebar(M, N)$ is, in addition to being a covariant functor of $N$, also a contravariant functor of $M$ (where $M$ is considered an object in a suitable category of manifolds over $N$), one can apply embedding calculus~\cite{WeissEmb} to this functor, and also to the functor  
$\D_n\Sigma^\infty\Ebar(M, N)$. Propositioin~\ref{P: EmbeddingOfOrthogonal} in fact describes the Taylor tower in the sense of embedding calculus of $\D_n\Sigma^\infty\Ebar(M, N)$. In particular, it shows that $\D_n\Sigma^\infty\Ebar(M, N)$ is, as a functor of $M$, a polynomial functor of degree $2n$, and it has layers in degrees between $n+1$ and $2n$.
\end{remark}

Notice that if $\Lambda$ is a partition of a set with $n+i$ elements, then $\Sigma_\Lambda$ is a subgroup of the symmetric group $\Sigma_{n+i}$. $\Sigma_{n+i}$ is acting freely on the complement of $\Delta^{n+i}M$ in $M^{n+i}$, and so the action of $\Sigma_\Lambda$ is free as well. Since $(M^{n+i},\Delta^{n+i}M)$ is equivariantly equivalent to a finite relative CW complex, it follows that   
\begin{multline*}
\Gamma\left(\left(M^{n+i},\Delta^{n+i} M\right);\map_*\left(T_\Lambda,\Sigma^\infty \Omega^\Lambda_{h\alpha} N_+\wedge S^{d\excess(\Lambda)}\right)\right)^{\Sigma_{\Lambda}}\simeq \\ \simeq\Gamma\left(\left(M^{n+i},\Delta^{n+i} M\right);\map_*\left(T_\Lambda,\Sigma^\infty \Omega^\Lambda_{h\alpha} N_+\wedge S^{d\excess(\Lambda)}\right)\right)_{\homotopy\Sigma_{\Lambda}}
\end{multline*}
by a standard Adams-isomorphism type argument. In this form, it is easy to see that the functor is homogeneous of degree $\excess(\Lambda)$. It follows that the functor $$\widetilde{\nat}_{\calE_n}\left((M^\Lambda,\Delta^\Lambda M); \map_*\left(T_\Lambda,\Sigma^\infty \Omega^\Lambda_{h\alpha} N_+\wedge S^{d\excess(\Lambda)}\right)\right)$$
is homogeneous of degree $\excess(\Lambda)$, because if a functor resolves into a finite tower of fibrations where each fiber is homogeneous of degree $n$, then the functor itself is homogeneous of degree $n$.

\begin{remark}
The description of the layers of the functor $\Sigma^\infty \Ebar(M,N)$ extend automatically to a description of the layers of any functor of the form $N\mapsto E\wedge \Ebar(M,N)$, where $E$ is a spectrum. One just replaces ``$\Sigma^\infty$'' with ``$E\wedge $'' everywhere in the formula. For example, if we take $E=\operatorname{HR}$, the Eilenberg-Mac Lane spectrum of our favorite coefficient ring, then the Taylor tower of the functor $\operatorname{HR}\wedge \Ebar(M,N)$ gives rise to a spectral sequence for $\HH_*( \Ebar(M,N); R)$ (the Taylor tower is known to converge if $2\dim(M) + 1<\dim(N)$). The formula for the layers in this case is 
$$\widetilde{\nat}_{\calE_n}\left((M^\Lambda,\Delta^\Lambda M); \map_*\left(T_\Lambda,\operatorname{HR}\wedge  \Omega^\Lambda_{h\alpha} N_+\wedge S^{d\excess(\Lambda)}\right)\right),$$
which in turn breaks up into pieces of the form
$$\Gamma\left((M^{n+i},\Delta^{n+i} M) ;\map_*\left(T_\Lambda, \operatorname{HR}\wedge\Omega^\Lambda_{h\alpha} N_+\wedge S^{d\excess(\Lambda)}\right)\right)^{\Sigma_{\Lambda}}.$$
This homotopy groups of this spectrum are given, roughly speaking, by the equivariant cohomology of the space $M^{n+i}/\Delta^{n+i} M\wedge T_\Lambda$ with certain local coefficients. The upshot of all this is that while this is not a trivial calculation, quite a lot is known about the equivariant cohomology of the spaces $M^{n+i}/\Delta^{n+i} M$ and  $T_\Lambda$ (references \cite{SchwartzGenus, AK} contain some information about the latter), and so it might be possible to do some interesting calculations with these formulas. We intend to come back to this in a future paper.
\end{remark}

If one wants to, one can rewrite the functor given in Theorem~\ref{T: MainTheorem} 
in the canonical form of a homogenous functor, represented
by a spectrum with an action of $\OO(n)$. To do this, proceed as follows. Let $\Lambda$ be an object of $\calE_n$. Thus $\Lambda$ is a partition of some set $s$. We call $s$ the support of $\Lambda$, and when we want to underscore the dependence of $s$ on $\Lambda$ we write $\supp(\Lambda)$ for the support of $\Lambda$. Recall that $\Cyl_\Lambda$ is the mapping cylinder of the surjection representing $\Lambda$ and that $\supp(\Lambda)\subset \Cyl_\Lambda$. For any $\Lambda$ in $\calE_n$, $\HH^1(\Cyl_\Lambda,\supp(\Lambda);\R)$ is isomorphic to $\R^{\excess(\Lambda)}\cong \R^n$. Let $\iso\left(\R^n,\HH^1(\Cyl_\Lambda,\supp(\Lambda);\R)\right)$ be the space of vector space isomorphisms. This is a space homeomorphic to $\GL_n(\R)$ with a canonical action of $\OO(n)$. Let us denote this space by $\GL(\Lambda)$. Thus the assignment $\Lambda\mapsto \GL(\Lambda)$ defines a contravariant functor from $\calE_n$ to spaces with an action of $\OO(n)$. Now define the semi-direct category in the evident way
$$\calE_n\ltimes \left(\GL(\Lambda)\times M^\Lambda\right).$$
This is a category with an action of $\OO(n)$. Standard manipulations with limits yield the following proposition as a corollary of Theorem~\ref{T: MainTheorem}.
\begin{proposition}\label{P: DescriptionOfDerivatives} 
There is a natural equivalence between $\D_n\Sigma^\infty \Ebar(M,N)$ and 
$$\widetilde{\nat}_{\calE_n}\left(\left(\GL(\Lambda)\times (M^\Lambda,\Delta^\Lambda M)\right); \map_*\left(T_\Lambda,\Sigma^\infty \Omega^\Lambda_{h\alpha} N_+\wedge S^{d\excess(\Lambda)}\right)\right)^{\OO(n)}.$$
Equivalently, the $n$-th derivative of the functor is the spectrum with an action of $\OO(n)$
$$\widetilde{\nat}_{\calE_n}\left(\left(\GL(\Lambda)\times (M^\Lambda,\Delta^\Lambda M)\right); \map_*\left(T_\Lambda,\Sigma^\infty \Omega^\Lambda_{h\alpha} N_+\wedge S^{d\excess(\Lambda)}\right)\right)\wedge S^{\ad_n}$$
where $\ad_n$ is the adjoint representation of $\OO(n)$.
\end{proposition} 
\begin{remark}
The twist by the adjoint sphere arises because of the way Adams' isomorphism works for compact Lie groups. For a spectrum $C_n$ with an action of $\OO(n)$, there is a natural zig-zag of maps 
$$(C_{n} \wedge S^{\ad_n})_{\homotopy\! \OO(n)}\longrightarrow C_n^{\homotopy\!\OO(n)}.$$
This map is a weak equivalence if $C_n$ is a finite spectrum with a free action of $\OO(n)$.
\end{remark}

Next, we make a remark about the homotopy invariance of the layers. It is easy to see that our description of the layers only depends on the stable homotopy types of $M$, $N$ and the basepoint map $h\colon M\to N$. More precisely, we have the following corollary.
\begin{corollary}\label{C: Homology Invariance}
Let $E$ be a spectrum. Suppose that $N$ and $N_1$ are stably parallelizable manifolds of the same dimension. Suppose that the embeddings $h\colon M\to N$ and $h_1\colon M_1\to N_1$ are related by a diagram of the form
$$
\begin{CD}
M @>>> N \\
@V{\cong_{E_*}} VV @VV{\cong_{E_*}}  V \\
X @>>> Y\\
@A{\cong_{E_*}}  AA @AA{\cong_{E_*}}  A \\
M_1 @>>> N_1
\end{CD}
$$
where $X$ and $Y$ are spaces and all the vertical maps induce an isomorphism in $E_*$-homology. Then the layers $\D_nE\wedge \Ebar(M,N)$ and $\D_n E\wedge \Ebar(M_1, N_1)$ are homotopy equivalent. 
\end{corollary}
\begin{remark}
We believe that a version of the corollary holds without the parallelizability assumption, but in this case one would need to demand that $N$ and $N_1$ are tangentially homotopy equivalent in a suitable sense.
\end{remark}
An important special case is when $N$ is a Euclidean space. In this case, our custom is to use $V$ instead of $N$, so we are looking at the functor $\Sigma^\infty\Ebar(M,V)$, where $V$ is a Euclidean space. In this case, $\Omega^\Lambda_{h\alpha} V$ is contractible, and $\Omega^\Lambda_{h\alpha} V_+\wedge S^{d\excess(\Lambda)}$ can be simplified to $S^{V\excess(\Lambda)}$. Moreover, it is not difficult to see that in this case, all the different spheres $S^{V\excess(\Lambda)}$ can be identified with a single sphere $S^{nV}$, and the various section spaces become equivariant mapping spaces. Theorem~\ref{T: MainTheorem} specializes to the following
\begin{corollary}\label{C: EuclideanCase}
$\D_n\Sigma^\infty\Ebar(M,V)$ is equivalent to
\begin{multline*}
\map_*\left(GL(\Lambda)_+\wedge M^{[\Lambda]} \otimes_{\calE_n} T_\Lambda,\Sigma^\infty S^{nV}\right)^{\OO(n)}\simeq \\ \simeq \map_*\left(GL(\Lambda)_+\wedge M^{[\Lambda]} \otimes_{\calE_n} T_\Lambda,\Sigma^\infty S^{\ad_n}\wedge S^{nV}\right)_{\homotopy\OO(n)}
\end{multline*}
where $M^{[\Lambda]}=M^\Lambda/\Delta^\Lambda M$, considered as a contravariant functor on $\calE_n$. The symbol $\otimes_{\calE_n}$ stands for the coend of a contravariant functor and a covariant functor.
\end{corollary}
And corollary~\ref{C: Homology Invariance} specializes to the following.
\begin{corollary}\label{C: HomotopyInvariance}
Suppose that $M_1$ and $M_2$ are related by a zig-zag of maps that induce equivalence in $E_*$-homology.
Then 
$$\D_nE\wedge \Ebar(M_1,V)\simeq \D_n E\wedge \Ebar(M_2,V)$$
for any choice of basepoint embeddings $M_1\hookrightarrow V$ and $M_2\hookrightarrow V$.
\end{corollary}
This corollary was used in~\cite{ALV}.

Now let us discuss the proof of Theorem~\ref{T: MainTheorem}. A good way to start understanding the functor $\Sigma^\infty \Ebar(M, N)$ and other similar functors is to first investigate the case when $M$ is the disjoint union of finitely many open balls (one could say that Embedding Calculus~\cite{WeissEmb} is about reducing the case of general $M$ to this case). So let us suppose that $M\cong D^m\times k$, where $D^m$ is the unit open ball in $\R^m$ and $k$ is a finite set. Let $\config(k, N)$ be the familiar configuration space of $k$-tuples of disjoint points in $N$. $\config(k, N)$ can be identified with $\Emb(k, N)$, where $k$ is viewed as a zero-dimensional manifold. Similarly, $N^k$ may be identified with $\Imm(k, N)$. Let $\overline{\config}(k, N)$ (which could also be called $\Ebar(k, N)$) be the homotopy fiber of the inclusion map $\config(k, N)\longrightarrow N^k$. There are maps $\Emb(D^m \times k, N)\longrightarrow \config(k, N)$ and $\Imm(D^m \times k, N)\longrightarrow N^k$ defined by evaluation at the center of each ball. There is a commutative square, where the horizontal maps are the evaluation maps, and the vertical maps are inclusions
$$\begin{CD}
\Emb(D^m \times k, N)@>>> \config(k, N) \\
@VVV @VVV \\
\Imm(D^m \times k, N) @>>> N^k \end{CD}.$$
It is well-known that this square is a homotopy pullback. It follows that the induced map on vertical homotopy fibers $\Ebar(M, N) \longrightarrow \overline{\config(k, N)}$ is a homotopy equivalence (these definitions depend on a choice of a basepoint $h\in \Imm(D^m\times k, N)$, which we suppress for now). It follows that there is a homotopy equivalence, natural in $N$, $\Sigma^\infty \Ebar(M, N)\longrightarrow \Sigma^\infty\overline{\config}(k, N)$. Therefore the two functors have homotopy equivalent Taylor towers.

To analyze the Taylor tower of $\Sigma^\infty\overline{\config}(k, N)$, let us consider the functor $\config(k, N)$ first. One can identify $\config(k, N)$ with $N^k \setminus \Delta^kN$, where $\Delta^kN\subset N^k$ is the fat diagonal. $\Delta^k N$ can be identified with the union (= colimit) 
of the spaces $N^{\comp(\Lambda)}$ where $\Lambda$ ranges over the category of non-discrete partitions of $k$. For the purposes of the introduction, let us denote this category by $\calP\sp{0}_k$. It follows that $\config(k, N)$ can be identified with the intersection (= limit) of the spaces $N^{k}\setminus N^{\comp(\Lambda)}$. Stabilizing, one obtains a natural map
$$\Sigma^\infty\config(k, N) \longrightarrow \underset{\Lambda\in \calP\sp{0}_k}{\holim}\Sigma^\infty N^k\setminus N^{\comp(\Lambda)}$$
and we prove that this map is a homotopy equivalence (Proposition~\ref{P: ModelForStableConfig}). Similarly, we obtain a homotopy equivalence
\begin{equation}\label{E: Model}
\Sigma^\infty\overline{\config}(k, N) \longrightarrow \underset{\Lambda\in \calP\sp{0}_k}{\holim}\Sigma^\infty \hofiber\left(N^k\setminus N^{\comp(\Lambda)}\longrightarrow N^k\right).\end{equation}
It turns out that $\Sigma^\infty \hofiber\left(N^k\setminus N^{\comp(\Lambda)}\longrightarrow N^k\right)$ is, as a functor of $N$, homogeneous of degree $\excess(\Lambda)$. In fact, there is a natural equivalence $$\Sigma^\infty \hofiber\left(N^k\setminus N^{\comp(\Lambda)}\longrightarrow N^k\right)\simeq \Sigma^\infty \Omega^\Lambda_{h} N_+ \wedge S^{d\excess(\Lambda)},$$ which should explain why the latter functor serves as a basic building block in our main theorem. Thus, Equation~\eqref{E: Model} essentially provides a model for the Taylor tower of $\Sigma^\infty\Ebar(M, N)$ when $M$ is a finite set. It is not hard to conclude from \eqref{E: Model} that Theorem~\ref{T: MainTheorem} is correct in this case. Moreover, $\Ebar(M, N)$ and the model of $\D_n\Sigma^\infty\Ebar(M, N)$ given by the theorem are both homotopy invariant under replacing $M$ with $M\times W$, where $W$ is a Euclidean space (or equivalently an open ball). It follows that the theorem is correct when $M=D^m\times k$. This constitutes a significant advance towards proving the theorem for a general $M$. In fact, it would almost amount to a proof, except in order to make the required leap it is necessary to have a model for $\D_n\Sigma^\infty\Ebar(D^m\times k, N)$ that is (contravariantly) functorial with respect to embeddings in the variable $D^m\times k$. However, the model given by \eqref{E: Model} really depends on identifying $M$ with a finite set, and is only functorial with respect to inclusions of finite sets. We are  missing embeddings between manifolds of the form $D^m \times k$ that are not injective on $\pi_0$. Thus, we need a model for $\Ebar(D^m\times k, N)$ that has more functoriality in the first variable than $\config(k, N)$.

At this point we pause to make some simplifying assumptions on $M$ and $N$. An easy calculus argument shows that it is enough to prove the main theorem in the case when $N$ is a nice submanifold of a Euclidean space, and the background map $h\colon M\longrightarrow N$ is a codimension zero embedding. We say that a submanifold $N$ of a Euclidean space $W$ is ``nice'', if $N$ is the interior of a closed codimension zero submanifold $\overline{N}\subset W$ with a boundary, and there exists an $\epsilon>0$ such that subset of $\overline N$ consisting of points that have distance $\le \epsilon$ from the boundary $\partial\overline N$ deformation retracts onto $\partial \overline N$. For example, if $N$ is a tubular neighborhood of a compact manifold without boundary then $N$ is nice. Note also that if $N$ is a nice submanifold of $W$, and $W_1$ is another Euclidean space, then $N\times W_1$ is a nice submanifold of $W\times W_1$. We can assume that $N$ is nice, because if $N$ is a tame stably parallelizable manifold then there exist Euclidean spaces $W_0, W$ such that $N\times W_0$ embeds as a nice open submanifold  of $W$. But for any functor $F$, the derivatives (and therefore the layers) of the functor $F_{W_0}(V):= F(W_0\times V)$ determine the derivatives (and therefore the layers) of the functor $F(V)$. Indeed, if we denote the former functor by $F_{W_0}$, then $\partial_n F\simeq \Omega^{nW_0}\partial _n F_{W_0}$. The reason we may assume that $h$ is a codimension zero embedding is that, first, we can assume that $h$ is an embedding after crossing $N$ with a Euclidean space, and second, the homotopy type of $\Ebar(M,N)$ does not change if $M$ is replaced with its tubular neighborhood in $N$. 

So, let us assume that $N$ is a nice submanifold of a Euclidean space $W$, and $M$ is homeomorphic to a disjoint union of open balls in $W$. When $N$ is an open subset of $W$ and $U$ is an open ball in $W$, let us define the space of standard embeddings from $U$ to $W$, denoted $\sEmb(U,W)$, to be the space of embeddings that differ from the inclusion by a translation (thus $\sEmb(U, W)\cong W$). Define $\sEmb(U, N)$ to be the subset of $\sEmb(U, W)$ consisting of those standard embeddings whose image lies in $N$. Finally, for $M$ a disjoint union of open balls, define $\sEmb(M, N)$ to be the space of embeddings that are standard on each connected component of $M$. Suppose that $M$ has $k$ connected components. Let $\config(k,N)$ be the configuration space of $k$-tuples of disjoint points in $N$. There is an evident map $\sEmb(M,N)\to \config(k,N)$ given by evaluation at the center of each ball. It is well-known, and not hard to show, that if the balls making up $M$ are small enough then this evaluation map is a homotopy equivalence (here one has to use the assumption that $N$ is nicely embedded). We will always assume that this condition is satisfied. Similarly, we may define $\sImm(M,N)$ to be the space of all maps from $M$ to $N$ that are standard on each ball (such maps are automatically immersions, which justifies the notation). Clearly, when the connected components of $M$ are small enough, evaluation at centers induces an equivalence $\sImm(M,N) \longrightarrow N^k$. These maps factor through $\Emb(M, N)$ and $\Imm(M, N)$, so we have a commutative diagram
$$
\begin{CD}
\sEmb(M,N) @>>> \Emb(M,N) @>>>\config(k, N) \\
@VVV @VVV @VVV \\
\sImm(M,N) @>>> \Imm(M,N) @>>> N^k
\end{CD}
$$
where the composed horizontal maps are weak equivalences (assuming, automatically, that the components of $M$ are small enough). Let $\sEbar(M, N)$ be the homotopy fiber of the map $\sEmb(M, N)\longrightarrow \sImm(M, N)$. We have an equivalence $\sEbar(M,N)\to \Ebar(M,N)$. Next, we construct an explicit model of the Taylor polynomials of the functor $\Sigma^\infty\sEbar(M,N)$. This model is a direct generalization of the model we gave previously for $\Sigma^\infty{\overline \config}(k, N)$. It is based on the category of partitions of $M$. But now it is functorial with respect to standard embeddings in the variable $M$. From our model of the Taylor polynomials of $\Sigma^\infty\sEbar(M,N)$, we derive a model for $\D_n \Sigma^\infty\sEbar(M,N)$. En route, we observe that our model for $\D_n \Sigma^\infty\sEbar(M,N)$ is, as a functor of $M$, functorial not only with respect to standard embeddings in the variable $M$, but with respect to all maps in this variable. This is the basic reason that $\D_n \Sigma^\infty\sEbar(M,N)$ is essentially a (contravariant) homotopy functor of $M$. Our model of $\D_n \Sigma^\infty\sEbar(M,N)$ is equivalent, functorially in $M$, to the model give by Theorem~\ref{T: MainTheorem} in the case when $M$ is the disjoint union of open balls. This is enough to conclude the main theorem.

{\bf Section by section outline}: In Section~\ref{S: Preliminaries} we deal with preliminaries. After specifying the exact meaning of ``space'', ``spectrum'' and ``manifold'', we review some notions of fiberwise homotopy theory, both stable and unstable. In particular, we recall the definition of a homotopy bundle of spectra, and how it represents a generalized cohomology with local coefficients. In Section~\ref{S: CategoriesAndLimits} we review some definitions having to do with topological categories (where the space of objects itself may have a topology) and the notion of a functor from a topological category to the category of spaces. There certainly is nothing new in this section, but we could not find a good reference for the definition of a functor in the general setting that we needed, so we wrote one down. We also include a discussion of homotopy limits of functors on a topological category, and spaces of natural transformations. At this point, we introduce the notion of spaces of twisted natural transformations, and study some of their properties. 

In Section~\ref{S: Partitions} we discuss partitions and various kinds of morphisms between them, as well as other notions about partitions that we need. We introduce the ``space of non-locally constant maps'' from a partition $\Lambda$ to a manifold $N$. This space can be denoted by $N^\Lambda\setminus N^{\comp(\Lambda)}$, and it serves as a basic building block for much that is done in this paper.  
In Section~\ref{S: Calculus} we begin to study the orthogonal tower in some baby cases related to $\Sigma^\infty\Ebar(M,N)$. In particular, we write down a model for $\D_n\Sigma^\infty\overline{\config}(k,N)$ (Proposition~\ref{P: ModelForStableConfig}). We also study the functor $\Sigma^\infty\hofiber(N^\Lambda\setminus N^{\comp(\Lambda)}\longrightarrow N^\Lambda)$, and show that it is equivalent to $\Omega\Sigma^\infty(\Omega^\Lambda_\alpha N)_+\wedge S^{d\excess(\Lambda)}$ (Lemma~\ref{L: BasicHomogeneousFunctor}). In Section~\ref{S: StandardEmbeddings} we introduce spaces of standard embeddings, denoted $\sEbar(M, N)$, and show how our analysis of $\config(k, N)$ extends to $\Ebar(M, N)$. In Section~\ref{S: DerivativesOfStandard} we write down a functorial model for $\D_n\Sigma^\infty\sEbar(M,N)$ (where $M$ and $N$ are manifolds for which $\sEbar(M, N)$ is defined). 
In the last section we complete the proof of Theorem~\ref{T: MainTheorem} and Proposition~\ref{P: EmbeddingOfOrthogonal}.


\section{Preliminaries}\label{S: Preliminaries}
\subsection{Spaces and spectra}
Throughout this paper, {\em space}, means a compactly generated, weak Hausdorff topological space. A {\em pointed space} is a space with a chosen basepoint. We let $\Top$ and $\Top_*$ denote the categories of spaces and pointed spaces respectively. It is well-known that $\Top$ and $\Top_*$ are closed symmetric monoidal categories. 
For two (pointed) spaces $X$ and $Y$, we let $\map(X,Y)$ ($\map_*(X,Y)$) denote the 
(pointed) space of (pointed) maps from $X$ to $Y$.

We take a very old fashioned approach to spectra. For us, a spectrum is a sequence of pointed spaces $\{E_n\}_{n=0}^\infty$, together with maps $S\sp{1}\wedge E_n \longrightarrow E_{n+1}$. An important role in this paper will be played by suspension spectra, including suspension spectra of unpointed spaces. For a pointed space $X$, the suspension spectrum of $X$, denoted by $\Sigma^\infty X$, is the spectrum $\{S^n\wedge X\}_{n=0}^\infty$. This construction can be extended to unpointed spaces as follows. For any space $X$ let $\tilde\Sigma^\infty X$ be the homotopy fiber (defined level wise) of the map $\Sigma^\infty X_+ \longrightarrow \Sigma^\infty S^0$ that is induced by the map $X\longrightarrow *$ ($X_+$ denotes the space $X$ with a disjoint basepoint added). There is a natural map $\tilde\Sigma^\infty X \longrightarrow \Sigma^\infty X$ that is a homotopy equivalence for any choice of a non-degenerate basepoint in $X$. Therefore, the construction $\tilde\Sigma^\infty X$ is a way to construct a ``suspension spectrum'' for an unpointed space. Subsequently, we will just use the notation $\Sigma^\infty X$ to mean the usual suspension spectrum when $X$ is pointed, and to mean $\tilde\Sigma^\infty X$ when $X$ is unpointed, or when we do not want to commit to a specific choice of basepoint.

\subsection{Fiberwise homotopy theory} \label{S: Fiberwise} We will need some notions from the theory of  fiberwise spaces and spectra. 

A fiberwise space over $B$ is the same thing as a map of spaces $f\colon E\longrightarrow B$.
One thinks of $f$ as defining a family of spaces, parametrized by points of $B$. We will often use the letter $F$ to denote a generic fiber in a fiberwise space, and use $F_b$ to denote $f^{-1}(b)$, where $b\in B$. 
A {\em fiberwise pointed space} over $B$ (often called an {\em ex-space} in the literature) consists of a map $f\colon E\longrightarrow B$, together with a section $s\colon B\longrightarrow E$. The section $s$ endows each fiber of $f$ with a basepoint. 


Let $E$ be a fiberwise pointed space over $B$, and let $K$ be a pointed space. The fiberwise wedge sum of $K$ and $E$, denoted by $K\vee_B E$ is defined to be the pushout of the diagram
$$K\times B \longleftarrow B \stackrel{s}{\longrightarrow} E.$$ $K\vee_B E$ is a pointed space over $B$, whose fiber at every point $b\in B$ is homeomorphic to $K\vee F_b$. The fiberwise product of $K$ and $E$ is, simply, the space $K\times E$, which is obviously a space over and under $B$. Finally, the fiberwise smash product of $K$ and $E$, denoted $K\wedge_B E$, is defined to be the pushout of the diagram 
$$B \longleftarrow K\vee_B E \longrightarrow K\times E
$$
$K\wedge_B E$ is a fiberwise pointed space over $B$, whose fiber at $b$ is $K\wedge F_b$.

Similarly, $\map_B(K,E)$ is defined to be the fiberwise pointed space over $B$, whose fiber at $b$ is $\map_*(K,F_b)$. More precisely, it is defined to be the subspace of the full mapping space $\map(K,E)$, consisting of maps $\alpha\colon K\longrightarrow E$ such that the image of $\alpha$ lies entirely inside one fiber $F_b$, for some $b\in B$, and such that the map $\alpha: K\longrightarrow F_b$ is a pointed map, where the basepoint of $F_b$ is $s(b)$.

With these constructions, the category pointed spaces over $B$ is tensored and cotensored over $\Top_*$. There is an adjunction, where $F(-,-)$ denotes, just for the moment, the space of maps in  the category of pointed spaces over $B$.
$$F(K\wedge_B E_1,E_2)\cong F(E_1, \map_B(K, E_2))$$
(see \cite{MS}, Proposition 1.2.8). 

\begin{remark}\label{R: Homeomorphism}
It is not difficult to see that there is a homeomorphism
$$\Gamma(\map_B(K, E))\cong \map_*(K,\Gamma(E)),$$
Where $\Gamma(E)$ is the space of sections of $E$, considered as a pointed space. 
\end{remark}

Let $I=[0,1]$. Having the construction of fiberwise product with $I$ (and fiberwise smash product with $I_+$) enables one to define the notions of a (pointed) fiberwise homotopy between fiberwise maps, and a fiber homotopy equivalence between fiberwise spaces. 
Note that the condition that a fiberwise map $f\colon E_1\longrightarrow E_2$ is a fiber homotopy equivalence is in general much stronger than the condition that $f$ induces a homotopy equivalence on all fibers. Note also that a fiber homotopy equivalence induces a homotopy equivalence of spaces of sections.

Another important example of fiberwise smash product and fiberwise mapping space is provided by the fiberwise suspension and the fiberwise loop space, denoted $S^1\wedge_B E$ and $\Omega_B E$ respectively. These constructions enable us to define the notion of a fiberwise spectrum.

\begin{definition}
A fiberwise spectrum over $B$ is a sequence $\{E_n\}_{n=0}^\infty$ of fiberwise pointed spaces over $B$, together with fiberwise pointed maps 
$$S^1\wedge_B E_n \longrightarrow E_{n+1}$$
(or, equivalently
$$ E_n \longrightarrow \Omega_B E_{n+1})$$
for each $n=0,1,2,\ldots$.
\end{definition}
An example of a fiberwise spectrum is the fiberwise suspension spectrum $\Sigma^\infty_B E$ of a fiberwise pointed space $E$. The fiberwise spectrum $\Sigma^\infty_B E$ is defined by the sequence of pointed spaces over $B$, $\{S^n\wedge_B E\}_{n=0}^\infty$. The structure maps are constructed using the evident isomorphisms
$$S^1\wedge_B \left(S^n\wedge_B E\right) \stackrel{\cong}{\longrightarrow} \left(S^1\wedge S^n\right)\wedge_B E\stackrel{\cong}{\longrightarrow} S^{n+1}\wedge_B E.$$

An important construction associated with fiberwise spaces and spectra is the space of sections. For a fiberwise (pointed) space $f\colon E\longrightarrow B$, let $\Gamma(f)$, or $\Gamma(E)$, be the (pointed) space of sections of $f$, topologized as a subspace of the space of maps $\map(B,E)$. Similarly, for a fiberwise spectrum $E$ over $B$, we define the ``spectrum of sections'' as follows.
\begin{definition}\label{D: SectionsSpectrum}
Let $E=\{E_n\}_{n=0}^\infty$ be a fiberwise spectrum over $B$. Define the spectrum $\Gamma(E)=\{\Gamma(E)_n\}_{n=0}^\infty$ by
$$\Gamma(E)_n:= \Gamma(E_n).$$
The structure map $\Gamma(E)_n\longrightarrow \Omega\Gamma(E)_{n+1}$ is defined as the composite
$$\Gamma(E_n)\longrightarrow \Gamma(\Omega_B E_{n+1}) \stackrel{\cong}{\longrightarrow}\Omega \Gamma(E_{n+1}).$$
Here the last map is the homeomorphism of Remark~\ref{R: Homeomorphism}.
\end{definition}

We think of $\Gamma(E)$ as defining a kind of twisted generalized cohomology of the space $B$, or more generally, on the category of spaces over $B$. However, one usually needs to make some assumptions on $E$ to ensure that $\Gamma(E)$ really has the properties that one expects of cohomology theory, such as homotopy invariance and excision. Perhaps the most obvious condition on $E$ that ensures that $\Gamma(E)$ has the correct behavior is that the map $E\longrightarrow B$ is a fibration, or fiber homotopy equivalent to a fibration. In the pointed case, one needs to require, in addition, that $E$ is {\em well-pointed} in the sense that the section $s\colon B\hookrightarrow E$ is a fiberwise cofibration (see \cite{CJ}). In this case, $\Gamma(E)$ has the expected (pointed) homotopy invariance and excision properties, and so one has some control over the homotopy type of the space of sections. However, the class of fibrations is somewhat inconvenient for the purposes of fiberwise stable homotopy theory, because the fiberwise suspension of a fibration may not be a fibration. On the other hand, fiberwise suspension tends to preserve local triviality conditions. We will now follow in the footsteps of~\cite{CJ}, in that we are going to identify a convenient class of fiberwise pointed spaces that on one hand is closed under most reasonable constructions that one can perform on fiberwise spaces (including fiberwise suspensions) and on the other hand has good homotopic properties. This is the class of  ``pointed homotopy bundles''. 

\begin{definition}
Let $p\colon E\longrightarrow B$ be a fiberwise space over $B$. We say that $E$ is a {\em homotopy bundle} over $B$, if every point $x\in B$ has an open neighborhood $U$, such that the space $p^{-1}(U)$ over $U$ is fiber homotopy equivalent to a product bundle $F\times U\longrightarrow U$.

Similarly, if $E$ is a fiberwise pointed space over $B$, then we say that $E$ is a {\em pointed homotopy bundle} if 
\begin{itemize}
\item
$E$ is {\em homotopy well pointed} in the sense that the section of $E$ is a fiberwise homotopy cofibration (\cite{CJ}, Definition II.1.15) 
\item
$E$ is a homotopy bundle and all the trivializations $E_U\simeq F\times U$ in the definition of a homotopy bundle can be chosen to be {\em pointed} fiber homotopy equivalences.
\end{itemize}
\end{definition}

\begin{example}
If $E$ is a fiberwise space over $B$, one can construct out of it a fiberwise pointed space by adjoining a fiberwise disjoint basepoint. The new total space $E\coprod B$. This space is always homotopy well-pointed. It is a pointed homotopy bundle over $B$ if $E$ is a homotopy bundle over $B$.
 \end{example}

As is pointed out in \cite{CJ}, for the theory of homotopy fiber bundles to work properly, it is desirable
to assume that the base space $B$ is an ENR. So from now on, we always assume that $B$ is an ENR. In the examples that we will consider, $B$ will actually be a manifold, so certainly an ENR. In some of the results of \cite{CJ} that we will quote, it is further assumed that $B$ is compact. In our examples, $B$ is a manifold $M$ that may not be compact. However, we always assume that $M$ is the interior of a compact manifold with boundary $\overline M$. We can always construct our bundles over $\overline M$ and then restrict to $M$. All our constructions will be homotopy invariant with respect to the inclusion $M\hookrightarrow \overline M$.

As mentioned above, the advantage of working with (pointed) homotopy bundles is twofold. The first one is that they have good closure properties. For example, it is easy to see that they are closed under fiberwise suspension and under pullbacks over a map of base spaces $B'\longrightarrow B$. Furthermore, a key fact here is that any fiberwise (pointed) map between homotopy fiber bundles is locally (pointed) fiber homotopic to a product map (\cite{CJ}, Propositions II.1.2 and II.1.25). It follows that (pointed) homotopy fiber bundles are closed under such constructions as fiberwise homotopy pushouts and fiberwise homotopy pullbacks ([op. cit.], Lemmas II.2.3 and II.2.7).

The second attractive property of homotopy bundles is that they are close enough to honest fibrations that spaces of sections have good homotopic properties. More precisely, we have the following proposition, which summarizes the bottom half of \cite{CJ}, page 153, starting with Proposition 1.27.
\begin{proposition}\label{P: EssentiallyFibration}
Let $E$ be a (pointed) homotopy bundle over $B$. There exists a homotopy well-pointed space $E^f$ over $B$, which can be constructed functorially in $E$, such that the map $E^f\longrightarrow B$ is a fibration, and $E^f$ is naturally (pointed) fiber homotopy equivalent to $E$.
\end{proposition}
It follows from the proposition that if $E$ is a (pointed) homotopy bundle, then $\Gamma(E)\simeq \Gamma(E^f)$.

Before we discuss the homotopy properties of section spaces of homotopy bundles in more detail, let us expand the discussion to include ``homotopy bundles of spectra''.
\begin{definition}
Let $E=\{E_n\}_{n=0}^\infty$ be a fiberwise spectrum over $B$. We say that $E$ is a {\em homotopy bundle of spectra} if each $E_n$ is a pointed homotopy fiber bundle over $B$.
\end{definition}
For example, the fiberwise suspension spectrum of a pointed homotopy bundle over $B$ is a homotopy bundle of spectra. As in the non-parametrized case, we will sometimes want to work with {\it a-priori} unpointed fiberwise spaces.
\begin{definition}\label{D: FiberwiseSuspension}
Let $E$ be an unpointed fiberwise space. Let $E \coprod B$ the fiberwise pointed space constructed out of $E$ in the obvious way. Let $\tilde\Sigma^{^\infty}_B E$ be the fiberwise spectrum obtained by taking the fiberwise homotopy fiber of the map 
$$\Sigma^{^\infty}_B (E \coprod B) \longrightarrow \Sigma^{^\infty}_B (B\coprod B).$$
\end{definition}

If $E$ is a homotopy bundle of spaces, then $\Sigma^{^\infty}_B E$ is a homotopy bundle of spectra. For any choice of a well-pointed section $B\hookrightarrow E$, there is a fiberwise homotopy equivalence $$\tilde\Sigma^{^\infty}_B E\stackrel{\simeq}{\longrightarrow} \Sigma^{^\infty}_B E$$
(that this map is a fiberwise homotopy equivalence follows from Lemma~\ref{L:  FiberInvariant} below). 
In particular, this shows that the construction $\Sigma^{^\infty}_B E$ is, to a large extent, independent of the section. As in the case of ordinary suspension spectra, we will subsequently use $\Sigma^{^\infty}_B E$ in lieu of $\tilde\Sigma^{^\infty}_B E$ when $E$ is unpointed, or when we do not want to commit ourselves to a particular section.

Let $E$ be a fiberwise space, or spectrum over $B$. Recall that we often use $F$ to denote a generic fiber of  $E$, and $F_b$ to denote the fiber over $b$. We will some times denote the space (or spectrum) of sections of $E$ by $\widetilde\map(B,F)$ (or $\widetilde\map(B, F_b)$, when we want to emphasize the dependence of $F$ on $B$). The notation is meant to suggest that the space of sections of $E$ is a kind of a twisted space of maps from $B$ to $F$. We will also need a relative version. If $B_0\subset B$ is a sub-ENR, we use $\widetilde\map((B, B_0), F)$ to denote the space (or spectrum) of sections that agree with the basepoint section on $B_0$ (so this is only defined for fiberwise pointed spaces and fiberwise spectra). 

Let us now list the relevant homotopy properties of $\widetilde\map\left(B, F\right)$, when $E\to B$ is a (pointed) homotopy bundle.
First, there is homotopy invariance, in $B$ and in $F$. The proof of the following lemma is contained in the proof of \cite{CJ}, Proposition II.2.8. It consists of a straightforward application of Proposition~\ref{P: EssentiallyFibration} and the homotopy lifting property.
\begin{lemma}
Let $B, B'$ be ENRs. Let $E$ be a homotopy bundle of either pointed spaces or spectra over $B$, with fiber homotopy equivalent to $F$. Let $\alpha_0, \alpha_1 \colon B'\longrightarrow B$ be homotopic maps. Then the pullbacks $\alpha_0^*(E)$ and $\alpha_1^*(E)$ are fiber homotopy equivalent, and the induced maps of section spaces
$$\alpha_0^*, \alpha_1^* \colon \widetilde\map\left(B, F\right)\longrightarrow  \widetilde\map\left(B', F\right)$$ are homotopic via this equivalence.
\end{lemma}
This lemma has the following corollary
\begin{corollary} \label{C: BaseInvariant}
Let $B, B', B$ be pairs of ENRs. Let $h\colon B'\longrightarrow B$ be a homotopy equivalence of pairs. Let $E$ be, as usual, a homotopy bundle over $B$ with fiber homotopy equivalent to $F$. Then the pullback $h^*(E)$ is a homotopy bundle over $B'$ with fiber homotopy equivalent to $F$. The following induced map is a homotopy equivalence
$$ \widetilde\map\left(B, F\right) \longrightarrow  \widetilde\map\left(B', F \right).$$
\end{corollary}

The homotopy invariance in the $F$ variable is given in the following lemma. It is an immediate consequence of Dold's theorem (\cite{CJ}, Theorems II.1.12 and II.1.29).
\begin{lemma}\label{L: FiberInvariant}
Let $f\colon E \longrightarrow E'$ be a fiberwise map between (pointed) homotopy bundles over $B$, such that $f$ induces a (pointed) homotopy equivalence on each fiber. Then $f$ induces a (pointed) fiber homotopy equivalence. 
\end{lemma}

Next, there is the Meyer-Vietoris property. The following lemma is a special case of [op. cit.], Propositions II.2.11 and II.2.14.
\begin{lemma}\label{L: Excision}
Let $E$ be a homotopy bundle of pointed spaces or spectra over an ENR $B$.
Let $B_1, B_2$ be closed sub-ENRs of $B$, such that $B=B_1\cup B_2$ and $B_0=B_1\cap B_2$ is an ENR. Then there is 
\begin{enumerate}
\item A homotopy fibration sequence, where the second map is the restriciton map
$$\widetilde\map\left((B,B_0), F\right)\longrightarrow \widetilde\map(B,F) \longrightarrow  \widetilde\map(B_0,F).$$
\item A homotopy pullback square, where the maps are restriction maps.
$$\begin{array}{ccc}
\widetilde\map(B,F) &\longrightarrow & \widetilde\map(B_1,F) \\
\downarrow & & \downarrow \\
\widetilde\map(B_2,F) &\longrightarrow & \widetilde\map(B_0,F) \end{array}$$
\end{enumerate}
\end{lemma}

\subsubsection{Stratified homotopy fibrations.} We will need to consider somewhat more general fiberwise objects than homotopy bundles. Suppose that $E_1$ is a homotopy bundle over an ENR $B$. Let $B_0\subset B$ be a sub-ENR. Let $E_0$ be the restriction of $E_1$ to $B_0$. Suppose that $E_0'\subset E_0$ is a homotopy bundle over $B_0$ that is fiberwise a subspace of $E_0$. Then we may define a pointed space $E$ over $B$ to be the subspace of $E_1$ consisting of all points over $B\setminus B_0$, together with all points of $E_0'$ ($E$ is given the subspace topology of $E_1$). There is a square of spaces of sections, which is both a strict pullback and a homotopy pullback. 
$$\begin{CD}
\Gamma(E) @>>> \Gamma(E_1)\\
@VVV @VVV \\
\Gamma(E_0') @>>> \Gamma(E_0) \end{CD}.$$
The following lemma is an easy consequence
\begin{lemma}\label{L: Deformation}
If in the above construction the inclusion $B_0\hookrightarrow B$ is a homotopy equivalence, then the restriction map $\Gamma(E)\longrightarrow \Gamma(E_0')$ is a homotopy equivalence.
\end{lemma}
A special case of the above construction is when $E_1$ is a pointed homotopy bundle over $B$ with fiber $F$, and $E_0'=B_0$. In this case, $\Gamma(E)=\widetilde\map\left((B, B_0), F\right)$.


\section{Topological categories and homotopy limits.} \label{S: CategoriesAndLimits}
\subsection{The definition.} The data for defining a (small) topological category consists of a pair of topological spaces $(\Ob,\Mor)$ (the space of objects and the space of morphisms respectively) together with structure maps (the identity map, the source and target maps, and the composition map)

$$i:\Ob\longrightarrow \Mor$$
$$s,t:\Mor\longrightarrow \Ob$$
$$c:\Mor\times_{^t\Ob^s}\Mor\longrightarrow \Mor$$

satisfying the evident identities. Here, by $\Mor\times_{^t\Ob^s}\Mor$ we mean the pullback of the diagram
$$\Mor\stackrel{t}{\longrightarrow}\Ob \stackrel{s}{\leftarrow} \Mor.$$
Given a topological category $\calC=(\Ob,\Mor)$, we define the nerve of $\calC$ to be the simplicial space $\Nerve\calC_\bullet$, where $\Nerve\calC_0=\Ob$, $\Nerve\calC_1=\Mor$, and for $i\ge 2$, $\Nerve\calC_i:=\Mor\times_{^t\Ob^s} \Nerve\calC_{i-1}= \Nerve\calC_{i-1}\times_{^t\Ob^s}\Mor$ is defined to be the pullback of the diagram
$$\Mor\stackrel{t}{\longrightarrow} \Ob \stackrel{s}{\leftarrow} \Nerve\calC_{i-1}$$
where the right map is the source map on the ``leftmost'' $\Mor$ factor of $\Nerve\calC_{i-1}$.
 The face and degeneracy maps in the simplicial nerve are defined in the familiar way, using the structure maps in the definition of topological category. As usual, the geometric realization of $\Nerve\calC_\bullet$ is called the classifying space of $\calC$, and we will denote it by $|\calC|$.

We will need to consider two types of {\em functors} on topological categories: functors between small topological categories, and functors from a small topological category to the category of spaces (or pointed spaces, or spectra). Rather than developing a unifying framework that encompasses all types of functors that we need, we proceed in a somewhat ad-hoc manner, and define the two types of functors separately.

A functor from one small category to another $\alpha:(\Ob_1,\Mor_1)\longrightarrow(\Ob_2,\Mor_2)$ consists, simply, of a pair of maps $\alpha_\Ob:\Ob_1\longrightarrow\Ob_2$ and $\alpha_\Mor:\Mor_1\longrightarrow\Mor_2$ that commute with all the structure maps. We say that $\alpha$ is a homotopy equivalence, if it induces a degree-wise equivalence of simplicial nerves.

Next, we need to consider the notion of a functor from a small topological category to the category of spaces, or pointed spaces, or spectra. Since it is important to us to discuss all three variants, and since we would like to avoid, as much as possible, repeating everything three times, let us, for the duration of this section, use the word {\em object} to signify either a space, or a pointed space, or a spectrum. Thus, a category of objects is one of the three categories $\Top$, $\Top_*$, or $\Spectra$. The definition of a functor from a topological category to a category of objects is long, but straightforward.

\begin{definition}\label{D: TopFunctor}
A functor from a small topological category $(\Ob,\Mor)$ to a category of objects consists of the following data.

(a) A fiberwise object $F$ over $\Ob$ (the fiber of $F$ at a point $x\in \Ob$ is what one would normally call the value of $F$ at $x$).

(b) A fiberwise map of objects over $\Mor$
$$ \alpha\colon s^*(F) \longrightarrow t^*(F)$$
where $s^*(F)$ and $t^*(F)$ are pullbacks of $F$ from $\Ob$ to $\Mor$ along the source and the target map respectively.

Moreover, the following needs to hold.

(i) (Unicity). Recall that $i:\Ob\longrightarrow \Mor$ is the map sending each object to its corresponding identity morphism. By definition of category, $s\circ i = t \circ i = \Id_{\Ob}$. Thus, $F=i^*(s^*(F))=i^*(t^*(F))$. Using these identifications, we require that the fiberwise map of objects over $\Ob$
$$i^*(\alpha)\colon i^*(s^*(F)) \longrightarrow i^*(t^*(F))$$ is the identity map from $F$ to $F$.

(ii) (Composition law). Let $\Mor \times_{^t\Ob^s} \Mor$ be the pullback that we defined above. There are three natural maps from $\Mor \times_{^t\Ob^s} \Mor$ to $\Ob$. Let us call them $s_l$, $t_r$, and $m$. $s_l$ is the source map applied to the ``left'' $\Mor$, $t_r$ is the target map applied to the ``right'' $\Mor$, and $m$ is the ``middle projection'' map. $m$ can be described as either $t_l$ or $s_r$. Consider the pull-backs of $F$ along these maps, $s_l^*(F), m^*(F),$ and $t_r^*(F)$. These are fiberwise objects over $\Mor \times_{^t\Ob^s} \Mor$. There are 
fiberwise maps
$$\alpha\times_{^t\Ob^s} \Id_{\Mor}\colon s_l^*(F) \longrightarrow m^*(F)$$
$$\Id_{\Mor} \times_{^t\Ob^s} \alpha \colon m^*(F) \longrightarrow t_r^*(F)$$
Notice also that $s_l=s\circ c$ and $t_r=t\circ c$, as maps from $\Mor \times_{^t\Ob^s} \Mor$ to $\Ob$. It follows that there are equalities
$$c^*(s^*(F))= s_l^*(F)$$ and
$$c^*(t^*(F))= t_r^*(F)$$
Finally, observe that there is a fiberwise map
$$c^*(\alpha)\colon c^*(s^*(F))\longrightarrow c^*(t^*(F))=t_r^*(F)$$
Now we are ready to state the composition law: we require that the following square diagram of fiberwise objects over $\Mor \times_{^t\Ob^s} \Mor$ commutes.
$$\begin{array}{ccc}
s_l^*(F)&  \longrightarrow&  m^*(F) \\
 =\downarrow \mbox{  }& & \downarrow \\
c^*(s^*(F)) &\longrightarrow &  t_r^*(F), \end{array}$$
where the maps are as defined above. This ends Definition~\ref{D: TopFunctor}.
\end{definition}

\begin{definition}\label{D: NatTrans}
Let $F$, $G$ be two functors defined in the sense of Definition \ref{D: TopFunctor}. A {\em natural transformation} from $F$ to $G$ is a fiberwise map (over $\Ob$) from $F$ to $G$ that is compatible with the structure maps $\alpha$. \end{definition}

The space of natural transformations from $F$ to $G$ will be denoted $\nat(F,G)$, or $\nat_{\calC}(F,G)$, or $\nat_{x\in\Ob}(F(x),G(x))$. It is topologized as a subspace of the space of all fiberwise maps from $F$ to $G$.

Of course, when $\calC$ is a discrete category, Definitions \ref{D: TopFunctor} and \ref{D: NatTrans} are equivalent to the standard ones.

\subsection{Homotopy limits.} Let $\calC=(\Ob,\Mor)$ be a topological category. Let $F$ be a fiberwise object over $\Ob$, defining a functor from $\calC$ to $\Top, \Top_*$, or $\Spectra$. Recall that $\Nerve\calC_i$ is the space of $i$-dimensional simplices in the simplicial nerve of $\calC$. Let $t_r(i):\Nerve\calC_i\longrightarrow \Ob$ be the target map applied to the ``rightmost'' $\Mor$ factor of $\Nerve\calC_i$. Let $t_r(i)^*(F)$ be the pull-back of $F$ along $t_r(i)$. Let $\Gamma_i:=\Gamma(t_r(i)^*(F))$ be the space (or spectrum) of sections of $t_r(i)^*(F)$, as $i=0,1,\ldots$. 
It is not difficult to see that that spaces (resp. spectra) $\Gamma_i$, $i=0,1, 2, \ldots$ fit into a cosimplicial space (resp. spectrum), which in the case when $\calC$ is a discrete category is the standard cosimplicial object that is used to define the homotopy limit of $F$ \cite{BK}. For a general topological category $\calC$, we define $\underset{\calC}{\holim}F$ to be the (derived) total space (or spectrum) of this cosimplicial object.


In practice, we will want to make some assumptions that will guarantee that homotopy limits behave homotopically as they should. Thus, we will assume that the spaces $\Nerve\calC_i$ are disjoint unions of compact ENRs. For a functor $F$ from $\calC$ to a category of objects, we will want to assume that the associated fiberwise object over $\Ob$ is a homotopy bundle, or, at worst, a stratified homotopy bundle. The following two lemmas follow from the corresponding invariance results on spaces of sections of homotopy bundles (Corollary~\ref{C: BaseInvariant}, Lemma~\ref{L: FiberInvariant} and Lemma~\ref{L: Excision}), and the homotopy properties of totalizations of cosimplicial spaces.

\begin{lemma}\label{L: CategoryInvariant}
Let $\calC=(\Ob(\calC),\Mor(\calC))$, $\calD=(\Ob(\calD),\Mor(\calD))$ be two small topological categories such that all the spaces in their simplicial nerves are disjoint unions of compact ENRs. Let $F$ be a functor from $\calD$ to $\Top, \Top_*$, or $\Spectra$. Suppose that as a fiberwise object over $\Ob(\calD)$, $F$ is a homotopy bundle. Let $\phi\colon \calC\longrightarrow \calD$ be a functor that induces a degree-wise homotopy equivalence of simplicial nerves. Then the induced map
$$\underset{\calD}{\holim} F\longrightarrow \underset{\calC}{\holim} F\circ \phi$$
is a homotopy equivalence.
\end{lemma}

\begin{lemma}\label{L: FunctorInvariant}
Let $\calC$ be a topological category and let $F, G$ be two functors from $\calC$ to one of the three standard categories. Make the same assumptions on $\calC$, $F$, and $G$ as in Lemma~\ref{L: CategoryInvariant}. Let $\psi\colon F\longrightarrow G$ be a natural transformation such that the associated fiberwise map over $\Ob(\calC)$ is a homotopy equivalence on each fiber. Then the induced map
$$\underset{\calC}{\holim} F\longrightarrow \underset{\calC}{\holim} G$$
is a homotopy equivalence.
\end{lemma}

\subsection{Spaces of natural transformations.} \subsubsection{Spaces of homotopy natural transformations}
An important case of homotopy limits is the space of homotopy natural transformations between two functors. In this context, we will only need to consider functors defined on discrete categories, so let $\calC$ be a small discrete category. Let $F$ be a functor from $\calC$ to $\Top$, and let $G$ be a functor from $\calC$ to either $\Top$, $\Top_*$, or $\Spectra$. We have already discussed the space (or spectrum) of natural transformations $\nat(F,G)$. However it is well-known that the construction $\nat(F,G)$ does not have good homotopic properties, and one often wants to replace it with a ``derived'' space of natural transformation, which we will denote $\hnat(F,G)$. One standard way to construct $\hnat(F,G)$ uses the machinery of model categories. One constructs a Quillen model structure on the category of functors, and one defines $\hnat(F,G):=\nat(QF, RG)$, where $Q$ and $R$ denote cofibrant and fibrant replacements respectively. But in this paper, we use a more direct construction, proposed by  Dwyer and Kan in ~\cite{DK}. Dwyer and Kan defined the space of homotopy natural transformation as the homotopy limit over the ``twisted arrow category'' $^a\calC$. This is a category whose objects are morphisms $x\rightarrow y$ in $\calC$, and where a morphism $(x\rightarrow y)\longrightarrow (x_1\rightarrow y_1)$ is given by a ``twisted'' commutative square diagram
$$
\begin{CD}
x @<<< x_1 \\
@VVV @VVV \\
y @>>> y_1 \end{CD}.$$

Clearly, there is a covariant functor from the category $^a\calC$ to the target category of $G$ given on objects by $$(x\rightarrow y)\mapsto \map(F(x),G(y)).$$
\begin{definition} The space (or spectrum) of homotopy natural transformations is defined by the following formula
$$\hnat_{\calC}(F,G):= \underset{x\rightarrow y\in ^a\calC}{\holim}\map(F(x),G(y)).$$
\end{definition}
\begin{remark}
The space of strict natural transformations can be defined as the strict inverse limit
$$\operatorname{Nat}_{\calC}(F,G):=  \underset{x\rightarrow y\in ^a\calC}{\lim}\map(F(x),G(y)).$$
There is a natural map 
$$\operatorname{Nat}_{\calC}(F,G) \longrightarrow \hnat_{\calC}(F,G).$$
This map is an equivalence if $F$ is cofibrant and $G$ is fibrant in any model structure that one can put on the category of functors, but in general it is not an equivalence.

When we speak of an individual natural transformation, we will usually mean a strict natural transformation.
\end{remark}

Here is a yet another way to present the space of (homotopy) natural transformations as a (homotopy) limit. Let $\calC$ be a (discrete) category, and let $F:\calC\longrightarrow \operatorname{Top}$ be a functor from $\calC$ to spaces. Recall from the introduction that we define the topological category $\calC\ltimes F$ as follows: an object of $\calC\ltimes F$ is an ordered pair $(c,x)$ where $c$ is an object of $\calC$ and $x\in F(c)$. A morphism $(c,x)\longrightarrow (d,y)$ consists of a morphism $\alpha:c\longrightarrow d$ in $\calC$, such that $F(\alpha)(x)=y$. The space of objects of $\calC\ltimes F$ is topologized as the disjoint union of the spaces $F(c)$, as $c$ ranges over objects of $\calC$. The space of morphisms is topologized as the disjoint union of the spaces $F(c_0)$, where the union is taken over all morphisms $c_0\longrightarrow c_1$ in $\calC$. There is an obvious projection functor $\calC\ltimes F\longrightarrow \calC$. It follows that any functor $G$ on $\calC$ can be thought of as a functor on $\calC\ltimes F$ by composition with this projection functor.
We view the following lemma as providing an alternative definition of the space of natural transformations. 
\begin{lemma}
There is a natural equivalence
$$\hnat_{\calC}(F,G)\simeq \underset{\calC\ltimes F}{\holim} \, G$$
\end{lemma}
\begin{proof}
It is easy to check that the cosimplicial space defining the left hand side is the edge-wise subdivision of the cosimplicial space defining the right hand side. For the definition of edge-wise subdivision, and a proof of the fact that edge-wise subdivision preserves totalization, see \cite{GKW}, Section 2, especially Propositions 2.2 and 2.3.
\end{proof}

It is also clear by inspection that there is a homeomorphism
$$\nat_{\calC}(F,G)\simeq \underset{\calC\ltimes F}{\lim} \, G$$

\subsubsection{Spaces of twisted natural transformations.} \label{S: TwistedNat}
We would like to extend the notion of space of natural transformations in the following way. Suppose that $F\colon \calC \longrightarrow \Top$ is a functor from $\calC$ to spaces. Let $G:\calC\ltimes F\longrightarrow \calD$ be a functor from $\calC\ltimes F$ to spaces (or spectra), that does not necessarily factor through the projection functor $\calC\ltimes F \longrightarrow \calC$. We define the {\em space (or spectrum) of twisted natural transformations} from $F$ to $G$, denoted $\widetilde{\nat}_{\calC}(F,G)$, by the following formula.
$$\widetilde{\nat}_{\calC}(F,G):= \underset{\calC\ltimes F}{\lim} \, G.$$
Similarly, we define the homotopy version as follows
$$\widetilde{\hnat}_{\calC}(F,G):= \underset{\calC\ltimes F}{\holim} \, G.$$

In order to ensure good homotopical behavior, we will usually assume that the functor $F$ takes values in ENRs, and that the functor $G:\calC\ltimes F\longrightarrow \calD$ defines a homotopy bundle over the space of objects of $\calC\ltimes F$ (or, at worst, a stratified homotopy bundle). 

It is possible to describe $\widetilde{\hnat}_{\calC}(F,G)$, as a homotopy limit over $^a\calC$, analogously to the definition of the ordinary space of homotopy natural transformations. The difference is that instead of a homotopy limit of a diagram of mapping spaces, one obtains a homotopy limit of a diagram of spaces of sections, or twisted mapping spaces. Indeed, $G$ being a functor on $\calC\ltimes F$ means that for every object $c$ of $\calC$, there is a homotopy bundle (of spaces, pointed spaces, or spectra) $G|_c$ over $F(c)$. Let us assume that all the fibers of this bundle are homotopy equivalent  to a fixed object $G(c)$ (this holds automatically if $F(c)$ is path connected). This assumption is not at all essential, but anyway it will be satisfied in all the examples we will consider later. 
For a morphism, $h\colon c\rightarrow d$ in $\calC$, let $h^*(G|_d)$ be the pullback of $G|_d$ along $F(h)\colon F(c)\longrightarrow F(d)$. $h^*(G|_d)$ is a fiberwise object over $F(c)$, whose fibers are homotopy equivalent to $G(d)$. $G$ being a functor on $\calC\ltimes F$ implies, in particular, that there is a fiberwise map $G|_c\rightarrow h^*(G|_d)$, which is natural in $h$ in an evident sense. In fact, it is easy to see that $G$ being a functor on $\calC\ltimes F$ implies the existence of a functor from the twisted arrow category $^a\calC$ to $\calD$ given on objects by 
$$(c\rightarrow d)\mapsto \widetilde\map\left(F(c),G(d)\right),$$ such that
$$\widetilde{\hnat}_{\calC}(F,G)\simeq \underset{(c\to d)\in ^a\calC}{\holim} \widetilde\map\left(F(c),G(d)\right).$$

There is an analogous formula for the strict version. Namely, there is a homeomorphism 
$$\widetilde{\nat}_{\calC}(F,G)\cong \underset{(c\to d)\in ^a\calC}{\lim} \widetilde\map\left(F(c),G(d)\right).$$

\subsubsection{The relative version.} We will want to allow for a slightly more general version of spaces of twisted natural transformations. Let $F$ be a functor from $\calC$ to pairs of spaces. Let's say $F(c)=(F_1(c),F_0(c))$, where $F_0, F_1$ are functors from $\calC$ to $\Top$, taking values in ENRs, such that $F_0(c)$ is a closed sub-ENR of $F_1(c)$. Let $G$ be a functor from $\calC\ltimes F_1$ to pointed spaces, or spectra. We define 
$$\widetilde{\nat}_{\calC}\left((F_1,F_0); G \right):=\operatorname{fiber}\left(\widetilde{\nat}_{\calC}\left(F_1, G \right)\longrightarrow \widetilde{\nat}_{\calC}\left(F_0, G \right)\right).$$
Similarly, we define
$$\widetilde{\hnat}_{\calC}\left((F_1,F_0); G \right):=\hofiber\left(\widetilde{\hnat}_{\calC}\left(F_1, G \right)\longrightarrow \widetilde{\hnat}_{\calC}\left(F_0, G \right)\right).$$
\label{SSS: RelativeVersion} 
The following homotopy invariance lemma is a consequence of Lemma~\ref{L: CategoryInvariant}, Lemma~\ref{L: FunctorInvariant}, and excision (Lemma~\ref{L: Excision}).
\begin{lemma}\label{L: hNatInvariant}
\begin{enumerate}
\item Let $F=(F_1, F_0)$ and $F'=(F_1', F_0')$ be functors from $\calC$ to pairs of ENRs. Let $h\colon F'\longrightarrow F$ be a natural transformation. Suppose that the square diagram of functors
$$\begin{CD}
F_0' @>>> F_1' \\
@VhVV @VVhV \\
F_0 @>>> F_1 
\end{CD}$$
is a homotopy pushout. 
 Let $G$ be a functor from $\calC\ltimes F_1$ to pointed spaces, or spectra. $G$ can also be considered a functor from $\calC\ltimes F_1'$ via composition with $h$. Under these assumptions, the induced map
$$\widetilde{\hnat}_{\calC}\left((F_1,F_0), G\right)\longrightarrow \widetilde{\hnat}_{\calC}\left((F_1',F_0'), G \right)
$$
is a weak equivalence.
\item Let $=(F_1, F_0)$ be as above, and let $G_1$, $G_2$ be functors from $\calC\ltimes F_1$ to pointed spaces, or spectra. Let $h\colon G_1\longrightarrow G_2$ be a weak equivalence of functors. Then the induced map
$$\widetilde{\hnat}_{\calC}\left((F_1,F_0), G_1\right)\longrightarrow \widetilde{\hnat}_{\calC}\left((F_1,F_0), G_2\right)
$$
is a weak equivalence.
\end{enumerate}
\end{lemma}
\subsubsection{Well filtered categories.} We will now study limits and homotopy limits of functors on $^a\calC$ for discrete categories $\calC$ that admit a filtration with certain nice properties. This will tell us something about spaces of (homotopy, twisted) natural transformations between functors on such categories.

Let us first consider the following basic situation. 
\begin{definition}\label{D: NiceExtension}
Let $\calC$ be a category, and let $\calC_0\subset \calC$ be a subcategory. We say that $\calC$ is a nice extension of $\calC_0$ if 
\begin{enumerate}
\item there are no morphisms from objects of $\calC\setminus \calC_0$ to objects of $\calC_0$, and
\item all morphisms between objects of $\calC\setminus \calC_0$ are isomorphism.
\end{enumerate}
\end{definition}
The main point that we want to make is contained in the following elementary proposition.
\begin{proposition}\label{P: ExtensionLimits}
Let $\calC$ be a nice extension of $\calC_0$. Let $H$ be a functor from $^a\calC$ to the category of spaces, or pointed spaces, or spectra. Then there is a pullback square
$$\begin{CD}
\underset{^a\calC}{\lim} H @>>> \underset{[d]\in[\calC\setminus \calC_0]}{\prod} H(d\stackrel{=}\to d)^{\Sigma(d)} \\
@VVV @VVV \\
\underset{^a\calC_0}{\lim} H @>>> \underset{[d]\in[\calC\setminus \calC_0]}{\prod}\left(\underset{(c\to d)\in(\calC_0\downarrow d)^{\op}}{\lim} H(c\to d)\right)^{\Sigma(d)}
\end{CD}.$$
where the products in the right column are indexed over the same set of representatives of isomorphism classes of $\calC\setminus \calC_0$, $\calC_0\downarrow d$ is the evident category of arrows of the form $c\to d$, and $\Sigma(d)$ is the group of automorphisms of $d$ in $\calC$, which is also the group of automorphisms of the identity arrow $d\stackrel{=}{\to} d$ in $^a\calC$.

There also is a diagram of the following form, where the square is a homotopy pullback, and the arrows marked by $\simeq$ are homotopy equivalences.
 $$\begin{CD}
\underset{^a\calC}{\holim} H @>>>  Y@>\simeq>> \underset{[d]\in[\calC\setminus \calC_0]}{\prod} H(d\stackrel{=}\to d)^{h\Sigma(d)} \\
@VVV @VVV \\
X @>>> \underset{[d]\in[\calC\setminus \calC_0]}{\prod}\left(\underset{(c\to d)\in(\calC_0\downarrow d)^{\op}}{\holim} H(c\to d)\right)^{h\Sigma(d)} \\
@V\simeq VV \\
\underset{^a\calC_0}{\holim} H
\end{CD}.$$
\end{proposition}
\begin{proof}
Let $\calC_0\downarrow \calC$ be the fulI subcategory of $^a\calC$ consisting of arrows $c\to d$ where $c\in \calC_0$. Let $\calC\downarrow \calC^-$ be the full subcategory consisting of arrows $c\to d$ where $d\in \calC\setminus \calC_0$. Let $\calC_0\downarrow \calC^- =\calC_0\downarrow \calC\cap\calC\downarrow \calC^-$. We have a commutative square of categories
$$\begin{CD}
\calC_0\downarrow \calC^- @>>> \calC_0\downarrow \calC \\
@VVV @VVV\\
\calC\downarrow \calC^- @>>> ^a\calC \end{CD}.$$
It is not difficult to see that since $\calC$ is a nice extension of $\calC_0$, the above square induces both a degree-wise pushout and a degree-wise homotopy pushout square of simplicial nerves. It follows that there is a strict pullback square
$$\begin{CD}
\underset{^a\calC}{\lim} H @>>> \underset{\calC\downarrow \calC^- }{\lim} H \\
@VVV @VVV \\
\underset{ \calC_0\downarrow \calC }{\lim} H @>>> \underset{\calC_0\downarrow \calC^-}{\lim} H
\end{CD}$$
and an analogous homotopy pullback square, obtained by replacing all the limits by homotopy limits.
We now proceed to analyze three corners of these two squares. Let us begin with the corner defined by limit over $\calC_0\downarrow \calC$. Clearly, there is an inclusion of categories $^a\calC_0\hookrightarrow \calC_0\downarrow \calC$. Let $c\to d$ be an object of $\calC_0\downarrow \calC$. Let $^a\calC_0\downarrow (c\to d)$ be the over category.
\begin{lemma} 
For every object $c\to d$ of $\calC_0\downarrow \calC$, the category $^a\calC_0\downarrow (c\to d)$ has a contractible nerve.
\end{lemma}
\begin{proof}
An object of $^a\calC_0\downarrow (c\to d)$  is the same as a commutative square of objects of $\calC$
$$\begin{CD}
c' @<<< c \\
@VVV @VVV \\
d' @>>> d \end{CD}$$
where $d'$ is an object of $\calC_0$.
Consider the full subcategory of $^a\calC_0\downarrow (c\to d)$ consisting of objects of the form 
$$\begin{CD}
c @<=<< c \\
@VVV @VVV \\
d' @>>> d \end{CD}.$$
It is easy to see that this is a final subcategory of $^a\calC_0\downarrow (c\to d)$, and so its nerve is homotopy equivalent to the nerve of the entire category. On  the other hand, this subcategory has an initial object, namely $c\stackrel{=}\to c$, and so its nerve is contractible. This proves the lemma.
\end{proof}
It follows that the map $$\underset{\calC_0\downarrow \calC}{\lim} H\longrightarrow \underset{^a\calC_0}{\lim} H$$
is a homeomorphism, while the analogous map of homotopy limits is a homotopy equivalence. This takes care of the left half of the proposition.

It remains to analyze the categories $\calC_0\downarrow \calC^-$, $\calC\downarrow \calC^-$, and (homotopy) limits over them. Without loss of generality, we may assume that $\calC\setminus\calC_0$ has one object in each isomorphism class, because the categories obtained from $\calC_0\downarrow \calC^-$ and $\calC\downarrow \calC^-$ by restricting to just one representative in each isomorphism class of $\calC\setminus\calC_0$ are equivalent to the whole category. Since both categories $\calC_0\downarrow \calC^-$ and $\calC\downarrow \calC^-$ split as a disjoint union indexed by (isomorphism classes of) objects $\calC\setminus\calC_0$, we may as well assume that $\calC\setminus\calC_0$ has just one object. Let us call it $d$. Let $\calC_0 \downarrow d$ and $\calC \downarrow d$ be, once again, the over categories. Thus, for example, an object of $\calC_0\downarrow d$ is an arrow $c\to d$, and a morphism is a commuting triangle of an evident kind. It is easy to see that the group $\Sigma(d)$ acts on $\calC_0 \downarrow d$ and $\calC \downarrow d$, and in fact there are isomorphisms of categories
$$\calC_0\downarrow \calC^-\cong \left(\calC_0\downarrow d\right)^{\op}\rtimes \Sigma(d)$$
and
$$\calC\downarrow \calC^-\cong \left(\calC\downarrow d\right)^{\op}\rtimes \Sigma(d).$$
It follows that for any functor $H$ on $^a\calC$, there are homeomorphisms
$$\underset{\calC_0\downarrow \calC^-}{\lim} H\cong \left(\underset{(\calC_0\downarrow d)^{\op}}{\lim} H\right)^{\Sigma(d)}$$
and 
$$\underset{\calC\downarrow \calC^-}{\lim}\, H\cong \left(\underset{(\calC\downarrow d)^{\op}}{\lim} H\right)^{\Sigma(d)}$$
as well as analogous homotopy equivalences obtained by replacing all limits with homotopy limits. Since the identity morphism $d\stackrel{=}{\to}d$ is a final object of ${\calC\downarrow d}$, the proposition follows.
\end{proof}
\begin{corollary}\label{C: InductiveModel}
Let $\calC$ be a nice extension of $\calC_0$. Let $F=(F_1, F_0)$ be a functor from $\calC$ to pairs of topological spaces. Let $G$ be a functor from $\calC\ltimes F_1$ to the category of pointed spaces, or spectra, so that the spaces (or spectra) $\underset{\calC}{\widetilde{\nat}}(F, G)$ and $\underset{\calC}{\widetilde{\hnat}}(F, G)$ are defined. Then there is a strict pullback diagram
$$\begin{CD}
\underset{\calC}{\widetilde{\nat}} (F, G) @>>> \underset{[d]\in[\calC\setminus \calC_0]}{\prod} \widetilde{\map}(F(d), G(d))^{\Sigma(d)} \\
@VVV @VVV \\
\underset{\C_0}{\widetilde{\nat}} (F, G) @>>> \underset{[d]\in[\calC\setminus \calC_0]}{\prod}{\widetilde{\map}}\left(\underset{(c\to d)\in(\calC_0\downarrow d)}{\operatorname{colim}} F(c), G( d)\right)^{\Sigma(d)}
\end{CD}$$
where the left vertical map is the evident restriction map, and the right vertical map is induced by maps, for each $d$
$$\underset{(c\to d)\in(\calC_0\downarrow d)}{\operatorname{colim}} F(c)\longrightarrow F(d).$$

There also is a homotopy pullback square, obtained essentially be replacing all limits and colimits by homotopy limits and homotopy colimits.
\end{corollary}

Let us now generalize the situation of a nice extension as follows.
\begin{definition}\label{D: NicelyFiltered}
Let $\calC$ be a small discrete category. We say that $\calC$ is nicely filtered, if there exists a chain of full subcategories
$$\emptyset=\calC_0\subset\calC_1\subset\cdots\subset \calC_n=\calC.$$
such that $\calC_i$ is a nice extension of $\calC_{i-1}$ for all $1\le i\le n$.
\end{definition}
Let $F=(F_1, F_0)$ be a functor from $\calC$ to pairs of spaces, and let $G$ be a functor from $\calC\ltimes F_1$ to pointed spaces, or spectra. Clearly, if $\calC$ is nicely filtered, then there is a tower of spaces (or spectra) of natural transformations
$$\underset{\calC}{\widetilde{\nat}}(F, G)\longrightarrow \cdots\longrightarrow \underset{\calC_i}{\widetilde{\nat}}(F, G)\longrightarrow \cdots$$
as well as an analogous tower of spaces of homotopy natural transformations. Corollary~\ref{C: InductiveModel} provides an inductive description of these towers.

Recall that a map of pairs $h\colon (A_1, A_0)\longrightarrow (B_0, B_1)$ is called a cofibration (resp. a weak equivalence) if the induced map $$\tilde h\colon\operatorname{colim}(A_1\leftarrow A_0 \rightarrow B_0)\longrightarrow B_1$$ is a cofibration (resp. a weak equivalence). By the cofiber of the map $h$ we will mean the cofiber of the map $\tilde h$.
\begin{definition}\label{D: EssentiallyCofibrant}
Let $\calC$ be a nicely filtered category. Let $F=(F_1, F_0)$ be a functor from $\calC$ to pairs of spaces. We say that $F$ is {\em essentially cofibrant} if for every object $d$ of $\calC_i\setminus\calC_{i-1}$
\begin{enumerate}
\item the map 
$$\underset{(c\to d)\in \calC_{i-1}\downarrow d}{\operatorname{colim}}F(c)\longrightarrow F(d)$$
is a cofibration
\item The group $\Sigma(d)$ acts freely (in the pointed sense) on the cofiber of the above map.
\end{enumerate}
\end{definition}
The following lemma follows from Corollary~\ref{C: InductiveModel} and induction
\begin{lemma}\label{L: NatFromCofibrant}
Let $\calC$ be a nicely filtered subcategory. Let $F$ be an essentially cofibrant functor from $\calC$ to pairs of spaces. Let $G$ be, as usual, a functor from $\calC\ltimes F$ to pointed spaces, or spectra. Then for all subcategories $\calC_i\subset \calC$ involved in the nice filtration of $\calC$, the natural map
$$ \underset{\calC_i}{\widetilde{\nat}}(F, G)\longrightarrow \underset{\calC_i}{\widetilde{\hnat}}(F, G)$$
is a homotopy equivalence.
\end{lemma}

\section{Partitions}\label{S: Partitions}
Formally speaking, a {\em partition} is an ordered pair $(s,\Lambda)$, where $s$ is a finite set, and $\Lambda$ is an equivalence relation on $s$. In this situation, we say that $\Lambda$ is a partition of $s$, and $s$ is the {\em support} of $\Lambda$. We will use $\supp(\Lambda)$ to denote the support of $\Lambda$, and we let
$\comp(\Lambda)$ denote the set of equivalence classes (which we will call {\em components}) of $\Lambda$. Often, when the support of a partition does not need to be identified explicitly, we will denote the partition $(s,\Lambda)$ simply by $\Lambda$.

A few terms that we will use: we will say that a partition is {\em irreducible} if none of its components is a singleton. Every set $s$ has two obvious extreme partitions: the {\em discrete} partition (which we will also refer to as the {\em trivial} partition), and the {\em indiscrete} one. We leave it to the reader to guess what they are.

There is a close relationship between partitions and surjective maps of sets. A surjection $\alpha:s\twoheadrightarrow c$ determines a partition of $s$ (the components are the inverse images of elements of $c$), and conversely a partition of $\Lambda$ of $s$  determines a surjection $\alpha:s\twoheadrightarrow c$ uniquely up to an isomorphism of $c$. We will say that $\alpha$ represents $\Lambda$. 

In this paper we will need consider several kinds of morphisms between partitions. Let us begin with the familiar construction of the poset of partitions of a fixed set. Let $s$ be a fixed set, and  
suppose $\Lambda_1$ and $\Lambda_2$ are two partitions of the same set $s$. We say that $\Lambda_2$ is a {\em refinement} of $\Lambda_1$ if every component of $\Lambda_2$ is a subset of a component of $\Lambda_1$. It is obvious that refinement is a partial order relation, and we say that $\Lambda_1\le \Lambda_2$ if $\Lambda_2$ is a refinement of $\Lambda_1$ (in this case we also say that $\Lambda_1$ is a coarsening of $\Lambda_2$). The collection of partitions of a fixed set $s$ forms a poset under refinement, and we will some times consider it as a category, in the usual way: there is a (unique) morphism $\Lambda_1\longrightarrow\Lambda_2$ in this category if and only if $\Lambda_1\le \Lambda_2$. 

Refinements have a characterization in terms of surjections, given in the following elementary lemma.
\begin{lemma}\label{L: Refinements}
Let $\Lambda, \Lambda'$ be partitions of $s$. Then $\Lambda\le\Lambda'$ if and only if for any surjections $\alpha: s \twoheadrightarrow  c$ and $\alpha: s \twoheadrightarrow  c'$, representing $\Lambda$ and $\Lambda'$ respectively,  there exist a commutative square
$$\begin{array}{ccc}
s & \stackrel{\alpha}{\twoheadrightarrow} & c \\
\!\!\!\!=\downarrow & & \uparrow \\
s & \stackrel{\alpha'}{\twoheadrightarrow} & c' \end{array}$$
where the left vertical map is the identity. Note that the right vertical map is pointing in the ``wrong'' direction.
\end{lemma}

\begin{definition}\label{D: PartitionPosets} 
Let $\Lambda$ be a partition of $s$. Define $\calP(\Lambda)$ to be the poset of all partitions $\Lambda'$ of $s$ such that $\Lambda\le \Lambda'$. 

If $\Lambda$ the the initial (i.e., the indiscrete) partition of $s$, then $\calP(\Lambda)$ is the familiar poset of partitions of $s$. In this case we will use notation $\calP_s:=\calP(\Lambda)$. If, furthermore, $s$ is the standard set with $n$ elements, we will denote $\calP_s$ by $\calP_n$. 
\end{definition}

For a small category $\calC$, we use $|\calC|$ to denote the geometric realization of the simplicial nerve of $|\calC|$. Since $\calP(\Lambda)$ has both an initial and a final object, $|\calP(\Lambda)|$ is a contractible complex. $|\calP(\Lambda)|$ can be thought of as polyhedron with a boundary. In fact, $\partial|\calP(\Lambda)|$ can be defined explicitly as the realization of the simplicial subset of the nerve of $\calP(\Lambda)$ consisting of simplices that do not contain both the initial and the final object of $\calP(\Lambda)$ as a vertex.
\begin{definition}\label{D: SpaceOfTrees}
Let $T_\Lambda:=|\calP(\Lambda)|/\partial|\calP(\Lambda)|$. In the case when $\Lambda$ is the indiscrete partition of the standard set with $n$ elements, we use $T_n$ to denote $T_\Lambda$.
\end{definition}

It is well-known that $$T_n\simeq \bigvee_{(n-1)!} S^{n-1}.$$
If $\Lambda$ is a partition with components $(n_1,\ldots,n_i)$, it is easy to see that there is an isomorphism of posets
$$\calP(\Lambda)\cong \calP_{n_1}\times\cdots\times  \calP_{n_i}$$
From here, in turn, it follows easily that there is a homeomorphism
$$T_\Lambda \cong T_{n_1}\wedge\cdots\wedge T_{n_i}$$
This determines the homotopy type of $T_\Lambda$ for a general $\Lambda$. Note in particular that $T_\Lambda$ is homotopy equivalent to a wedge sum of spheres of dimension $n-i$, where $n$ is the size of the support, and $i$ is the number of components of $\Lambda$. The number $n-i$ is an important invariant of $\Lambda$, and deserves a separate definition. 
\begin{definition}\label{D: Excess}
Let $\Lambda$ be a partition of $s$. Let $c$ be the set of  components of $\Lambda$. The number $|s|-|c|$ is called the {\em excess} of $\Lambda$. We will denote it by $\excess(\Lambda)$.
\end{definition}

Now we will introduce another notion of morphisms between partitions, one that may involve partitions of different sets. 
Suppose $\Lambda$ is a partition of $s$, and let $f:s\longrightarrow s'$ be a map of sets. Together, $f$ and $\Lambda$ determine a relation $G_f$ on $s'$ in an obvious way: $(x',y')\in G_f$ if there exist $x\in f^{-1}(x')$ and $y\in f^{-1}(y')$ such that $x$ and $y$ are in the same equivalence class of $\Lambda$. Clearly, $G_f$ is reflexive and symmetric, but in general not transitive, thus is not an equivalence relation. We use $f(\Lambda)$ to denote the equivalence relation on $s'$ generated by $G_f$. The following elementary lemma provides a convenient characterization of $f(\Lambda)$.

\begin{lemma}\label{L: Pushout}
Let $\Lambda, \Lambda'$ be partitions of $s,s'$ respectively. Let $f:s\longrightarrow s'$ be a map of sets. Then $f(\Lambda)=\Lambda'$ if and only if for some (and therefore any) surjections $\alpha\colon s \twoheadrightarrow c$, $\alpha'\colon s' \twoheadrightarrow c'$, representing $\Lambda$ and $\Lambda'$ respectively, there exists a map $h:c\longrightarrow c'$ such that the following is a strict pushout square
$$\begin{array}{ccc}
s & \stackrel{\alpha}{\twoheadrightarrow} & c \\
\!\!\!\! f \! \downarrow & & \,\,\,\,\downarrow \!h \\
s' & \stackrel{\alpha'}{\twoheadrightarrow} & c'
\end{array}
$$
\end{lemma}

\begin{definition}
Let $(s,\Lambda)$ and $(s',\Lambda')$ be two partitions. A {\em fusion} of $\Lambda$ into $\Lambda'$, denoted $f:\Lambda\longrightarrow\Lambda'$, is a map of sets (which we will denote with the same letter) $f:s\longrightarrow s'$ such that $f(\Lambda)=\Lambda'$. 
\end{definition}

It is obvious that composition of fusions is a fusion, and thus fusions of partitions form a category. Let $\calF$ denote the category of all partitions and fusions between them. It is also easy to see that fusions and refinements satisfy the following property
\begin{lemma}
Let $f\colon s\to s'$ be a map of sets. Let $\Lambda, \Lambda'$ be partitions of $s$. If $\Lambda\le \Lambda'$ then $f(\Lambda)\le f(\Lambda')$.
\end{lemma}
It follows that there is a functor from the category of sets to the category of small categories that is defined on objects by the formula $s\mapsto \calP_s$, and is defined on morphisms by the formula $\Lambda \to f(\Lambda)$ where $f\colon s\to s'$ is a map of sets, and $\Lambda$ is an object of $\calP_s$. There a well known construction due to Grothendieck that encodes this situation in a single category. We call the resulting category $\calS$. Explicitly, $\calS$ is defined as follows.
\begin{definition}\label{D: CategoryOfPartitions}
$\calS$ is the category whose objects are all partitions $(s, \Lambda)$. A morphism $$f\colon (s,\Lambda)\longrightarrow (s',\Lambda')$$ in $\calS$ consists of a map of sets (which we denote with the same letter) $f\colon s\to s'$ such that $f(\Lambda)\le \Lambda'$.
\end{definition}
We refer of $\calS$ as ``the'' category of partitions. We often denote an object of $\calS$ by $\Lambda$, rather than $(s, \Lambda)$, $s$ being understood. Any morphism $f:\Lambda_1\longrightarrow \Lambda_2$ can be factored canonically as a fusion followed by a refinement. Indeed, the factorization is
$$\Lambda_1 \longrightarrow f(\Lambda_1) \longrightarrow \Lambda_2.$$ For a general morphism $f:\Lambda_1\longrightarrow\Lambda_2$, we say that $f$ is a fusion if $\Lambda_2=f(\Lambda_1)$, and we say that $f$ is a refinement if the underlying map of sets is the identity map.

The notion of excess can now be refined as follows. Let $V$ be a Euclidean space. Let $\Lambda$ be a partition, represented by a surjection $\alpha:s\twoheadrightarrow c$. $\alpha$ induces an injective linear homomorphism $V^c\hookrightarrow V^s$. Clearly, the quotient, which we will denote by $V^s/V^c$, has dimension $\excess(\Lambda)\dim(V)$. Moreover, suppose $f:\Lambda\longrightarrow\Lambda'$ is a fusion. Let $\Lambda$ and $\Lambda'$ be represented by surjections $s\twoheadrightarrow c$ and $s'\twoheadrightarrow c'$ respectively. By Lemma \ref{L: Pushout}, there is a pushout square 
$$\begin{array}{ccc}
s & \stackrel{\alpha}{\twoheadrightarrow} & c \\
\!\!\!\! f \! \downarrow & & \,\,\,\,\downarrow \!h \\
s' & \stackrel{\alpha'}{\twoheadrightarrow} & c'
\end{array}
$$
Taking maps into $V$, we obtain a pullback square of vector spaces
\begin{eqnarray}
V^{c'} & \hookrightarrow & V^{s'} \nonumber \\
\downarrow & &\downarrow \label{D: PullBack}\\
V^{c} & \hookrightarrow & V^{s}\nonumber \end{eqnarray}
Taking cokernels in the horizontal direction, and denoting them $V^{\excess(\Lambda')}$ and $V^{\excess(\Lambda)}$ respectively, we obtain a homomorphism $f^*:V^{\excess(\Lambda')}\longrightarrow V^{\excess(\Lambda)}$. Since the diagram \eqref{D: PullBack} is a pullback square, it follows that $f^*$ is a monomorphism. In fact,  we have constructed a {\em contravariant functor} from the category $\calF$ to the category of vector spaces and monomorphisms, given on objects by $\Lambda \mapsto V^{\excess(\Lambda)}$. In particular, if $f:\Lambda\longrightarrow\Lambda'$ is a fusion, then $\excess(\Lambda)\ge \excess(\Lambda')$. 

Notice also that if $\Lambda \le \Lambda'$ then $\excess(\Lambda)\ge \excess(\Lambda')$. It follows that for any morphism $f\colon \Lambda\to \Lambda'$ in $\calS$, $\excess(\Lambda)\ge \excess(\Lambda')$.
However, a refinement $\Lambda \le \Lambda'$ does not naturally give rise to an injective homomorphisms $V^{\excess(\Lambda')}\longrightarrow V^{\excess(\Lambda)}$. Instead, it gives rise to a surjective
homomorphism $V^{\excess(\Lambda)}\longrightarrow V^{\excess(\Lambda')}$ (this follows easily from Lemma~\ref{L: Refinements}).

A certain subcategory of fusions plays an important role in this paper. 
\begin{definition}
Let $f:\Lambda\longrightarrow \Lambda'$ be a fusion. We say that $f$ is a {\em strict fusion} if the induced monomorphism $f^*:V^{\excess(\Lambda')}\longrightarrow V^{\excess(\Lambda)}$ is an isomorphism for some (and hence any) non-zero Euclidean space $V$.\footnote{This definition differs slightly from the one given in the introduction, where we used $\Z$ instead of a Euclidean space $V$, but it is easy to see that the two definitions are equivalent}
\end{definition}
In particular, there can only be a strict fusion between partitions of the same excess. 
We need to record some properties of strict fusions. Clearly, the composition of strict fusions is again a strict fusion. This statement has the following converse.
\begin{lemma}\label{L: StrictFusionsHereditary}
Let $f:\Lambda\longrightarrow \Lambda'$,  and $g:\Lambda'\longrightarrow \Lambda''$ be fusions. If $g\circ f$ is a strict fusion,  then so are $f$ and $g$.
\end{lemma}
\begin{proof}
Let $V$ be a Euclidean space, and consider the induced homomorphisms
$$V^{e(\Lambda'')}\stackrel{g^*}{\hookrightarrow}V^{e(\Lambda')} \stackrel{f^*}{\hookrightarrow} V^{e(\Lambda)}$$
Since these homomorphisms are monomorphisms, and the composed homomorphism is, by our assumption, an isomorphism, it follows that $g^*$ and $f^*$ are isomorphisms.
\end{proof}
Let $f:s\longrightarrow s'$ be a map of sets. We say that $f$ is {\em elementary} if $f$ brings together two elements $x,y\in s$, and is otherwise injective. Suppose that $\Lambda$ is a partition of $s$, and $x,y\in s$. Let $f:s\longrightarrow s'$ be an elementary map identifying $x$ and $y$. It is easy to see that the induced morphism of partitions $\Lambda\longrightarrow f(\Lambda)$ is a strict fusion if and only if $x$ and $y$ belong to different components of $\Lambda$. We call a strict fusion induced by an elementary map, {\em an elementary strict fusion}. 
The following lemma gives another convenient characterization of strict fusions.
\begin{lemma}\label{L: StrictFusionsDecompose}
Let $f:\Lambda\longrightarrow \Lambda'$ be a fusion. Then $f$ can be written (not uniquely) as a composition of elementary fusions. $f$ is a strict fusion if and only if it can be written as a composition of elementary strict fusions.
\end{lemma}
\begin{proof}
First consider $f$ as a map of sets $f:s\longrightarrow s'$. Clearly, $f$ can be written as a composition of elementary maps, say $f=f_i\circ f_{i-1} \circ \cdots \circ f_1$. Each $f_j$ gives rise to an elementary fusion, and this in fact gives a decomposition of $f$ as a product of elementary fusions. From Lemma~\ref{L: StrictFusionsHereditary} we know that if $f$ is a strict fusion, each of the $f_j$-s has to be a strict fusion as well.
\end{proof}

The following proposition provides some ``topological'' characterizations of strict fusions.
\begin{proposition}\label{P: HomotopyPushout}
Let $f:\Lambda\longrightarrow \Lambda'$ be a fusion. By Lemma \ref{L: Pushout}, $f$ can be represented by a strict pushout square
\begin{equation}\label{Eq: Square}
\begin{array}{ccc}
s & \stackrel{\alpha}{\twoheadrightarrow} & c \\
\downarrow & & \downarrow \\
s' & \stackrel{\alpha'}{\twoheadrightarrow} & c'
\end{array}.
\end{equation}
Also, let $\Cyl_\alpha$, $\Cyl_{\alpha'}$ be the mapping cylinders of $\alpha, \alpha'$ respectively. Clearly, \eqref{Eq: Square} gives rise to a map of quotient spaces $\Cyl_\alpha/s \longrightarrow \Cyl_{\alpha'}/s'$. Note that $\Cyl_\alpha/s$ is equivalent to a wedge sum of $\excess(\Lambda)$ copies of $S^1$. The following are equivalent:
\begin{enumerate}
\item $f$ is a strict fusion.
\item The diagram~\eqref{Eq: Square} above is a homotopy pushout square.
\item The induced map $\Cyl_\alpha/s \longrightarrow \Cyl_{\alpha'}/s'$ is a homotopy equivalence.
\end{enumerate}
\end{proposition}
\begin{proof}
First, let us prove the implication (1) $\Rightarrow$ (2). Suppose that $f$ is a strict fusion. Then $f$ can be written as a composition of elementary strict fusions, and so it is enough to prove that (3) holds in the case when $f$ is an elementary stict fusion. This is a very easy calculation and is left to the reader.
The implication (2) $\Rightarrow$ (3) follows by elementary properties of homotopy pushouts. It remains to prove that (3) $\Rightarrow$ (1). We will prove the contrapositive statement.  Suppose that $f$ is not a strict fusion. By definition, this means that the induced homomorphism 
$V^{\excess(\Lambda')} \longrightarrow V^{\excess(\Lambda)}$ is not an isomorphism. But $V^{\excess(\Lambda)}$ is naturally isomorphic to $\HH^1(\Cyl_{\alpha}/s; V)$, so it the map in part (3) of the proposition does not induce an isomorphism in cohomology, thus is not a homotopy equivalence.
\end{proof}
\begin{corollary}\label{C: SimplyConnected}
Let $\Lambda$ be a partition represented by a surjection $\alpha\colon s\twoheadrightarrow c$. Let $f\colon\Lambda\longrightarrow\Lambda'$ be a fusion. Then $f$ is a strict fusion if and only if the homotopy pushout of the diagram
\begin{equation}\label{Eq: Graph}
s'\stackrel{f}{\leftarrow} s\stackrel{\alpha}{\twoheadrightarrow} c\end{equation}
is simply connected.
\end{corollary}
\begin{proof}
Let $X$ be the homotopy pushout of \eqref{Eq: Graph}, and let $Y$ be the strict pushout of the same diagram. There is a natural map $X\longrightarrow Y$, which is clearly a bijection on $\pi_0$. By Proposition~\ref{P: HomotopyPushout}, $f$ is a strict fusion if and only this map is an equivalence. Since $X$ is clearly a $1$-dimensional complex, and $Y$ is a $0$-dimensional complex, this map is an equivalence if and only if $X$ is simply connected. 
\end{proof}

We will also need the following couple of lemmas about the interaction of strict fusions and refinements.
\begin{lemma} \label{L: RefinementPreservesStrictFusion}
Suppose $\Lambda\le \Delta$ and $f:\Lambda\longrightarrow\Lambda'$ is a strict fusion. Then the induced morphism
$\Delta\longrightarrow f(\Delta)$ is a strict fusion.
\end{lemma}
\begin{proof}
Using Lemma~\ref{L: StrictFusionsDecompose} and induction, it is easy to see that it is enough to consider the case when $f$ is an elementary strict fusion. In this case, $f$ brings together two elements $x,y\in \supp(\Lambda)$, belonging to different components of $\Lambda$. Since $\Delta$ is a refinement of $\Lambda$, it follows that $x$ and $y$ belong to different components of $\Delta$, and therefore the induced map $\Delta\longrightarrow f(\Delta)$ is a strict fusion.
\end{proof}

It follows from the above lemma that while fusions preserve refinements, strict fusions preserve strict refinements. More precisely, we have the following lemma.

\begin{lemma}\label{L: PreserveStrictRefinements}
Suppose $\Lambda < \Delta$ and a morphism $f:\Lambda\longrightarrow f(\Lambda)$ is a strict fusion. Then $f(\Lambda)<f(\Delta)$.
\end{lemma}
\begin{proof}
It is easy to see that if $\Lambda\le \Delta$ then $\Lambda < \Delta$ if and only if $\excess(\Lambda)>\excess(\Delta)$. By Lemma~\ref{L: RefinementPreservesStrictFusion}, the morphism $\Delta\longrightarrow f(\Delta)$ is a strict fusion. 
We have the following relations between the excesses of these partitions
$$\excess(f(\Lambda))=\excess(\Lambda) > \excess(\Delta)=\excess (f(\Delta)).$$
Here the first equality holds because $f$ is a strict fusion, the middle inequality holds because $\Delta$ is a strict refinement of $\Lambda$, and the third equality holds because by Lemma~\ref{L: RefinementPreservesStrictFusion} there is a strict fusion $\Delta\longrightarrow f(\Delta)$. It follows that there is a strict inequality $\excess(f(\Lambda)) > \excess(f(\Delta))$, which means that $f(\Delta)$ has to be a strict refinement of $f(\Lambda)$. 
\end{proof}

The following definition plays a crucial role in the formulation of the main theorem.
\begin{definition}\label{D: TheCategoryE_n}
Let $\calE_n$ be the category of irreducible partitions of excess $n$ and strict fusions between them.
\end{definition}
\begin{example}\label{Ex: TheCategoryE_n}
$\calE_1$ is a groupoid, equivalent to a category with one object: the partition $(1,2)$. The group of automorphisms  of an object of $\calE_1$ is $\Sigma_2$. 

$\calE_2$ is equivalent to a category with two objects: the partitions $(1,2)(3,4)$ and $(1,2,3)$. The monoid of self-morphisms of the first object is the wreath product group $\Sigma_2\wr\Sigma_2$. The monoid of self-morphisms of the second object is the group $\Sigma_3$. There are $4$ morphisms from the first object to the second one, corresponding to the $4$ different ways to glue together two points from the two components of $(1,2)(3,4)$. There are no morphisms from $(1,2,3)$ to $(1,2)(3,4)$.
\end{example}

\begin{remark}\label{R: StructureOfE_n}
For a general $n$, $\calE_n$ can be filtered by full subcategories 
\begin{equation}\label{Eq: Filtration}
\calE_n^1\subset \calE_n^2\subset\cdots\subset \calE_n^{n}=\calE_n
\end{equation}
where $\calE_n^i$ is the full sub-category of $\calE_n$ consisting of partitions whose number of components is between $1$ and $i$. Equivalently, $\calE_n^i$ consists of partitions whose support has  between $n+1$ and $n+i$ elements. For example, the category $\calE_n^1$ is equivalent to the category with one object, corresponding to the indiscrete partition of $\underline{n+1}$, whose group of automorphisms is, naturally, $\Sigma_{n+1}$.
On the other extreme, $\calE_n^{n}\setminus \calE_n^{n-1}$ consists of partitions of type $\underbrace{2-\cdots-2}_n$. The group of automorphisms of an object of this type is the wreath product $\Sigma_2\wr\Sigma_n$. This latter type of partitions corresponds to the ``chord diagrams'' in the Vassiliev spectral sequence for the homology of knot spaces. 

Note in particular that the size of support of objects of $\calE_n$ ranges between $n+1$ and $2n$. It is not hard to see that morphisms in $\calE_n$ are surjective as maps of sets. It follows that there are no morphisms from objects of $\calE_n^{i-1}$ to objects of $\calE_n^i\setminus \calE_n^{i-1}$. Moreover, every morphism between objects of $\calE_n^i\setminus \calE_n^{i-1}$ is an isomorphism. In other words, the category $\calE_n^{\op}$ is nicely filtered, in the sense of Definition~\ref{D: NicelyFiltered}.
\end{remark}

Let $f:\Lambda\longrightarrow \Lambda'$ be a fusion. Clearly, $f$ gives rise to a map of partition posets ${\calP}(\Lambda)\longrightarrow {\calP}(\Lambda')$, and therefore it induces a map of spaces $|{\calP}(\Lambda)|\longrightarrow|{\calP}(\Lambda')|$. If $f$ is a strict fusion, then more is true.
\begin{lemma} \label{L: PartitionPosetsFunctorial}
Let $f:\Lambda\longrightarrow \Lambda'$ be a strict fusion. Then the induced map $f:|{\calP}(\Lambda)|\longrightarrow|{\calP}(\Lambda')|$ restricts to a map
$f:\partial|{\calP}(\Lambda)|\longrightarrow\partial|{\calP}(\Lambda')|$.
Therefore it also induces a map of spaces $T_\Lambda\longrightarrow T_{\Lambda'}$. 
\end{lemma}
\begin{proof}
Recall that $\partial|{\calP}(\Lambda)|$ is defined as the subcomplex of the simplicial nerve of $|{\calP}(\Lambda)|$ consisting of chains of partitions that do not contain both the initial and the final partition. It follows from Lemma~\ref{L: PreserveStrictRefinements} that $f$, being a strict fusion, preserves chains with this property.
\end{proof}
\begin{remark}\label{R: DefinedFunctors}
In view of Lemma~\ref{L: PartitionPosetsFunctorial}, we have constructed a covariant functor from $\calE_n$ to pairs of spaces, given on objects by $$\Lambda\mapsto \left(|{\calP}(\Lambda)|,\partial|{\calP}(\Lambda)|\right)$$
Taking quotients, we obtain a functor from $\calE_n$ to pointed spaces 
$\Lambda\mapsto T_\Lambda$.
\end{remark}
The following variation of the definition of $T_\Lambda$ will be useful to us at one point. Recall that $\calP(\Lambda)$ is the poset of refinements of $\Lambda$. $\calP(\Lambda)$ has an initial and a final object - $\Lambda$ itself, and the discrete partition. Let $\hat 1$ denote the discrete partition. Let $\calP(\Lambda)\setminus \{\Lambda, \hat 1\}$ be the poset obtained from $\calP$ by removing the initial and final objects. Similarly, define the poset $\calP(\Lambda)\setminus \{ \hat 1\}$. Clearly, there is a natural covariant functor from $\calE_n$ to pairs of categories, given on objects by
$$\Lambda\mapsto \left(\calP(\Lambda)\setminus \{ \hat 1\}, \calP(\Lambda)\setminus \{\Lambda, \hat 1\}\right).$$
The proof of the following elementary lemma is left to the reader.
\begin{lemma}\label{L: AlternativeForT}
There is a natural equivalence of functors on $\calE_n$
$$T_\Lambda\simeq S^1\wedge|\calP(\Lambda)\setminus \{ \hat 1\}|/|\calP(\Lambda)\setminus \{\Lambda, \hat 1\}|.$$
\end{lemma}

\subsection{Good and bad partitions.}\label{S: GoodnessBadness} In this subsection, we consider partitions of a fixed set $s$. Let $\Lambda$ be a partition of $s$, represented by a surjection $\alpha:s\twoheadrightarrow \comp(\Lambda)$, and let $\Delta$ be another partition of $s$, represented by a surjection $\beta:s  \twoheadrightarrow \comp(\Delta)$. 
We may form a pushout square, in which all maps are surjective
$$\begin{array}{ccc}
s & \stackrel{\alpha}{\twoheadrightarrow} & \comp(\Lambda) \\
\!\!\!\!\!\beta\downarrow & & \downarrow\\
\comp(\Delta) &\twoheadrightarrow & \comp(\Lambda\wedge\Delta) \end{array}.$$
Here $\Lambda\wedge\Delta$ is the ``finest common coarsening'' of $\Lambda$ and $\Delta$. It is easy to see that the above square is indeed a pushout square. From this square, we can see that $\Lambda$ and $\Delta$ induce a partition of $\comp(\Delta)$, whose set of components is $\comp(\Lambda\wedge\Delta)$. Let us denote this partition by $(\comp(\Delta),\comp(\Lambda\wedge\Delta))$. From the above square, we see that there is a fusion $\Lambda\longrightarrow (\comp(\Delta),\comp(\Lambda\wedge\Delta))$.

\begin{definition}\label{D: GoodnessAndBadness}
We say that {\em $\Delta$ is good relative to $\Lambda$} if the morphism $\Lambda\longrightarrow (\comp(\Delta),\comp(\Lambda\wedge\Delta))$ is a strict fusion. Otherwise, $\Delta$ is {\em bad} relative to $\Lambda$.
\end{definition}
In view of Corollary~\ref{C: SimplyConnected}, we have the following characterization of goodness (and badness).
\begin{lemma}\label{L: CharacterizationOfGoodness}
$\Delta$ is good relative to $\Lambda$ if and only if the homotopy pushout
of the following diagram is simply connected
$$\comp(\Lambda)\twoheadleftarrow s \twoheadrightarrow \comp(\Delta).$$
\end{lemma}
%


We observe that the property of being bad relative to a fixed partition $\Lambda$ is preserved under taking coarsenings.
\begin{lemma}\label{L: BadnessHereditary}
Suppose $\Delta$ is bad relative to $\Lambda$ and $\Delta'\le \Delta$. Then $\Delta'$ is bad relative to $\Lambda$.
\end{lemma}
\begin{proof} The assumption that $\Delta'\le \Delta$ means, in terms of representing surjections, that there is a factorization
$$s\twoheadrightarrow \comp(\Delta)\twoheadrightarrow \comp(\Delta').$$ We have the following diagram, in which maps are surjective, and all squares are pushout squares
$$\begin{array}{ccc}
s & \stackrel{\alpha}{\twoheadrightarrow} & \comp(\Lambda) \\
\downarrow & & \downarrow\\
\comp(\Delta) &\twoheadrightarrow & \comp(\Lambda\wedge\Delta) \\
\downarrow & & \downarrow\\
\comp(\Delta') &\twoheadrightarrow & \comp(\Lambda\wedge\Delta') 
\end{array}.$$
This diagram can be interpreted as a sequence of fusions
$$\Lambda \longrightarrow (\comp(\Delta), \comp(\Lambda\wedge\Delta))\longrightarrow (\comp(\Delta'), \comp(\Lambda\wedge\Delta')).$$
The assumption that $\Delta$ is bad relative to $\Lambda$ is equivalent to saying that the first fusion is not strict. It follows by Lemma~\ref{L: StrictFusionsHereditary} that the composed fusion is not strict, which is equivalent to saying that $\Delta'$ is bad relative to $\Lambda$.
\end{proof}
The proof of the following easy lemma is left to the reader.
\begin{lemma}\label{L: RefinementsAreBad}
Let $\Lambda$ be a partition and let $\Delta$ be a non-discrete refinement of $\Lambda$ then $\Delta$ is bad relative to $\Lambda$.
\end{lemma}
We will also need the following lemma, which says that ``badness is preserved by strict fusions''.
\begin{lemma}\label{L: StrictFusionsPreserveBadness}
Suppose $f:\Lambda\longrightarrow \Lambda'$ is a strict fusion, and $\Delta$ is bad relative to $\Lambda$. Then $\Delta':=f(\Delta)$ is bad relative to $\Lambda'$.
\end{lemma}
\begin{proof}
It is easy to see that we have a square diagram of fusions.
$$\begin{array}{ccc}
\Lambda & \longrightarrow &(\comp(\Delta),\comp(\Lambda\wedge\Delta)) \\
\downarrow & & \downarrow \\
\Lambda' & \longrightarrow& (\comp(\Delta'),\comp(\Lambda'\wedge\Delta'))\end{array}.$$
Since $\Delta$ is bad relative to $\Lambda$, the top horizontal fusion is not strict. It follows that the composed fusion from the upper left to the lower right corner is not strict. By our assumption, the left vertical fusion is strict. It follows that the bottom horizontal fusion is not strict, which means that $\Delta'$ is bad relative to $\Lambda'$.
\end{proof}

\subsection{Non locally constant maps.} Let $M$ be a manifold. Let $(s, \Lambda)$ be a partition. We use the notation $M^\Lambda$ to mean $M^s$. This is the space of maps from $\Lambda$ to $M$. Now suppose that $\Lambda$ is represented by a surjection $s \twoheadrightarrow \comp(\Lambda)$. This induces an inclusion $M^{\comp(\Lambda)}\hookrightarrow M^\Lambda$. Note that the image of this inclusion map does not depend on the choice of surjection $s\twoheadrightarrow \comp(\Lambda)$. We identify $M^{\comp(\Lambda)}$ with its image in $M^\Lambda$. Consider the space
$$M^\Lambda \setminus M^{\comp(\Lambda)}.$$
Thins construction plays an important role in the paper. We refer to the image of $M^{\comp(\Lambda)}$ in $M^\Lambda$ as the space of locally constant maps from $\Lambda$ to $M$, for its points correspond precisely to the maps from $\Lambda$ to $M$ that are constant on each component of $\Lambda$. Accordingly, $M^\Lambda\setminus M^{\comp(\Lambda)}$ is the space of maps from $\Lambda$ to $M$ that are not locally constant. We will some times denote it $\nlc(\Lambda, M)$. We claim that this construction defines a contravariant functor on the category of partitions. More precisely, let $f\colon (s_1, \Lambda_1)\to (s_2, \Lambda_2)$ be a morphism in $\calS$. Clearly, the underlying map of sets gives rise to a map
$f^{\sharp}\colon M^{\Lambda_2}\longrightarrow M^{\Lambda_1}$. The following lemma refers to this notation.
\begin{lemma} \label{L: BasicFunctoriality}
If $f$ is a fusion then 
$$(f^\sharp)^{-1}(M^{\comp(\Lambda_1)})= M^{\comp(\Lambda_2)}$$

For any morphism of partitions $f$, 
$$f^\sharp(M^{\Lambda_2}\setminus M^{\comp(\Lambda_2)})\subset M^{\Lambda_1}\setminus M^{\comp(\Lambda_1)}.$$
\end{lemma}
\begin{proof}
By Lemma \ref{L: Pushout}, $f$ fits into a pushout square
$$\begin{array}{ccc}
s_1 & \stackrel{\alpha_1}{\twoheadrightarrow} & c_1 \\
 f\!\downarrow \mbox{ }& & \downarrow \\
s_2 & \stackrel{\alpha_2}{\twoheadrightarrow} & c_2 \end{array}.$$
Taking functions into $M$, we obtain a pullback square
$$\begin{array}{ccc}
M^{\comp(\Lambda_2)} & \hookrightarrow & M^{\Lambda_2} \\
\downarrow & &\mbox{ } \downarrow { f^\sharp} \\
M^{\comp(\Lambda_1)} & \hookrightarrow & M^{\Lambda_1} \end{array} $$
This being a pullback square is equivalent to the first statement of the lemma. It follows that the second claim is satisfied if $f$ is a fusion. Since every morphism in $\calS$ factors as a fusion followed by a refinement, it remains to prove the second claim in the case when $f$ is a refinement. This follows easily from Lemma~\ref{L: Refinements}.
\end{proof}

\section{Orthogonal calculus - some toy examples.} \label{S: Calculus}
We already explained in the introduction how we would like to apply M. Weiss' orthogonal calculus to the study of covariant isotopy functors on manifolds. Recall that a functor $F(N)$ is polynomial (homogeneous, etc.) in our sense, if the associated functor on vector spaces $V\mapsto F(N\times V)$ is polynomial (homogeneous, etc.) in the sense of Weiss. We now proceed to discuss some examples. The following example is, one could say, the basic building block in our main construction.

Let $\Lambda$ be a non-discrete partition. Let  $h:\supp(\Lambda)\longrightarrow N$ be a point in $N^\Lambda$. Consider the functor
$$N\mapsto \Sigma^\infty \hofiber_h(N^\Lambda\setminus N^{\comp(\Lambda)} \longrightarrow N^\Lambda)$$
The definition of this functor depends on a choice of a basepoint $h\in N^\Lambda$. On the other hand, since we can use a basepoint-free version of $\Sigma^\infty$,  we do not need to choose a basepoint in ${N^\Lambda\setminus N^{\comp(\Lambda)}}$. When $h$ happens to be an element of $N^\Lambda\setminus N^{\comp(\Lambda)}$ we may use any version of $\Sigma^\infty $.

Recall that $\Cyl_\Lambda$ is the mapping cylinder of the surjection $\supp(\Lambda)\twoheadrightarrow \comp(\Lambda)$ representing $\Lambda$. There are natural inclusions (the second of which is a homotopy equivalence) $\supp(\Lambda)\hookrightarrow \Cyl_{\Lambda}$ and  $\comp(\Lambda)\hookrightarrow \Cyl_\Lambda$. 
\begin{definition}
Let $\Omega^\Lambda_h N$ be the space of maps $\alpha\colon\Cyl_\Lambda\longrightarrow N$ which satisfy $\alpha|_{\supp(\Lambda)}=h$. In other words, $\Omega^\Lambda_h N$ is defined by means of a pullback square
$$\begin{array}{ccc}
\Omega^\Lambda_h N & \longrightarrow & N^{\Cyl_\Lambda}\\
\downarrow & & \downarrow \\
* &\stackrel{h}{\longrightarrow} & N^{\supp(\Lambda)} \end{array}$$
\end{definition}
Sometimes, we may omit the subscript $h$ from the notation, and write $\Omega^\Lambda N$ to mean  ``$\Omega^\Lambda_h N$ for some choice of $h$''. To motivate the notation $\Omega^\Lambda N$ notice that if, for example, $N$ is connected then there are equivalences 
$$\Omega^\Lambda N\simeq \map_*(\Cyl_\Lambda/\supp(\Lambda), N)\cong (\Omega N)^{\excess(\Lambda)}$$ where $\Omega N$ is the ordinary loop space of $N$ (for some choice of basepoint). This is because the quotient $\Cyl_\Lambda/\supp(\Lambda)$ is homotopy equivalent to a wedge of $\excess(\Lambda)$ circles. 

The inclusion $\comp(\Lambda)\hookrightarrow \Cyl_\Lambda$ gives rise to the composite map
$\Omega^\Lambda N\longrightarrow N^{\Cyl_\Lambda}\longrightarrow N^{\comp(\Lambda)}.$
Recall that there is a natural vector bundle of dimension $\excess(\Lambda)d$ over $N^{\comp(\Lambda)}$. Namely, the normal bundle of $N^{\comp(\Lambda)}$ in $N^\Lambda$. This pulls back to a vector bundle over $\Omega^\Lambda N$ via the above map. Let $\left(\Omega^\Lambda N\right)_+\tilde\wedge S^{\excess(\Lambda)d}$
be the Thom space of this bundle.

\begin{lemma}\label{L: BasicHomogeneousFunctor}
There is a natural weak equivalence 
$$ \Sigma^\infty \hofiber_h({N^\Lambda\setminus N^{\comp(\Lambda)}} \longrightarrow N^\Lambda)
\simeq \Omega \Sigma^\infty  \left(\Omega^\Lambda_h N\right)_+\tilde\wedge S^{\excess(\Lambda)d}.$$
In particular, the functor $ \Sigma^\infty \hofiber_h({N^\Lambda\setminus N^{\comp(\Lambda)}} \longrightarrow N^\Lambda)$ is homogeneous of degree $\excess(\Lambda)$.
\end{lemma}
\begin{proof}
To analyze the functor $\Sigma^\infty \hofiber(N^\Lambda\setminus N^{\comp(\Lambda)} \longrightarrow N^\Lambda)$, it is convenient to replace $N^\Lambda\setminus N^{\comp(\Lambda)}$ with a suitable compactification. Let us $\overline{N^\Lambda\setminus N^{\comp(\Lambda)}}$ be the spherical blowup of $N^\Lambda$ at $N^{\comp(\Lambda)}$. It is a manifold with boundary, obtained by replacing $N^{\comp(\Lambda)}$ with the sphere bundle of its normal bundle in $N^\Lambda$. 
Then the map $N^\Lambda\setminus N^{\comp(\Lambda)} \longrightarrow N^\Lambda$ factors as
$$N^\Lambda\setminus N^{\comp(\Lambda)} \longrightarrow \overline{N^\Lambda\setminus N^{\comp(\Lambda)}} \longrightarrow N^\Lambda$$
where the first map is a weak equivalence and the second map is the blow-up map. Therefore, we may analyze $\Sigma^\infty \hofiber(\overline{N^\Lambda\setminus N^{\comp(\Lambda)}} \longrightarrow N^\Lambda)$
instead of $\Sigma^\infty \hofiber({N^\Lambda\setminus N^{\comp(\Lambda)}} \longrightarrow N^\Lambda)$.

There is a commutative square, which is both a pushout and a homotopy pushout.
$$
\begin{CD}
\partial\left(\overline{N^\Lambda\setminus N^{\comp(\Lambda)}}\right) @>>> \overline{N^\Lambda\setminus N^{\comp(\Lambda)}} \\
@VVV @VVV \\
N^{\comp(\Lambda)} @>>> N^\Lambda 
\end{CD}.
$$
Moreover, $\partial\left(\overline{N^\Lambda\setminus N^{\comp(\Lambda)}}\right)$ is homeomorphic to $N^{\comp(\Lambda)}\tilde\times S^{\excess(\Lambda)d-1}$, the sphere bundle of the normal bundle of $N^{\comp(\Lambda)}$ in $N^\Lambda$. 
The square diagram above is in fact a diagram of spaces over $N^\Lambda$. It is well-known that taking homotopy fibers of maps into a fixed space commutes with homotopy pushouts. Therefore, we have a homotopy pushout square
$$\begin{CD}
\hofiber(N^{\comp(\Lambda)}\tilde\times S^{\excess(\Lambda)d-1}\longrightarrow N^\Lambda) 
@>>> 
\hofiber(\overline{N^\Lambda\setminus N^{\comp(\Lambda)}}\longrightarrow N^\Lambda) \\
@VVV @VVV \\
\hofiber(N^{\comp(\Lambda)}\longrightarrow N^\Lambda) @>>> * \end{CD}.
$$
For the moment, let $X$ be the homotopy cofiber of any of the vertical maps in the square. Clearly, $X$ is a space with a basepoint, and there is a natural equivalence 
$$\Sigma^\infty \hofiber(\overline{N^\Lambda\setminus N^{\comp(\Lambda)}}\longrightarrow N^\Lambda) \simeq \Omega \Sigma^\infty X.$$
It remains to show that $X\simeq \left(\Omega^\Lambda N\right)_+\tilde\wedge S^{\excess(\Lambda)d}$.
Let us begin by analyzing the homotopy fiber of the map $N^{\comp(\Lambda)}\longrightarrow N^\Lambda=N^{\supp(\Lambda)}$. Since the map $\supp(\Lambda)\twoheadrightarrow \comp(\Lambda)$ factors as 
$$\supp(\Lambda)\hookrightarrow \Cyl_\Lambda \stackrel{\simeq}{\longrightarrow} \comp(\Lambda)$$
where the first map is a cofibration and the second map is a homotopy equivalence, it follows that
the map $N^{\comp(\Lambda)}\longrightarrow N^\Lambda$ factors as 
$$N^{\comp(\Lambda)}\stackrel{\simeq}{\longrightarrow} N^{\Cyl_\Lambda} \longrightarrow N^\Lambda ,$$ where the first map is a homotopy equivalence and the second map is a fibration. Therefore, the homotopy fiber of the map $N^{\comp(\Lambda)}\longrightarrow N^\Lambda=N^{\supp(\Lambda)}$
can be identified with the strict fiber of the map $N^{\Cyl_\Lambda} \longrightarrow N^\Lambda$ (over the basepoint map $h:\supp(\Lambda)\longrightarrow N$). This space is exactly $\Omega^\Lambda_h N$. Similarly, the space $\hofiber(N^{\comp(\Lambda)}\tilde\times S^{\excess(\Lambda)d-1}\longrightarrow N^\Lambda) $
can be identified as $(\Omega^\Lambda_h N)\tilde\times S^{\excess(\Lambda)d-1}$, which is the pullback of the sphere bundle over $N^{\comp(\Lambda)}$ along the map
$\Omega^\Lambda N \longrightarrow N^{\comp(\Lambda)}$. The map $\hofiber(N^{\comp(\Lambda)}\tilde\times S^{\excess(\Lambda)d-1}\rightarrow N^\Lambda) \longrightarrow \hofiber(N^{\comp(\Lambda)}\rightarrow N^\Lambda) $ can now be identified with a sphere bundle  projection. The homotopy cofiber of a sphere bundle projection is the Thom space of the underlying vector bundle. Which is what we wanted to prove.

To see that this functor is indeed homogeneous of degree $\excess(\Lambda)$, note that it follows that there is a natural equivalence
$$ \Sigma^\infty \hofiber_h({(N\times V)^\Lambda\setminus (N\times V)^{\comp(\Lambda)}} \longrightarrow (N\times V)^\Lambda)
\simeq \Omega \Sigma^\infty  \left(\Omega^\Lambda_h N\right)_+\tilde\wedge S^{\excess(\Lambda)d}\wedge S^{\excess(\Lambda) V}.$$
\end{proof}
 
\begin{corollary}\label{C: ParallelizableCase}
If $N$ is stably parallelizable, then 
$$ \Sigma^\infty \hofiber_h({N^\Lambda\setminus N^{\comp(\Lambda)}} \longrightarrow N^\Lambda)
\simeq \Omega \Sigma^\infty  \left(\Omega^\Lambda_h N\right)_+\wedge S^{\excess(\Lambda)d}.$$
I.e., ``Thom space'' may be replaced with smash product.
\end{corollary}

Note that the homotopy type of $\Sigma^\infty\hofiber_h(N^\Lambda\setminus N^{\comp(\Lambda)}\longrightarrow N^\Lambda)$ is determined by excess of $\Lambda$ and the homotopy type of the basepoint map $h\colon \Lambda\longrightarrow N$. These are preserved by strict fusions. Indeed, it is not difficult to adapt the proof of Lemma~\ref{L: BasicHomogeneousFunctor} to prove (using Proposition~\ref{P: HomotopyPushout}) the following lemma.
\begin{lemma}\label{L: AlwaysEquivalence}
Let $\alpha\colon \Lambda\longrightarrow \Lambda_1$ be a strict fusion. Let $h\in N^{\Lambda_1}$. $\alpha$ induces a homotopy equivalence
$$\Sigma^\infty\hofiber_h(N^{\Lambda_1}\setminus N^{\comp(\Lambda_1)}\rightarrow N^{\Lambda_1}) \longrightarrow \Sigma^\infty\hofiber_{h\circ\alpha}(N^\Lambda\setminus N^{\comp(\Lambda)}\rightarrow N^\Lambda).$$
\end{lemma}

\subsection{Configuration Spaces.} 
Let $\config(k,N)=N^k\setminus \Delta^kN$ be the configuration space of ordered $k$-tuples of distinct points in $N$, or, in other words, the complement of the fat diagonal (which we denote $\Delta^kN$) in $N^k$. We would like to describe the homogeneous layers in the orthogonal tower of the functor $N\mapsto\Sigma^\infty\config(k,N)$, and also the functor $\Sigma^\infty \overline{\config}(k,N)$, where $\overline{\config}(k,N)=\hofiber(\config(k,N)\longrightarrow N^k)$. The definition of $\overline{\config}(k,N)$ depends on a choice of a point in $N^k$, which we will suppress from the notation except when we really have to mention it explicitly.

Recall that $\calP_k$ is the poset of partitions of $k$. In this section we will need to use a modified version of the poset.
\begin{definition}
Let $\calMP_k$ be the poset of non-discrete partitions of $k$, together with an adjoint initial object, denoted $\Omega$.
\end{definition}
\begin{example}
Here is a picture of the poset $\calMP_3$.
$$
\xymatrix{
& & (13)(2) \\
\Omega\ar[r] & (123) \ar[r]\ar[ru] \ar[rd] & (12)(3)\\
& & (23)(1) 
}
$$
\end{example}

There is an evident contravariant functor (by Lemma~\ref{L: BasicFunctoriality}) from $\calMP_k$ to spaces given on objects as follows  $$\Lambda\mapsto \left\{ \begin{array}{cl}
N^k\setminus N^{\comp(\Lambda)} &\mbox{If }\Lambda\ne \Omega \\
N^k & \mbox{if } \Lambda=\Omega \end{array}\right.$$
We will denote this functor simply by $N^k\setminus N^{\comp(\Lambda)}$, it being understood that if $\Lambda=\Omega$ then $N^{\comp(\Lambda)}=\emptyset$.
Applying $\Sigma^\infty$, we obtain the functor $\Lambda\mapsto \Sigma^\infty N^k\setminus N^{\comp(\Lambda)}$. 
\begin{example}
Here is a picture of this functor in the case $k=3$.
$$
\xymatrix{
& & \ar[ld]\Sigma^\infty N^3\setminus N^2 \\
 \Sigma^\infty N^3 &  \Sigma^\infty N^3\setminus N \ar[l] &  \ar[l]\Sigma^\infty N^3\setminus N^2\\
& &  \ar[lu]\Sigma^\infty N^3\setminus N^2 
}
$$
\end{example}
Clearly, $\Sigma^\infty \config(k,N)$ maps into this diagram, and so it maps into the homotopy inverse limit of this diagram. Similarly, there is an evident contravariant functor given on objects by $\Lambda \mapsto \Sigma^\infty \hofiber(N^k\setminus N^{\comp(\Lambda)}\longrightarrow N^k)$, and $\Sigma^\infty \overline{\config}(k,N)$ maps into the homotopy inverse limit of this functor.
\begin{proposition}\label{P: ModelForStableConfig}
The natural maps
$$\Sigma^\infty \config(k,N) \longrightarrow \underset{\Lambda\in\calMP_k}{\holim  }\, \Sigma^\infty N^k\setminus N^{\comp(\Lambda)}$$
and
$$\Sigma^\infty \overline{\config}(k,N) \longrightarrow \underset{\Lambda\in\calMP_k}{\holim  }\, \Sigma^\infty \hofiber\, (N^k\setminus N^{\comp(\Lambda)}\longrightarrow N^k)$$
are homotopy equivalences.
\end{proposition}
\begin{proof}
To begin with, we may ignore the initial object $\Omega$, for it is easy to see that the subposet $\calMP_k\setminus\{\Omega\}$ is final in $\calMP_k$ and therefore for any contravariant functor on $\calMP_k$, the restriction to $\calMP_k\setminus\{\Omega\}$ induces an equivalence on homotopy limits. So, for the duration of this proof only, we will use $\calMP_k$ to mean $\calMP_k\setminus \{\Omega\}$.
\begin{remark} The reader may wonder why we needed to introduce the extra object $\Omega$ to begin with. Indeed, it plays only a minor role, but it makes the construction more natural, and some things work better with it. More specifically, the object $\Omega$ plays a role in the proof of Lemma~\ref{L: NeedsOmega} below. 
\end{remark}
We will analyze functors on $\calMP_k$ by relating them to certain cubical diagrams, which we will now describe. For the duration of this subsection only, let $\calS$ be the poset of subsets of $k\choose 2$, where $k\choose 2$ is the set of unordered pairs of distinct elements of $k$. A functor on $\calS$ is a cubical diagram of dimension $k \choose 2$ (we often use the same symbol to denote both a finite set and the number of elements in it).
 
The poset $\calS$ can be identified with the poset of graphs with vertex set $k$. There is an order-reversing map of posets from $\calS$ to the poset of all partitions of $k$, which sends a subset $U\subset {k\choose 2}$ to the partition of $k$ given by the connected components of the graph corresponding to $U$. Let us denote the partition associated with $U$ by $\Lambda(U)$. 

Let $\calS\sp{1}$ be the poset of non-empty subsets of $k\choose 2$. Clearly, the map of posets $U\mapsto \Lambda(U)$ restricts to an (order reversing) map from $\calS\sp{1}$ to $\calMP_k$, the poset of non-trivial partitions of $k$. Let $F(-)$ be a contravariant functor from $\calMP_k$ to spaces or spectra. Then $F(\Lambda(-))$ is a covariant functor on $\calS\sp{1}$. There is a natural map
$$\underset{\Lambda\in \calMP_k}{\holim} \, F(\Lambda) \longrightarrow \underset{U\in \calS\sp{1}}{\holim} \,F(\Lambda(U)).$$
\begin{lemma}\label{L: EquivalentLimits}
The above map of homotopy limits is a weak equivalence for any contravariant functor $F$.
\end{lemma}
\begin{proof}
It is well-known that to prove the lemma it is enough to show that for every object $\Lambda$ of $\calMP_k$, the nerve of the over category $\Lambda\downarrow\calS\sp{1}$, consisting of all sets $U$ such that $\emptyset \ne U\subset {k\choose 2}$, and such that $\Lambda(U)$ is a refinement of $\Lambda$, is contractible. But it is easy to see that $\Lambda\downarrow\calS\sp{1}$ has a minimal object - the set of all pairs of elements of $k$ that are contained in some component of $\Lambda$ - and so its nerve is contractible.
\end{proof}

Now define a cubical diagram $\chi\colon \calS \longrightarrow \Top$ as follows.
$$
\chi(U)=\left\{
\begin{array}{ll}
N^k\setminus N^{\comp(\Lambda(U))} & \mbox{If } U\neq \emptyset \\
\config(k,N) & \mbox{If } U=\emptyset \end{array}\right..
$$
the maps in the diagram being the obvious inclusions.
The following lemma summarizes some elementary properties of the cubical diagram $\chi$. Its proof is left as an easy exercise to the reader.
\begin{lemma}\label{L: PropertiesOfChi}
\begin{enumerate}
\item For all $U\subset {k\choose 2}$, $\chi(U)$ is an open subset of $\chi\left({k\choose 2}\right)=N^k\setminus N$.
\item For every $\emptyset \neq U\subset {k\choose 2}$, 
$$\chi(U)=\bigcup_{x\in U} \chi\left(\{x\}\right),$$
where the union is taken in $\chi\left({k\choose 2}\right)$.
\item $\chi(\emptyset)=\underset{x\in {k \choose 2}}{\bigcap} \chi\left(\{x\}\right)$.
\end{enumerate}
\end{lemma}
The key fact about the cubical diagram $\chi$ is given by the following lemma.
\begin{lemma} \label{L: ModelForConfig}
$\chi$ is a homotopy pushout cube.
\end{lemma}
Lemma~\ref{L: ModelForConfig} follows immediately from Lemma~\ref{L: PropertiesOfChi} and Lemma~\ref{L: MoreGeneral} below.
\begin{lemma}\label{L: MoreGeneral}
Let $X$ be a topological space. Let $X_1,\ldots, X_k$ be open subsets of $X$ such that $X=\cup_{i=1}^k X_i$. Let $X_0=\cap_{i=1}^k X_i$. Define a $k$-dimensional cubical diagram $\chi$ by $\chi(\emptyset)=X_0$, and for $\emptyset \neq U \subset \{1,\ldots, k\}$ $$\chi(U)=\bigcup_{i\in U} X_i,$$
where the union is taken in $X$. Then $\chi$ is a homotopy pushout cube.
\end{lemma}
\begin{proof}
We prove by induction on $k$ beginning with $k=2$. In this case $\chi$ is a square diagram
$$\begin{array}{ccc}
X_0 & \longrightarrow & X_1 \\
\downarrow & & \downarrow \\
X_2 & \longrightarrow & X \end{array}$$
where all maps are inclusions of open subsets, $X_0=X_1\cap X_2$, and $X= X_1\cup X_2$. So $\chi$ is obviously a homotopy pushout square. Let us assume that the lemma holds for $k-1$, and let $\chi$ be $k$-dimensional. We can present $\chi$ as a map between $k-1$-dimensional cubical diagrams $\chi_1\longrightarrow \chi_2$, where $\chi_1$ is the subcube of $\chi$ indexed by all subsets of $\{1,\ldots,k\}$ that do not contain $k$, while $\chi_2$ is indexed by subsets that do contain $k$. Then $\chi$ is a homotopy pushout if and only if the map $\chi_1\longrightarrow \chi_2$ induces a homotopy equivalence of total homotopy cofibers. The cubes $\chi_1$ and $\chi_2$ almost satisfy the assumptions of the lemma, except the initial objects of these cubes do not equal the intersection of the rank $1$ objects. Let $\chi_1'$ and $\chi_2'$ be the cubes obtained from $\chi_1$ and $\chi_2$ respectively by replacing the initial objects with the intersection of rank $1$ objects. It is easy to see that there is a naturally defined square diagram of $k-1$-dimensional cubes
$$\begin{array}{ccc}
\chi_1 & \longrightarrow & \chi_2 \\
\downarrow & & \downarrow \\
\chi_1' & \longrightarrow & \chi_2' \end{array} $$
The cubes $\chi_1'$ and $\chi_2'$ are $k-1$-dimensional cubes that satisfy the assumptions of the lemma, so they are homotopy pushouts by the induction hypothesis. If we can show that the square diagram above induces a homotopy pushout square of total homotopy cofibers, it will follow that the map $\chi_1\longrightarrow \chi_2$ induces a homotopy equivalence of total cofibers. By definition, the vertical maps in the square of cubes, namely the maps $\chi_1(U)\longrightarrow \chi_1'(U)$ and $\chi_2(U)\longrightarrow \chi_2'(U)$ are equivalences for all $U$ except $U=\emptyset$, so the induced square diagrams of spaces are homotopy pushouts except possibly at the initial objects. So, we need to show the square diagram of spaces
$$\begin{array}{ccc}
\chi_1(\emptyset) & \longrightarrow & \chi_2(\emptyset) \\
\downarrow & & \downarrow \\
\chi_1'(\emptyset) & \longrightarrow & \chi_2'(\emptyset) \end{array} $$
is a homotopy pushout square. By definition, this square is equivalent to 
$$\begin{array}{ccc}
X_0 & \longrightarrow & X_k \\
\downarrow & & \downarrow \\
\bigcap_{i=1}^{k-1} X_i & \longrightarrow & \bigcap_{i=1}^{k-1} X_k\cup X_i
\end{array}. $$
Clearly, $\bigcap_{i=1}^{k-1} X_k\cup X_i=X_k\cup \left(\bigcap_{i=1}^{k-1} X_i\right)$, and 
$X_0 = X_k\cap \bigcap_{i=1}^{k-1} X_i$, so this is a homotopy pushout square.
\end{proof}

By Lemma~\ref{L: ModelForConfig}, the cubical diagram $\chi$ is a homotopy pushout. It follows immediately that the cubical diagram of spectra $\Sigma^\infty \chi$, which is defined by 
$$
\Sigma^\infty\chi(U)=\left\{
\begin{array}{ll}
\Sigma^\infty N^k\setminus N^{\comp(\Lambda(U))} & \mbox{If } U\neq \emptyset \\
\Sigma^\infty \config(k,N) & \mbox{If } U=\emptyset \end{array}\right..
$$
is a homotopy pushout cubical diagram. But, in the category of spectra, homotopy pushout cubes are also homotopy pullback cubes, and so $\Sigma^\infty\chi(U)$ is a homotopy pullback cube. In other words, the natural map
$$\Sigma^\infty\config(k,N) \longrightarrow \underset{U\in\calS\sp{1}}{\holim}\,\Sigma^\infty N^k\setminus N^{\comp(\Lambda(U))}$$
is a homotopy equivalence. The first part of Proposition~\ref{P: ModelForStableConfig} now follows from Lemma~\ref{L: EquivalentLimits}.

For the second part of the proposition, recall that $\chi$ is in fact a cube of subspaces of $N^k$, and let $\hofiber(\chi\longrightarrow N^k)$ be the cube obtained by taking, objectwise in $\chi$, the homotopy fiber of the inclusion map into $N^k$. It is well-known that taking homotopy fibers over a fixed space commutes with taking homotopy pushouts, and so $\hofiber(\chi\longrightarrow N^k)$ is a homotopy pushout cube. The proof of the proposition now proceeds in the same way as in the first part. 
\end{proof}

As a corollary, we obtain the following lemma, describing the Taylor polynomials of the functor $\Sigma^\infty \overline{\config}(k,N)$
\begin{lemma}\label{L: TaylorReducedConfigurations}
Let $\calMP_k^n$ be the sub-poset of $\calMP_k$ consisting of partitions of excess $\le n$. Then
$$\TP_n\Sigma^\infty \overline{\config}(k,N) \longrightarrow \underset{\Lambda\in\calMP_k^n}{\holim  }\, \Sigma^\infty \hofiber(N^k\setminus N^{\comp(\Lambda)}\longrightarrow N^k).$$
\end{lemma}
\begin{proof}
By the second part of Proposition~\ref{P: ModelForStableConfig}, 
$$\TP_n\Sigma^\infty \overline{\config}(k,N) \simeq \TP_n\left(\underset{\Lambda\in\calMP_k}{\holim  }\, \Sigma^\infty \hofiber(N^k\setminus N^{\comp(\Lambda)}\longrightarrow N^k)\right).$$
The operator $\TP_n$ is constructed using compact homotopy limits and filtered homotopy colimits. Therefore, it commutes (up to equivalence) with compact homotopy limits. Clearly, the category $\calMP_k$ is compact (indeed, finite). Therefore,
$$\TP_n\Sigma^\infty \overline{\config}(k,N) \simeq \underset{\Lambda\in\calMP_k}{\holim  }\, \TP_n\Sigma^\infty \hofiber(N^k\setminus N^{\comp(\Lambda)}\longrightarrow N^k).$$
By Lemma~\ref{L: BasicHomogeneousFunctor}, the functor $\Sigma^\infty \hofiber(N^k\setminus N^{\comp(\Lambda)}\longrightarrow N^k)$ is homogeneous of degree $\excess(\Lambda)$ (in case $\Lambda=\Omega$, this functor is equivalent to the constant one point functor). Therefore,
$$ 
\TP_n\Sigma^\infty \hofiber(N^k\setminus N^{\comp(\Lambda)}\longrightarrow N^k)\simeq
\left\{\begin{array}{ll}
\Sigma^\infty \hofiber(N^k\setminus N^{\comp(\Lambda)}\longrightarrow N^k) & \mbox{If } \Lambda\in\calMP^n_k \\
* & \mbox{Otherwise}. \end{array}\right. $$
It is easy to check that this functor from $\calMP_k$ to the category of spectra has the property that it is equivalent to the homotopy right Kan extension of its restriction to $\calMP_k^n$. It follows that the restriction of this functor from $\calMP_k$ to $\calMP_k^n$ induces an equivalence of homotopy limits. The lemma follows. 
\end{proof}

\section{Spaces of standard embeddings.} \label{S: StandardEmbeddings}
Let $W$ be a Euclidean space. Recall that an open submanifold $N$ of $W$ is called nice, if it is the interior of a closed submanifold with boundary $\overline N$, where $\overline N$ has the property that there is an $\epsilon>0$ such that the space of all points in $\overline N$ lying within distance less than $\epsilon$ of the boundary of $\overline N$ deformation retracts onto the boundary. 

For the rest of the section, we will assume that $N$ is a nice submanifold of $W$, and $M$ is a disjoint union of $k$ open balls in $W$. We will also use the letter $k$ to denote the set of path components of $M$. Let $q\colon M \to k$ be the quotient map. Recall from the introduction that $\sEmb(M,N)$ is the space of standard embeddings of $M$ into $N$, where an embedding is called standard if on each component of $M$ it differs by a translation from the inclusion. For every choice of section $k\to M$ of $q$ there is an obvious evaluation map $\sEmb(M,N)\longrightarrow \config(k, N)$. It is easy to show that if the balls making up $M$ are small enough, then any such evaluation map is a homotopy equivalence (It is for this step that one needs the assumption that $N$ is nicely embedded). We will generally assume that this condition holds, i.e., that the components of $M$ are small enough to guarantee that this evaluation map, and other similar ones that we will encounter, are equivalences. 

Let $\sEbar(M,N)$ be the homotopy fiber of the inclusion map $\sEmb(M,N)\longrightarrow \sImm(M,N)$ (as usual, this space depends on a choice of a basepoint in $\sImm(M,N)$). Clearly, there is a natural map $\sEbar(M,N)\longrightarrow \overline{\config}(k,N)$ that is a homotopy equivalence (assuming that the components of $M$ are small enough). The constructions $\sEmb(M,N)$ and $\sEbar(M,N)$ are contravariantly functorial with respect to standard embeddings of the $M$ variable. Our task in this section is to extend the description of the Taylor polynomials of $\Sigma^\infty\config(k, N)$ and $\Sigma^\infty\overline{\config}(k, N)$ to a description for $\Sigma^\infty\sEmb(M,N)$ and $\Sigma^\infty\sEbar(M,N)$ that is natural with respect to the $M$ variable.

\begin{definition}
Let $\calS_k$ be the category with the following objects 
\begin{itemize}
\item pairs $(\Lambda, \alpha)$ where $\Lambda$ is an object of $\calS$ (i.e., a partition), and $\alpha\colon \Lambda\to k$ is a map of sets from the support of $\Lambda$ to $k$. 
\item an additional initial object, denoted $\Omega$.
\end{itemize}
A morphism $(\Lambda, \alpha)\longrightarrow (\Lambda', \alpha')$ in $\calS_k$ consists of a morphism $\Lambda\to \Lambda'$ in $\calS$ such that the evident triangle commutes. In addition, there is a unique morphism from $\Omega$ to any other object of $\calS_k$
\end{definition}
We will refer to objects and morphisms in $\calS_k$ as partitions over $k$ and morphism over $k$ respectively. Next we are going to define a contravariant functor from $\calS_k$ to spaces. Recall from Lemma~\ref{L: BasicFunctoriality} that we have a contravariant functor from $\calS$ to spaces
$$\Lambda\mapsto M^\Lambda\setminus M^{\comp(\Lambda)}.$$
Let us denote this construction by $\nlc(\Lambda, M)$ (``nlc'' stands for ``non locally constant''). For an object $(\Lambda, \alpha)$ of $\calS_k$, define
$\nlc_\alpha(\Lambda, M)\subset \nlc(\Lambda, M)$ to be the space of non locally constant maps $$\tilde \alpha\colon \Lambda\longrightarrow M$$
such that $\alpha=q\circ\tilde\alpha$. We will call elements of $\nlc_\alpha(\Lambda, M)$ ``lifts'' of $\alpha$. By our convention, if $\Lambda$ is the initial object $\Omega$, then $\nlc_\alpha(\Lambda, M)$ is the one point space. The construction
$$(\Lambda,\alpha)\mapsto \nlc_\alpha(\Lambda, M)$$ is indeed a contravariant functor, because if $\Lambda\to \Lambda'$ is a morphism in $\calS$ over $k$, and we have a non locally constant morphism $\Lambda' \to M$ then the resulting morphism $\Lambda \to \Lambda' \to M$ is not locally constant (Lemma~\ref{L: BasicFunctoriality}).

\begin{remark}\label{R: PropertiesOfNlc}
If $\Lambda$ is a discrete partition, there are no non-locally constant maps from $\Lambda$ to $M$, so we only need to consider partitions that are not discrete in this context. Also, let us note that  
since each component of $M$ is, by assumption, an open ball in a fixed vector space $W$, it follows that $\nlc_\alpha(\Lambda, M)$ can be identified with a subset of $W^{\supp(\Lambda)}$. If $\alpha\colon \Lambda\to k$ is itself not locally constant then every map $\tilde\alpha\colon \Lambda \to M$ lifting $\Lambda$ is non-locally constant, and it follows that in this case $\nlc_\alpha(\Lambda, M)$ is homeomorphic to a  ball. More precisely, it can be identified with a product of components of $M$, where the component associated with $i\in k$ occurs $\alpha^{-1}(i)$ times. It also follows that if $(\Lambda', \alpha')$ is a partition over $k$ where the map $\alpha'\colon \Lambda' \to k$ is not locally constant, and $\Lambda\le \Lambda'$ (i.e., $\Lambda'$ is a refinement of $\Lambda$) then the induced map 
$$\nlc_{\alpha'}(\Lambda', M) \longrightarrow \nlc_\alpha(\Lambda, M)$$
is a homeomorphism.

On the other hand, if $\alpha$ is locally constant then $\nlc_\alpha(\Lambda, M)$ is homeomorphic to the compement of an equatorial disk of dimension $\comp(\Lambda)d$ in a disk of dimension $\supp(\Lambda)d$, and thus is homotopy equivalent to a sphere of dimension $\excess(\Lambda)d-1$.
\end{remark}

As usual, we can associate with this functor a topological category, denoted $\calS_k\ltimes \nlc_{(-)}(-, M)$. This category has an initial object, which we continue denoting by $\Omega$, and every other object consists of an ordered pair  $\left((\Lambda, \alpha),\tilde\alpha\right)$ where $(\Lambda,\alpha)$ is a partition over $k$, and $\tilde\alpha\colon \Lambda \to M$ is a non-locally constant lift of $\alpha$. A morphism $\left((\Lambda, \alpha),\tilde\alpha\right)\longrightarrow \left((\Lambda', \alpha'),\tilde\alpha'\right)$ consists of a morphism $h\colon\Lambda\to \Lambda'$ over $k$ such that $\tilde\alpha=\tilde\alpha'\circ h$. The space of objects is topologized as the disjoint union of spaces of the form $\nlc_\alpha(\Lambda, M)$. The space of morphisms is topologized accordingly.
\begin{remark}
Of course, $\tilde\alpha$ determines $\alpha$, so our notation has a built-in redundancy, but will be convenient to have it for a while. Essentally, $\calS_k\ltimes \nlc_{(-)}(-, M)$ is the (topological) category of partitions of $M$, with an initial object added. 
\end{remark}

Next we define a contravariant functor from $\calS_k\ltimes\nlc_{(-)}(-, M)$ to spaces. To begin with, the functor sends $\Omega$ to $\sImm(M,N)$, the space of standard immersions of $M$ into $N$. For any other object $\left((\Lambda, \alpha),\tilde\alpha\right)$ of $\calS_k\ltimes\nlc_{(-)}(-, M)$, let $\sImm(M, N)\setminus_{\tilde\alpha} N^{\comp(\Lambda)}$ be the space of standard immersions from $M$ to $N$ whose composition with $\tilde\alpha$ yields a non locally constant map from $\Lambda$ to $N$. The union of all spaces $\sImm(M, N)\setminus_{\tilde\alpha} N^{\comp(\Lambda)}$, where $(\Lambda, \alpha)$ is fixed and $\tilde\alpha$ ranges over all lifts of $\alpha$, is topologized as a subspace of $\nlc_\alpha(\Lambda, M)\times \sImm(M,N)$. Let us call this union $E_\alpha$. Thus there is a map $E_\alpha\to \nlc_\alpha(\Lambda, M)$, and the inverse image of $\tilde\alpha$ under this map is $\sImm(M, N)\setminus_{\tilde\alpha} N^{\comp(\Lambda)}$.
It is easy to check, using Lemma~\ref{L: BasicFunctoriality} again, that this construction defines a contravariant functor from $\calS_k\ltimes\nlc(-, M)$ to spaces. We will denote this functor simply by $\sImm(M, N)\setminus_{\tilde\alpha} N^{\comp(\Lambda)}$, it being understood that the value of the functor at $\Omega$ is $\sImm(M, N)$.
Similarly, one can define a contravariant functor from $\calS_k\ltimes\nlc_{(-)}(-, M)$ to spaces by the formula (which depends on a choice of a point in $\sImm(M,N)$)
$$\overline{\sImm(M, N)\setminus_{\tilde\alpha} N^{\comp(\Lambda)}}:=\hofiber\left(\left(\sImm(M, N)\setminus_{\tilde\alpha} N^{\comp(\Lambda)}\right)\longrightarrow \sImm(M, N)\right).$$

Our goal is to analyze the homotopy limits of the functors $\Sigma^\infty\sImm(M, N)\setminus_{\tilde\alpha} N^{\comp(\Lambda)}$ and ultimately $\Sigma^\infty\overline{\sImm(M, N)\setminus_{\tilde\alpha} N^{\comp(\Lambda)}}.$ 
Clearly, these homotopy limits are spaces of twisted natural transformations between contravariant functors on $\calS_k$, and therefore can be expressed as homotopy limits over the twisted arrow category of $\calS_k$ (see section~\ref{S: TwistedNat}). Let $^a\calS_k$ be the twisted arrow category of $\calS_k$. Let us recall the definition of this category (we will tweak it slightly to make it suitable for contravariant functors). An object of this category is an arrow of the form $\Omega\longrightarrow \Omega$, or $\Omega\longrightarrow (\Lambda, \alpha)$, or $(\Lambda,\alpha)\longrightarrow (\Lambda',\alpha')$. A morphism, say between the following objects
$$\left((\Lambda,\alpha)\longrightarrow (\Lambda',\alpha')\right)\longrightarrow \left((\Delta,\beta)\longrightarrow (\Delta',\beta')\right)$$ in $^a\calS_k$, consists of a commutative diagram in $\calS_k$ (note the customary twist)
$$
\begin{CD}
(\Lambda, \alpha) @<<< (\Delta, \beta) \\
@VVV @VVV \\
(\Lambda', \alpha') @>>> (\Delta', \beta').
\end{CD}$$
In what follows, we will usually spell out the details for objects of $^a\calS_k$ that do not involve $\Omega$, and leave it to the reader to supplement the missing details. 
For an object $(\Lambda,\alpha)\longrightarrow (\Lambda',\alpha')$ let us denote by $E_{\alpha\to \alpha'}$ the pullback of the corresponding diagram
$$\nlc_{\alpha'}(\Lambda', M)\longrightarrow \nlc_{\alpha}(\Lambda, M) \longleftarrow E_\alpha.$$
This is a space over $\nlc_{\alpha'}(\Lambda', M)$. Let $$\Gamma\left(\nlc_{\alpha'}(\Lambda', M); \Sigma^\infty \sImm(M, N)\setminus_{\tilde\alpha} N^{\comp(\Lambda)}\right)$$ be the spectrum of sections of the corresponding fiberwise suspension spectrum over $\nlc_{\alpha'}(\Lambda', M)$. This construction defines a covariant functor from $^a\calS_k$ to spectra, and there is a natural equivalence
$$\underset{\calS_k\ltimes\nlc_{(-)}(-, M)}{\holim} \Sigma^\infty\sImm(M, N)\setminus_{\tilde\alpha} N^{\comp(\Lambda)} \simeq \underset{^a\calS_k}{\holim} \Gamma\left(\nlc_{\alpha'}(\Lambda', M); \Sigma^\infty \sImm(M, N)\setminus_{\tilde\alpha} N^{\comp(\Lambda)}\right).$$
We need to analyze the spectrum of sections $\Gamma\left(\nlc_{\alpha'}(\Lambda', M); \Sigma^\infty \sImm(M, N)\setminus_{\tilde\alpha} N^{\comp(\Lambda)}\right)$ (Proposition~\ref{P: DescriptionSections} below) This spectrum depends on a morphism $f\colon (\Lambda, \alpha)\to (\Lambda', \alpha')$. Factor $f$ as a fusion followed by refinement in the usual way $$\Lambda \longrightarrow \Lambda'' \longrightarrow \Lambda'$$ where $\Lambda''=f(\Lambda)$ and we have suppressed that these are morphisms of partitions over $k$. Let $\alpha''\colon \Lambda'' \longrightarrow \Lambda' \stackrel{\alpha'}{\longrightarrow} k$ be the evident composed map. Let us assume first that $\alpha'\colon \Lambda'\to k$ is not locally constant. In this case the map $$\nlc_{\alpha'}(\Lambda', M) \longrightarrow \nlc_{\alpha''}(\Lambda'', M)$$ is a homeomorphism (remark~\ref{R: PropertiesOfNlc}) and thus our spectrum $$\Gamma\left(\nlc_{\alpha'}(\Lambda', M); \Sigma^\infty \sImm(M, N)\setminus_{\tilde\alpha} N^{\comp(\Lambda)}\right)$$ is equivalent to 
 $$\Gamma\left(\nlc_{\alpha''}(\Lambda'', M); \Sigma^\infty \sImm(M, N)\setminus_{\tilde\alpha} N^{\comp(\Lambda)}\right).$$
On the other hand, since the morphism $\Lambda\to \Lambda''$ is a fusion it follows (from the first part of Lemma~\ref{L: BasicFunctoriality}) that a map $\Lambda''\to N$ is not locally constant if and only if the composed map $\Lambda \to N$ is not locally constant. It follows that the inclusion map, for each lift $\tilde\alpha''$ of $\alpha''$, 
$$\sImm(M, N)\setminus_{\tilde\alpha''} N^{\comp(\Lambda'')} \longrightarrow \sImm(M, N)\setminus_{\tilde\alpha} N^{\comp(\Lambda)}$$
is a homeomorphism. Thus, our spectrum is really equivalent to 
 $$\Gamma\left(\nlc_{\alpha''}(\Lambda'', M); \Sigma^\infty \sImm(M, N)\setminus_{\tilde\alpha''} N^{\comp(\Lambda'')}\right)$$
which is the spectrum of sections associated with the map
$$E_{\alpha''}\to \nlc_{\alpha''}(\Lambda'', M).$$
To analyze this spectrum of sections, we need to analyze the inverse images of this map. Let $\tilde\alpha''\in  \nlc_{\alpha''}(\Lambda'', M)$. Thus $\tilde\alpha''\colon \Lambda'' \to M$ is a lift of $\alpha''$. Recall that it follows from our assumptions that $\alpha''$ is not locally constant, and therefore any lift is automatically not locally constant. We need to understand the corresponding space 
$$\sImm(M, N)\setminus_{\tilde\alpha''} N^{\comp(\Lambda'')}.$$ Recall once again that $M$ is a disjoint union of open balls in a fixed vector space $W$, so two points in the same component of $M$ can be subtracted, and the difference is an element of $W$. It follows that two lifts $\tilde\alpha_1''\colon \Lambda \to M$ and $\tilde\alpha_2''\colon \Lambda \to M$ of $\alpha''\colon \Lambda \to k$ can be subtracted, and the difference can be thought of as a map from $\Lambda$ to $W$. 
\begin{definition}
Let $\tilde\alpha_1''\colon \Lambda \to M$ and $\tilde\alpha_2''\colon \Lambda \to M$ be two lifts of $\alpha''\colon \Lambda \to k$ we say that $\tilde\alpha_1''$ and $\tilde\alpha_2''$ are {\em similar}, if their difference is a locally constant map from $\Lambda$ to $W$.
\end{definition} 
Clearly, similarity is an equivalence relation on $\nlc_{\alpha''}(\Lambda'', M)$. In fact, it is a restriction of the equivalence relation given by the cosets of $W^{\comp(\Lambda'')}$ in $W^{\supp(\Lambda'')}$ (recall from Remark~\ref{R: PropertiesOfNlc} that $\nlc_{\alpha''}(\Lambda'', M)$ can be idenitfied with a subspace of $W^{\supp(\Lambda'')}$).

 It is obvious that if $\tilde\alpha_1''$ and $\tilde\alpha_2''$ are two similar lifts of $\alpha''$, then their fibers
in $E_{\alpha''}$, that is the spaces 
$$\sImm(M, N)\setminus_{\tilde\alpha_1''} N^{\comp(\Lambda'')} \mbox{ and } \sImm(M, N)\setminus_{\tilde\alpha_2''} N^{\comp(\Lambda'')}$$
are the same. That is, they are the same subspace of $\sImm(M,N)$. This observation can be extended as follows. Since every component of $M$ is a subset of $W$, we can consider $M$ as a subset of a disjoint union of $k$ copies of $W$, or in other words a subspace of $k\times W$. We may consider lifts of $\alpha''\colon \Lambda'' \to k$ to $k\times W$, including lifts that do not necessarily land in $M$. For such a generalized lift $\tilde\alpha''$, define $\sImm(M, N)\setminus_{\tilde\alpha''} N^{\comp(\Lambda'')} $ to be the space of standard immersions of $k\times W$ to $W$ that take $M$ into $N$, and such that the composed map $\Lambda'' \stackrel{\tilde\alpha''}{\longrightarrow} k\times W \longrightarrow  W$ is not locally constant. Clearly, the invariance under similarity holds in this generalized setting.

Let us say that a lift $\tilde\alpha''\colon \Lambda'' \to M$ of $\alpha''$ is very good if any two points in the support of $\Lambda''$ that go to the same component of $M$ go to the same point of $M$. Equivalently, $\tilde\alpha''$ is very good if the map $\tilde\alpha''\colon \Lambda''\to M$ factors as $\Lambda''\stackrel{\alpha''}{\to} k\stackrel{s}{\to} M$, where $s$ is a section of the quotient map $q\colon M\to k$. We say that $\tilde\alpha''$ is good if it is similar to a very good lift. We have the following lemma
\begin{lemma}\label{L: NlcImpliesGood}
Suppose $\tilde\alpha''$ is not good. Then
$$\sImm(M, N)\setminus_{\tilde\alpha''} N^{\comp(\Lambda'')}= \sImm(M,N) $$
\end{lemma}
\begin{proof}
It is not hard to show that if $\tilde\alpha''\colon \Lambda \longrightarrow M$ is not good, i.e., is not similar to a very good lift, then $\tilde\alpha''$ is similar to a lift that takes two points in the same component of $\Lambda$ to two different points in the same component of $M$. But it is obvious that if a map $\Lambda \longrightarrow M$ has this property, and we compose it with a standard immersion from $M$ to $N$ then the resulting map $\Lambda \longrightarrow N$ is not locally constant, because a standard immersion is injective on each component of $M$.
\end{proof}
Let $G\subset \nlc_{\alpha''}(\Lambda'', M)$ be the subset of good lifts of $\alpha''$. It follows from Lemma~\ref{L: NlcImpliesGood} that there is a pullback diagram
\begin{equation}\label{D: SectionSpaces}
\begin{CD}
\Gamma\left(\nlc_{\alpha''}(\Lambda'', M); \Sigma^\infty \sImm(M, N)\setminus_{\tilde\alpha''} N^{\comp(\Lambda'')}\right) @>>> \map\left(\nlc_{\alpha''}(\Lambda'', M), \Sigma^\infty \sImm(M, N)\right) \\
 @VVV  @VVV \\
\Gamma\left(G; \Sigma^\infty \sImm(M, N)\setminus_{\tilde\alpha''} N^{\comp(\Lambda'')}\right) @>>> 
\map\left(G,  \Sigma^\infty \sImm(M, N)\right)
\end{CD}
\end{equation}
where the vertical maps are restriction maps. The following lemma is an elementary exercise
\begin{lemma}
If $\alpha''\colon \Lambda \longrightarrow k$ is not locally constant, that the space 
 $G$ of good lifts of $\alpha''$ is a closed contractible subspace of $\nlc_{\alpha''}(\Lambda'', M)$. 
\end{lemma} 
 It follows that the right vertical map in~\eqref{D: SectionSpaces} is a fibration and a weak homotopy equivalence, and therefore so is the left vertical map. Our next goal is to analyze the spectrum in the lower left corner of the above diagram. Recall once again that $\alpha''$ is a map $\alpha''\colon \Lambda'' \to k$, and let $\tilde\alpha''\in G$ be a good lift of $\alpha''$ into $M$. By definition, this means that $\tilde\alpha''$ is similar to a very good lift (possibly into $k\times W$ rather than $M$), and therefore the fiber at $\tilde\alpha''$ $\sImm(M, N)\setminus_{\tilde\alpha''} N^{\comp(\Lambda'')}$ (or its suspension spectrum) is the same as the fiber at a very good $\tilde\alpha''$. Let us suppose that $\tilde\alpha''$ is very good. By definition, this means that the map $\tilde\alpha''\colon \Lambda'' \longrightarrow M$ factors as a composition 
$$\Lambda'' \stackrel{\alpha''}{\to}{k}\stackrel{s}{\hookrightarrow}  M$$
where $s$ is some section of the quotient map $q\colon M\to k$. Let $\alpha''(\Lambda'')$ be the partition of $k$ that is the image of $\alpha''$. It is easy to see that in this case, assuming that the balls making up $M$ are small enough, $\sImm(M, N)\setminus_{\tilde\alpha''} N^{\comp(\Lambda'')}$ is naturally homotopy equivalent (in fact a deformation retract of) $N^k\setminus N^{\comp(\alpha''(\Lambda''))}$. In fact, we have the following lemma, whose proof is elementary.
\begin{lemma}\label{L: TrivialOverGood}
Assuming the balls that make up $M$ are small enough, the pullback of the diagram 
$$E_{\alpha''} \longrightarrow \nlc_{\alpha''}(\Lambda'', M) \hookleftarrow G$$ is a homotopy bundle over $G$, fiber homotopy equivalent to the product space $(N^{k}\setminus N^{\comp(\alpha(\Lambda''))})\times G$.
\end{lemma}
We are now ready to prove the following proposition, which describes the sections spectrum 
$\Gamma\left(\nlc_{\alpha'}(\Lambda', M); \Sigma^\infty \sImm(M, N)\setminus_{\tilde\alpha} N^{\comp(\Lambda)}\right)$ associated to a morphism $f\colon (\Lambda, \alpha)\to (\Lambda', \alpha')$ of partitions over $k$.
\begin{proposition}\label{P: DescriptionSections}
Let $f\colon (\Lambda, \alpha)\to (\Lambda', \alpha')$ be a morphism in $\calS_k$. 
There exists an $\epsilon>0$ such that if all the connected components of $M$ have diameter less than $\epsilon$ then the following holds.

 If the map $\alpha'\colon \Lambda' \to k$ is not locally constant, there is a homotopy equivalence 
$$\Gamma\left(\nlc_{\alpha'}(\Lambda', M); \Sigma^\infty \sImm(M, N)\setminus_{\tilde\alpha} N^{\comp(\Lambda)}\right)\stackrel{\simeq}{\longrightarrow} \Sigma^\infty N^k \setminus N^{\comp(\alpha(\Lambda))}$$
where the map is given by evaluating a section at the lift of of $\alpha'$ that sends each point  $x\in \supp(\Lambda')$ to the center of the corresponding connected component of $M$.

If the map $\alpha'\colon \Lambda' \to k$ is locally constant, then the spectrum of sections equals to the mapping spectrum
$$\map_*\left(\nlc_{\alpha'}(\Lambda', M)_+, \Sigma^\infty N^k\right).$$
\end{proposition}
\begin{proof}
Suppose first that $\alpha'$ is not locally constant. We saw above that the spectrum of sections that we are interested in is equivalent to the spectrum 
$$\Gamma\left(\nlc_{\alpha''}(\Lambda'', M); \Sigma^\infty \sImm(M, N)\setminus_{\tilde\alpha} N^{\comp(\Lambda)}\right)$$
where $\Lambda''=f(\alpha)$. During the proof we used the fact that given a lift $\tilde\alpha''$ of $\alpha''$, and the corresponding lift $\tilde\alpha$ of $\alpha$, the map 
$$ \sImm(M, N)\setminus_{\tilde\alpha''} N^{\comp(\Lambda'')}\longrightarrow  \sImm(M, N)\setminus_{\tilde\alpha} N^{\comp(\Lambda)}$$
is a homeomorphism, and so we do not need to distinguish between the two. We saw that the following restriction map is an equivalence
$$\Gamma\left(\nlc_{\alpha''}(\Lambda'', M); \Sigma^\infty \sImm(M, N)\setminus_{\tilde\alpha} N^{\comp(\Lambda)}\right)\longrightarrow \Gamma\left(G; \Sigma^\infty \sImm(M, N)\setminus_{\tilde\alpha} N^{\comp(\Lambda)}\right)$$
where $G\subset \nlc_{\alpha''}(\Lambda'', M)$ is the subset of good lifts. By Lemma~\ref{L: TrivialOverGood}, the latter spectrum of sections is equivalent to $\map_*\left(G_+, \Sigma^\infty N^k\setminus N^{\comp(\alpha(\Lambda))}\right)$, which, since $G$ is contractible, is equivalent to $\Sigma^\infty N^k\setminus N^{\comp(\alpha(\Lambda))}$.

Now suppose that $\alpha'$ is locally constant. Then any lift $\tilde\alpha'\colon \Lambda' \to M$ will have the property that it will send certain two elements in the same component of $\alpha'$ to two distinct points in the same components of $M$. The same will be true for the corresponding lift $\tilde\alpha$ of $\alpha$. It follows that for any standard immersion  $h\in \sImm(M,N)$, the composed map 
$h\circ \tilde\alpha$ is not locally constant, because standard immersions are injective on components of $M$. It follows that $ \sImm(M, N)\setminus_{\tilde\alpha} N^{\comp(\Lambda)}=\sImm(M,N)\simeq N^k$ for every lift of $\tilde\alpha'$ of $\alpha'$ (the last equivalence holds if the components of $M$ are small enough), and so
$$\Gamma\left(\nlc_{\alpha'}(\Lambda', M); \Sigma^\infty \sImm(M, N)\setminus_{\tilde\alpha} N^{\comp(\Lambda)}\right)\simeq \map_*\left(\nlc_{\alpha'}(\Lambda', M)_+, \Sigma^\infty N^k\right).$$
\end{proof}

Next we proceed to analyze the homotopy limit
$$ \underset{^a\calS_k}{\holim} \Gamma\left(\nlc_{\alpha'}(\Lambda', M); \Sigma^\infty \sImm(M, N)\setminus_{\tilde\alpha} N^{\comp(\Lambda)}\right).$$
Recall that $^a\calS_k$ is the twisted arrow category, whose objects are arrows $(\Lambda,\alpha)\to (\Lambda',\alpha')$ of partitions over $k$ (as well as some arrows involving $\Omega$). We consider the following two subcategories of this twisted arrow category.
\begin{definition}
Let $\calC_1\subset {^a\calS_k}$ be the full subcategory containing all arrows of the form $(\Lambda,\alpha)\to (\Lambda',\alpha')$ where $\alpha\colon \Lambda\to k$ is not locally constant (arrows involving $\Omega$ are included).

Let $\calC_0\subset \calC_1$ be the full subcategory whose objects are all arrows where both $\alpha$ and $\alpha'$ are not locally constant.
\end{definition}
\begin{lemma} \label{L: NeedsOmega}
The functor $\Gamma\left(\nlc_{\alpha'}(\Lambda', M); \Sigma^\infty \sImm(M, N)\setminus_{\tilde\alpha} N^{\comp(\Lambda)}\right)$, considered as a covariant functor from $^a\calS_k$ to spectra, is equivalent to the homotopy right Kan extension of its restriction to $\calC_1$.
\end{lemma}
\begin{proof}
For the rest of the section, let $\Theta$ denote a generic object of $^a\calS_k$, and let $\Gamma(\Theta)$ be the corresponding sections spectrum. The content of the lemma is that for every object $\Theta$ of $^a\calS_k$, the natural map 
$$\Gamma(\Theta)\longrightarrow \underset{\Theta\to\Theta'\in \Theta\downarrow\calC_1}{\holim}{\Gamma(\Theta')}$$
is a weak equivalence. Since $\calC_1$ is a full subcategory, it is enough to consider the case when $\Theta$ is not in $\calC_1$. Suppose $\Theta$ is the arrow $(\Lambda, \alpha)\to (\Lambda', \alpha')$, where $\alpha$ is locally constant. It follows that $\alpha'$ is locally constant as well (Lemma~\ref{L: BasicFunctoriality}). By the second part of Proposition~\ref{P: DescriptionSections},  
$$\Gamma(\Theta)\simeq\map\left(\nlc_{\alpha'}(\Lambda', M), \Sigma^\infty N^k\right).$$
Now let us turn our attention to $\holim \Gamma(\Theta')$. The objects of $\Theta\downarrow \calC_1$ are diagrams of the form
$$\begin{CD}
(\Lambda, \alpha) @<<< (\Delta, \beta) \\
@VVV @VVV \\
(\Lambda', \alpha') @>>> (\Delta', \beta') 
\end{CD}$$
where $(\Delta, \beta) \to (\Delta', \beta')$ is an object of $\calC_1$, i.e., $\beta$ is not locally constant (there also are the diagrams with $\Omega$ in place of $(\Delta,\beta)$).
It is easy to see that the full subcategory of $\Theta\downarrow \calC_1$ consisting of those objects for which the morphism $\Lambda'\to \Delta'$ is the identity is an initial subcategory of $\Theta\downarrow \calC_1$, so it is enough to consider the homotopy limit over this subcategory. The objects of this subcategory can be identified simply with arrows $(\Delta,\beta) \to (\Lambda, \alpha)$, where $\beta$ is not locally constant, together with the arrow $\Omega\to (\Lambda, \alpha)$, which is a final object of this category. We are looking at the homotopy limit of spectra
$$\Gamma\left(\nlc_{\alpha'}(\Lambda', M); \Sigma^\infty \sImm(M, N)\setminus_{\tilde\beta} N^{\comp(\Delta)}\right).$$
Since $\alpha'$ is locally constant, it follows, again by the second part of Proposition~\ref{P: DescriptionSections}, that all these spectra are equivalent to 
$$\map\left(\nlc_{\alpha'}(\Lambda', M), \Sigma^\infty N^k\right),$$ the diagram of spectra $\Gamma(\Theta')$ is equivalent to a constant diagram, and the maps $\Gamma(\Theta)\to \Gamma(\Theta')$ are equivalences for all $\Theta'$. Since the underlying diagram has a final object\footnote{It is for the sake of this step that we have introduced the object $\Omega$.} and therefore has a contractible nerve, this is enough to prove the lemma.
\end{proof}
It follows from the above lemma that the restriction map 
$$\underset{\Theta\in{^a\calS_k}}{\holim} \Gamma(\Theta) \longrightarrow 
\underset{\Theta\in\calC_1}{\holim} \Gamma(\Theta)$$ 
is a homotopy equivalence. The next lemma will imply that the restriction of the homotopy limit from $\calC_1$ to $\calC_0$ is an equivalence as well.
\begin{lemma}
For every object $\Theta$ of $\calC_1$, the category $\calC_0\downarrow \Theta$ has a final object.
\end{lemma}
\begin{proof}
It is enough to consider the case when $\Theta$ is in $\calC_1$ but not in $\calC_0$. Let $\Theta$ be the arrow $(\Lambda, \alpha) \to (\Lambda', \alpha')$ where $\alpha'$ is locally constant, but $\alpha$ is not. Let $\Lambda''$ be the partition with the same support as $\Lambda$ and where $x$ and $y$ are in the same component if and only if they are in the same component of $\Lambda$ and $\alpha(x)=\alpha(y)$. Clearly, $\Lambda''$ is a refinement of $\Lambda$,  the morphism $\alpha\colon\Lambda'' \to k$ is locally constant, and the arrow $(\Lambda, \alpha) \to (\Lambda'',\alpha)$ is initial among arrows in $\calS_k$ from $(\Lambda, \alpha)$ to partitions over $k$ whose map into $k$ is locally constant. In particular, the morphism $(\Lambda, \alpha) \to (\Lambda', \alpha')$ factors as $(\Lambda, \alpha) \to (\Lambda'', \alpha) \to (\Lambda', \alpha')$. It follows that the diagram 
$$\begin{CD}
(\Lambda'', \alpha) @<<< (\Lambda, \alpha) \\
@VVV @VVV \\
(\Lambda', \alpha') @>>> (\Lambda', \alpha') 
\end{CD}$$
represents a final object of $\calC_0\downarrow \Theta$. 
\end{proof}
\begin{corollary}\label{C: Reduction}
Assuming that the connected components of $M$ are small enough, the following restriction map is a weak equivalence
$$\underset{^a\calS_k}{\holim}\Gamma(\Theta) \longrightarrow \underset{\calC_0}{\holim}\Gamma(\Theta).$$
\end{corollary}
Next, we want to analyze $\underset{\calC_0}{\holim}\Gamma(\Theta).$ By the first part of Proposition~\ref{P: DescriptionSections}, there is a natural equivalence
$$\underset{\calC_0}{\holim}\Gamma(\Theta)\stackrel{\simeq}{\longrightarrow} \underset{\calC_0}{\holim}\Sigma^\infty N^k\setminus N^{\comp(\alpha(\Lambda))}.$$
Recall once again that $\calS_k$ is the category of partitions with maps into $k$, together with an initial object $\Omega$. Let $\nlc(\calS, k)$ be the full subcategory of $\calS_k$ consisting of partitions $(\Lambda, \alpha)$ such that the map $\alpha\colon \Lambda\longrightarrow k$ is not locally constant, together with $\Omega$. By definition, $\calC_0$ is the twisted arrow category of $\nlc(\calS, k)$. 
Therefore, $\underset{\calC_0}{\holim}\Sigma^\infty N^k\setminus N^{\comp(\alpha(\Lambda))}$ can be identified with the spectrum of homotopy natural transformations from the one point functor to the functor $\Sigma^\infty N^k\setminus N^{\comp(\alpha(\Lambda))}$. In other words,
$$\underset{\calC_0}{\holim}\Sigma^\infty N^k\setminus N^{\comp(\alpha(\Lambda))}\simeq \underset{\nlc(\calS, k)}{\holim}\Sigma^\infty N^k\setminus N^{\comp(\alpha(\Lambda))}.$$

Recall that $\calMP_k$ is the category of non-discrete paritions of $k$, together with an initial object $\Omega$. There is an evident functor from $\nlc(\calS, k)$ to $\calMP_k$ sending $(\Lambda, \alpha)$ to $\alpha(\Lambda)$ and $\Omega$ to $\Omega$. This induces a map
$$\underset{\Delta\in\calMP_k}{\holim}\Sigma^\infty N^k\setminus N^{\comp(\Delta)}\simeq \underset{(\Lambda, \alpha)\in\nlc(\calS, k)}{\holim}\Sigma^\infty N^k\setminus N^{\comp(\alpha(\Lambda))}.$$
It is easy to check that the functor $\nlc(\calS, k)\longrightarrow\calMP_k$ satisfies the usual criterion for inducing a weak equivalence on homotopy limits (of contravariant functors) and thus the above map is a homotopy equivalence. We are now ready to prove the main theorem of this section.
\begin{theorem} \label{T: StandardEmbeddingsModel}
Let $M$ be the disjoint union of small enough open balls in vector space $W$. Let $N\subset W$ be a nice submanifold. There are weak equivalences, given by the natural evaluation maps
$$\Sigma^\infty\sEmb(M,N)\longrightarrow \underset{\calS_k\rtimes\nlc_{(-)}(-, M)}{\holim}  \Sigma^\infty \sImm(M, N)\setminus_{\tilde\alpha} N^{\comp(\Lambda)}$$
and
$$\Sigma^\infty\sEbar(M,N)\longrightarrow \underset{\calS_k\rtimes\nlc_{(-)}(-, M)}{\holim}  \Sigma^\infty \overline{\sImm(M, N)\setminus_{\tilde\alpha} N^{\comp(\Lambda)}}.$$
\end{theorem}
\begin{proof}
Consider the commutative square
$$\begin{CD}
\Sigma^\infty\sEmb(M,N) @>>> \underset{\calS_k\rtimes\nlc_{(-)}(-, M)}{\holim}  \Sigma^\infty \sImm(M, N)\setminus_{\tilde\alpha} N^{\comp(\Lambda)}\\
@VVV @VVV\\
\Sigma^\infty \config(k, N) @>>> \underset{\calMP_k}{\holim}  \Sigma^\infty \sImm(M, N)\setminus N^{\comp(\Lambda)}.
\end{CD}$$
where the vertical maps are given by evaluation at the centers of the components of $M$, and the horizontal maps are evaluation maps. The left vertical map is an equivalence if the components of $M$ are small enough, the bottom horizontal map is an equivalence by Proposition~\ref{P: ModelForStableConfig}, and the right vertical map was just shown to be an equivalence. This proves the first part of the theorem. The second part can be proved in exactly the same way.
\end{proof}
\section{Derivatives of the standard embeddings functor.} \label{S: DerivativesOfStandard}
Our goal in this section is to describe the derivatives (in the sense of orthogonal calculus) of the functor $N\mapsto \Sigma^\infty\sEbar(M,N)$ (see Proposition~\ref{P: MainThmForStandard} below, which amounts to a special case of Theorem~\ref{T: MainTheorem}, namely the case of standard embeddings). It is crucial that the description is natural in $M$, where $M$ ranges over the category whose objects are disjoint unions of (small enough) balls in $W$ and whose morphisms are standard embeddings. In the previous section, we obtained a natural equivalence, induced by the evident evaluation map
$$\Sigma^\infty\sEbar(M,N)\longrightarrow \underset{\calS_k\ltimes\nlc_{(-)}(-, M)}{\holim}  \Sigma^\infty \overline{\sImm(M, N)\setminus_{\tilde\alpha} N^{\comp(\Lambda)}}.$$
We remind the reader that $\overline{\sImm(M, N)\setminus_{\tilde\alpha} N^{\comp(\Lambda)}}$ stands for the homotopy fiber of the inclusion map $$\sImm(M, N)\setminus_{\tilde\alpha} N^{\comp(\Lambda)} \longrightarrow \sImm(M,N).$$ 

Let us also recall that this homotopy limit can be expressed as a homotopy limit over the twisted arrow category of $\calS_k$. To wit, the evaluation map induces a homotopy equivalence
$$\Sigma^\infty\sEbar(M,N)\longrightarrow \underset{(\Lambda,\alpha)\to(\Lambda',\alpha')\in ^a\calS_k}{\holim}\Gamma\left(\nlc_{\alpha'}(\Lambda', M), \Sigma^\infty \overline{\sImm(M, N)\setminus_{\tilde\alpha} N^{\comp(\Lambda)}}\right).$$

It is not hard to show that while the above diagram is not finite, its homotopy limit is equivalent to the homotopy limit of a finite diagram, and therefore the operators $\TP_n$ and $\D_n$ commute with the above homotopy limit. Thus we have a natural equivalence
$$\TP_n\Sigma^\infty\sEbar(M,N)\longrightarrow \underset{(\Lambda,\alpha)\to(\Lambda',\alpha')\in ^a\calS_k}{\holim}\TP_n\left(\Gamma\left(\nlc_{\alpha'}(\Lambda', M), \Sigma^\infty \overline{\sImm(M, N)\setminus_{\tilde\alpha} N^{\comp(\Lambda)}}\right)\right).$$
%

%

%
By Proposition~\ref{P: DescriptionSections} the functor $\Gamma\left(\nlc_{\alpha'}(\Lambda', M), \Sigma^\infty {\sImm(M, N)\setminus_{\tilde\alpha} N^{\comp(\Lambda)}}\right)$ is equivalent to $N^k\setminus N^{\comp(\alpha(\Lambda))}$ if $\alpha'(\Lambda')$ is not discrete (equivalently, if $\alpha'\colon \Lambda'\to k$ is not locally constant), and is equialent to $\map_*\left(\nlc_{\alpha'}(\Lambda', M), \Sigma^\infty N^k\right)$ if $\alpha'(\Lambda')$ is discrete. It follows that the spectrum $$\Gamma\left(\nlc_{\alpha'}(\Lambda', M), \Sigma^\infty \overline{\sImm(M, N)\setminus_{\tilde\alpha} N^{\comp(\Lambda)}}\right)$$ is equivalent to $$\overline{N^k\setminus N^{\comp(\alpha(\Lambda))}}:=\hofiber\left({N^k \setminus N^{\comp(\alpha(\Lambda))}}\longrightarrow N^k.\right)$$ if $\alpha'(\Lambda')$ is not discrete, and is contractible if $\alpha'(\Lambda')$ is discrete. Recall that if $\alpha'(\Lambda')$ is not discrete $\alpha(\Lambda)$ is not discrete, and in this case the functor $\overline{N^k \setminus N^{\comp(\alpha(\Lambda))}}$ is, as a functor of $N$, homogeneous of degree $\excess(\alpha(\Lambda))$ (Lemma~\ref{L: BasicHomogeneousFunctor}). Let  Let $\widetilde{\calS_k^n}$ be the full subcategory of $\calS_k$ consisting of partitions $(\Lambda, \alpha)$ over $k$ such that $\excess(\alpha(\Lambda))\le n$ (this includes those $\Lambda$ for which $\alpha(\Lambda)$ is discrete, but on the other hand does not include the initial object $\Omega$, which is no longer needed). Let $\widetilde{\calS_k^n}\ltimes \nlc_{(-)}(-, M)$ be the corresponding subcategory of $\calS_k\ltimes \nlc_{(-)}(-, M)$. Clearly, if $((\Lambda,
\alpha),\tilde\alpha)$ is an object of ${\calS_k^n}\ltimes\nlc_{(-)}(-, M)$, but not of $\widetilde{\calS_k^n}\ltimes\nlc_{(-)}(-, M)$, then $\TP_n\Sigma^\infty\overline{\sImm(M,N)\setminus_{\tilde\alpha} N^{\comp(\Lambda)}}\simeq *$. Since morphisms of partitions never  increase excess, it follows that $\widetilde{\calS_k^n}$ is a final subcategory of $\calS_k^n$. Since we are dealing with a contravariant functor, it follows in turn that the homotopy limit of $\TP_n\Sigma^\infty\overline{\sImm(M,N)\setminus_{\tilde\alpha} N^{\comp(\Lambda)}}$ does not change homotopy type if we restrict the underlying category from $\calS_k\ltimes\nlc_{(-)}(-, M)$ to $\widetilde{\calS_k^n}\ltimes\nlc_{(-)}(-, M)$. Furthermore, let $\calS_k^n$ be the full subcategory of $\calS_k$ consisting of partitions $(\Lambda, \alpha)$ over $k$ such that $\Lambda$ is not discrete, and $\excess(\Lambda)\le n$. Clearly, $\calS_k^n$ is a subcategory of $\widetilde{\calS_k^n}$, because morphisms can only lower excess. Let $\calS_k^n\ltimes \nlc_{(-)}(-, M)$ be the corresponding subcategory of $\widetilde{\calS_k}\ltimes \nlc_{(-)}(-, M)$. Once again, it is not difficult to show, using standard manipulations of  homotopy limits and what we know about our functor that the restriction map 
$$\underset{\widetilde{\calS_k^n}\ltimes\nlc_{(-)}(-, M)}{\holim}  \TP_n\Sigma^\infty \overline{\sImm(M, N)\setminus_{\tilde\alpha} N^{\comp(\Lambda)}} \longrightarrow \underset{\calS_k^n\ltimes\nlc_{(-)}(-, M)}{\holim}  \TP_n\Sigma^\infty \overline{\sImm(M, N)\setminus_{\tilde\alpha} N^{\comp(\Lambda)}}$$
is a homotopy equivalence. It follows that the evaluation map
$$
\Sigma^\infty \sEbar(M,N)\longrightarrow \underset{\calS_k^n\ltimes\nlc_{(-)}(-, M)}{\holim}  \Sigma^\infty \overline{\sImm(M, N)\setminus_{\tilde\alpha} N^{\comp(\Lambda)}}
$$
yields a model for $\TP_n\Sigma^\infty\sEmb(M,N)$, and therefore also for $\D_n \Sigma^\infty\sEmb(M,N)$ (assuming, as  usual, that the components of $M$ are small enough).
This means, in particular that there is an equivalence
$$\D_n\Sigma^\infty \sEbar(M,N)\simeq \underset{\calS_k^n\ltimes\nlc_{(-)}(-, M)}{\holim}  \D_n\Sigma^\infty \overline{\sImm(M, N)\setminus_{\tilde\alpha} N^{\comp(\Lambda)}}.$$
Moreover, this equivalence is natural with respect to standard embeddings in the variable $M$. 
Let us review what we know about $\D_n\Sigma^\infty \overline{\sImm(M, N)\setminus_{\tilde\alpha} N^{\comp(\Lambda)}}.$ Recall once again that the functor $\Sigma^\infty \overline{\sImm(M, N)\setminus_{\tilde\alpha} N^{\comp(\Lambda)}}$ is homogeneous of degree $\excess(\alpha(\Lambda))$. Thus $\D_n\Sigma^\infty \overline{\sImm(M, N)\setminus_{\tilde\alpha} N^{\comp(\Lambda)}}$ is equivalent to $\Sigma^\infty \overline{\sImm(M, N)\setminus_{\tilde\alpha} N^{\comp(\Lambda)}}$ if $\excess(\alpha(\Lambda))=n$, and is contractible if $\excess(\alpha(\Lambda))\ne n$. Since in our category $\calS_k^n$, the objects $(\Lambda, \alpha)$ satisfy $\excess(\Lambda)\le n$, it follows that $\D_n\Sigma^\infty \overline{\sImm(M, N)\setminus_{\tilde\alpha} N^{\comp(\Lambda)}}$ is contractible, unless $\excess(\Lambda)=n$, and furthermore $\excess(\alpha(\Lambda))=n$, i.e., the fusion $\Lambda\to \alpha(\Lambda)$ is a strict fusion. Let us define a new contravariant functor on the category $\calS_k^n\ltimes \nlc_{(-)}(-, M)$ by the following formula
$$
F((\Lambda, \alpha),\tilde\alpha)=\left\{ 
\begin{array}{cl}
\Sigma^\infty \overline{\sImm(M, N)\setminus_{\tilde\alpha} N^{\comp(\Lambda)}} & \mbox{if }\excess(\Lambda) =\excess(\alpha(\Lambda))=n \\
* & \mbox{otherwise} 
\end{array}\right. 
$$
It is easy to see that the full subcategory of $\calS_k^n\ltimes \nlc_{(-)}(-, M)$ consisting of objects $((\Lambda, \alpha),\tilde\alpha)$ satisfying $\excess(\Lambda)=\excess(\alpha(\Lambda))=n$ has the property that there are no morphisms from outside this subcategory into this subcategory. Therefore, one can define $F$ on morphisms to be the restriction of the functor $\Sigma^\infty \overline{\sImm(M, N)\setminus_{\tilde\alpha} N^{\comp(\Lambda)}}$, and moreover there is a map, natural in $M$ and $N$,
$$\underset{\calS_k^n\ltimes\nlc_{(-)}(-, M)}{\holim}  F((\Lambda, \alpha),\tilde\alpha) \longrightarrow \underset{\calS_k^n\ltimes\nlc_{(-)}(-, M)}{\holim}  \Sigma^\infty \overline{\sImm(M, N)\setminus_{\tilde\alpha} N^{\comp(\Lambda)}}$$
 which induces an equivlalence on $D_n$. Thus, the source of this natural transformation provides a model for $\D_n\Sigma^\infty\sEbar(M, N)$, which is functorial in $M$.
 
The next step will show that the functoriality $\D_n\Sigma^\infty\sEbar(M, N)$ in the variable $M$ is in some sense much more robust than the functoriality of $\TP_n\Sigma^\infty\sEbar(M, N)$. While the latter functor is only functorial with respect to standard embeddings in the variable $M$, the former is functorial with respect to arbitrary maps in the variable $M$. Recall that $\nlc_\alpha(\Lambda, M)$ is the space of non locally constant lifts of $\alpha\colon\Lambda\to k$ to $M$. Let $\map_{\alpha}(\Lambda, M)$ be the space of all lifts of $\alpha$ to $M$. Thus $\map_\alpha(\Lambda, M)$ is the space of all maps $\tilde\alpha\colon \Lambda\to M$ such that $q\circ \tilde\alpha =\alpha$. There is a natural inclusion 
$$\nlc_\alpha(\Lambda, M)\hookrightarrow \map_\alpha(\Lambda, M)$$
which is a homeomorphism if $\alpha$ is itself not locally constant. Let $\calS_k^n\ltimes \map_{(-)}(-, M)$ be the corresponding topological category. It contains $\calS_k^n\ltimes \nlc_{(-)}(-, M)$ as a subcategory. We claim that the functor $F$ from $\calS_k^n\ltimes \nlc_{(-)}(-, M)$ to the category of spectra that we defined above extends to a functor on $\calS_k^n\ltimes \map_{(-)}(-, M)$ by defining 
$F((\Lambda, \alpha),\tilde\alpha)=*$ whenever $\tilde\alpha$ is in $\map_\alpha(\Lambda, M)\setminus \nlc_\alpha(\Lambda, M)$, i.e., whenever $\tilde\alpha$ is locally constant. This is well defined, because if $\tilde\alpha$ is locally constant then $\alpha$ is locally constant, and it follows that $F$ is already defined to be $*$ on the entire connected component of $((\Lambda,\alpha),\tilde\alpha)$. Moreover, if 
$$((\Lambda,\alpha),\tilde\alpha)\longrightarrow ((\Lambda',\alpha'),\tilde\alpha')$$ is a morphism in 
$\calS_k^n\ltimes \map_{(-)}(-, M)$, and $\tilde\alpha$ is locally constant, then so is $\tilde\alpha'$, which says that there is no problem defining the extension of $F$ on morphisms. Moreover we claim that the restriction map 
$$\underset{\calS_k^n\ltimes\map_{(-)}(-, M)}{\holim}  F((\Lambda, \alpha),\tilde\alpha)\longrightarrow \underset{\calS_k^n\ltimes\nlc_{(-)}(-, M)}{\holim}  F((\Lambda, \alpha),\tilde\alpha)$$ is an equivalence (in fact, an isomorphism). One can see this by rewriting the homotopy limit as a homotopy limit of sections spectra over the discrete category $^a\calS_n^n$. It is easy to see that all spaces of sections involving a locally constant $\alpha$, i.e., in all places in the diagrams where there might be a difference, are equal to $*$ on both sides.

We conclude that the homotopy limit
$$\underset{\calS_k^n\ltimes\map_{(-)}(-, M)}{\holim}  F((\Lambda, \alpha),\tilde\alpha)$$
provides a good model for $\D_n\Sigma^\infty\sEmb(M, N)$. By this we mean that there is a zig-zag of maps, natural in $M$ and $N$, connecting $\Sigma^\infty\sEmb(M, N)$ and the homotopy limit above, where all the maps induce an equivalence on the $n$-th layer $\D_n$.

By definition $F((\Lambda, \alpha),\tilde\alpha)$ is either $\Sigma^\infty \overline{\sImm(M, N)\setminus_{\tilde\alpha} N^{\comp(\Lambda)}}$ or $*$. Recall that $\sImm(M, N)\setminus_{\tilde\alpha} N^{\comp(\Lambda)}$ is the space of standard immersions from $M$ to $N$ whose composition with $\tilde\alpha\colon \Lambda\to M$ is not locally constant. It follows that there is a natural evaluation map
$$\sImm(M, N)\setminus_{\tilde\alpha} N^{\comp(\Lambda)}\longrightarrow N^\Lambda\setminus N^{\comp(\Lambda)}$$
and therefore also a natural map
$$\Sigma^\infty\overline{\sImm(M, N)\setminus_{\tilde\alpha} N^{\comp(\Lambda)}}\longrightarrow \Sigma^\infty\overline{N^\Lambda\setminus N^{\comp(\Lambda)}}.$$
Let us define a contravariant functor
$$G\colon \calS^n_k\ltimes \map_{(-)}(-, M)\longrightarrow \operatorname{Spectra}$$ by the formula
$$G((\Lambda, \alpha), \tilde\alpha)=\left\{
\begin{array}{cl}
\Sigma^\infty\overline{N^\Lambda\setminus N^{\comp(\Lambda)}} & \mbox{if } \excess(\Lambda)=\excess(\tilde\alpha(\Lambda))=n \\
* &\mbox{otherwise}
\end{array}\right. $$
By the discussion above, there is a natural transformation of functors $$F((\Lambda, \alpha), \tilde\alpha)\longrightarrow G((\Lambda, \alpha), \tilde\alpha).$$ This map is an equivalence, except in the case when $\alpha$ is locally constant while $\tilde\alpha$ is not. But, it is easy to see that for a locally constant map $\alpha\colon\Lambda \to k$, the space of locally constant lifts of $\alpha$ is a contractible subspace of the (also contractible) space of all lifts of $\alpha$. If follows that the natural transformation $F\to G$ induces an equivalence of homotopy limits
$$\underset{\calS_k^n\ltimes\map_{(-)}(-, M)}{\holim}  F((\Lambda, \alpha),\tilde\alpha) \longrightarrow \underset{\calS_k^n\ltimes\map_{(-)}(-, M)}{\holim}  G((\Lambda, \alpha),\tilde\alpha)$$
and so the target of this map provides a good model for $\D_n\Sigma^\infty\sEbar(M,N)$.

Let us now reconsider the category $\calS_k^n\ltimes\map_{(-)}(-, M)$. We have been writing a generic object in this category in the form $((\Lambda, \alpha),\tilde\alpha)$. But clearly, $\alpha$ is redundant, since it is determined by $\tilde\alpha$. In other words, the category $\calS_k^n\ltimes\map_{(-)}(-, M)$ is isomorphic to the category $\calS^n\ltimes\map(-, M)$. An object in this category is a pair $(\Lambda, \tilde\alpha)$, where $\Lambda$ is a non-discrete partition, and $\tilde\alpha$ is a map $\tilde\alpha\colon \Lambda\to M$. So, we can think of $G$ as a functor on the category $\calS^n\ltimes\map(-, M)$, and we will write $G(\Lambda,\tilde\alpha)$ instead of $G((\Lambda, \alpha),\tilde\alpha)$. We remind the reader once again that $G(\Lambda, \tilde\alpha)$ is $\Sigma^\infty \overline{N^\Lambda\setminus N^{\comp(\Lambda)}}$ if $\excess(\Lambda)=\excess(\tilde\alpha(\Lambda))$, and is $*$ otherwise. Thus, our current model for $\D_n\Sigma^\infty\sEbar(M, N)$ is given by the homotopy limit 
$$\underset{\calS^n\ltimes\map(-, M)}{\holim}  G(\Lambda, \tilde\alpha).$$
Next we want to rewrite this homotopy limit in a way that fits better with the main theorem. Recall that $\calE_n$ is the category of irreducible partitions of excess $n$ and strict fusions between them. Let $\widetilde\calP(\Lambda):=\calP(\Lambda)\setminus \hat1$ be the poset of non-discrete refinements of $\Lambda$. Recall that the assignment $\Lambda\mapsto \widetilde\calP(\Lambda)$ defines a covariant functor from $\calE_n$ to the category of small categories (this follows from Lemma~\ref{L: PreserveStrictRefinements}). Let $\calE_n\ltimes \widetilde\calP(-)$ be the associated Grothendieck cathegory. We will denote a generic object of $\calE_n\ltimes \widetilde\calP(-)$ by $(\Lambda, \Lambda')$, where $\Lambda$ is an object of $\calE_n$, and $\Lambda'$ is a non-discrete refinement of $\Lambda$. There is a forgetful functor $\calE_n\ltimes \widetilde\calP(-)\longrightarrow \calS_n$ that sends $(\Lambda, \Lambda')$ to $\Lambda'$. Let $\calE_n\ltimes \widetilde\calP(-)\ltimes \map(-, M)$ be the category whose objects are triples $(\Lambda, \Lambda', \tilde\alpha)$ where $\Lambda$ and $\Lambda'$ are the same as before, and $\tilde\alpha$ is a map from the support of $\Lambda$ (or of $\Lambda'$, since the two partitions have the same support) to $M$. There is an evident forgetful functor $\calE_n\ltimes \widetilde\calP(-)\ltimes \map(-, M)\longrightarrow \calS_n\ltimes\map(-, M)$ that sends $(\Lambda, \Lambda', \tilde\alpha)$ to $(\Lambda', \tilde\alpha)$. Therefore we have a map of homotopy limits
$$\underset{(\Lambda,\tilde\alpha)\in\calS^n\ltimes\map(-, M)}{\holim}  G(\Lambda, \tilde\alpha)\longrightarrow \underset{(\Lambda,\Lambda',\tilde\alpha)\in \calE_n\ltimes\widetilde\calP(-)\ltimes\map(-, M)}{\holim}  G(\Lambda', \tilde\alpha).$$
The following lemma is another exercise in manipulating homotopy limit. If one feels uncomfortable with working on the level of topological categories, one can reformulate and prove it in terms of functors on discrete categories $\calS^n$ and $\calE_n\ltimes\widetilde\calP(-)$.
\begin{lemma}
The functor $G$, considered as a functor on $\calS^n\ltimes\map(-, M)$ is equivalent to the homotopy right Kan extension of $G$ as a functor on $\calE_n\ltimes\widetilde\calP(-)\ltimes\map(-, M)$. Therefore, the above map of homotopy limits is a homotopy equivalence. 
\end{lemma}
It follows that $$ \underset{\calE_n\ltimes\widetilde\calP(-)\ltimes\map(-, M)}{\holim}  G(\Lambda', \tilde\alpha)$$
is a good model for $\D_n\Sigma^\infty\sEbar(M, N)$. We want to analyze this homotopy limit. To begin with, note again that the functor $\map(-, M)$ may be considered a functor on $\calE_n$, rather than $\calE_n\ltimes\widetilde\calP(-)$, because it only depends on the support of the partition, which is the same for $\Lambda$ and $\Lambda'$. It follows that the category $\calE_n\ltimes\widetilde\calP(-)\ltimes\map(-, M)$ may also be written as $\left(\calE_n\ltimes\map(-, M)\right)\ltimes\widetilde\calP(-)$, and there is a natural equivalence
$$ \underset{\calE_n\ltimes\widetilde\calP(-)\ltimes\map(-, M)}{\holim}  G(\Lambda', \tilde\alpha)\simeq
 \underset{(\Lambda,\tilde\alpha)\in \calE_n\ltimes\map(-, M)}{\holim} \left( \underset{\Lambda'\in\widetilde\calP(\Lambda)}{\holim} G(\Lambda', \tilde\alpha)\right).$$
Let us now switch back to $\widetilde\nat$ notation, and also write $M^\Lambda$ instead of $\map(\Lambda, M)$. Thus the above homotopy limit can be written as
$$\underset{\Lambda\in \calE_n}{\widetilde\nat}\left(M^\Lambda; \left( \underset{\Lambda'\in\widetilde\calP(\Lambda)}{\holim} G(\Lambda', \tilde\alpha)\right)\right).$$
Recall that $G(\Lambda', \tilde\alpha)=*$ if $\excess(\tilde\alpha(\Lambda'))\ne n$. In particular, if the morphism $\Lambda\to \tilde\alpha(\Lambda)$ is not a strict fusion, then $\excess(\tilde\alpha(\Lambda'))\le\excess(\tilde\alpha(\Lambda) < \excess(\Lambda)=n$, and $G(\Lambda', \tilde\alpha)=*$. To say that the morphism $\Lambda\to \tilde\alpha(\Lambda)$ is not a strict fusion is equivalent to saying that $\tilde\alpha\in \Delta^\Lambda M$. Therefore, the above spectrum of twisted natural transformations can be written in the following form
$$\underset{\Lambda\in \calE_n}{\widetilde\nat}\left((M^\Lambda, \Delta^\Lambda M); \left( \underset{\Lambda'\in\widetilde\calP(\Lambda)}{\holim} G(\Lambda', \tilde\alpha)\right)\right).$$
Let us analyze $\underset{\Lambda'\in\widetilde\calP(\Lambda)}{\holim} G(\Lambda', \tilde\alpha)$. Recall that $\widetilde\calP(\Lambda)$ is a category with an initial object (namely $\Lambda$), and $G$ is the contravariant functor whose value is $\Sigma^\infty \overline{N^\Lambda\setminus N^{\comp(\Lambda)}}$ on the initial object, and is $*$ on all other objects. Recall that $$\Sigma^\infty \overline{N^\Lambda\setminus N^{\comp(\Lambda)}}\simeq \Omega\Sigma^\infty \Omega^\Lambda N_+\wedge S^{d\excess(\Lambda)}$$
 It follows that there are natural equivalences, where the last one relies on Lemma~\ref{L: AlternativeForT} 
\begin{multline*}
\underset{\Lambda'\in\widetilde\calP(\Lambda)}{\holim} G(\Lambda', \tilde\alpha)\simeq\map_*(|\widetilde\calP(\Lambda)|/|\widetilde\calP(\Lambda)\setminus\{\Lambda\}|, \Omega\Sigma^\infty \Omega^\Lambda N_+\wedge S^{d\excess(\Lambda)})\simeq \\ \simeq \map_*(T_\Lambda, \Sigma^\infty \Omega^\Lambda N_+\wedge S^{d\excess(\Lambda)}).\end{multline*}
In conclusion, we have obtained the following proposition, which is the main result of this section.
\begin{proposition}\label{P: MainThmForStandard}
Let $M$ be a disjoint union of balls in $W$. Let $N$ be a nice open submanifold of $W$. Assuming that the balls making up $M$ are small enough, there is a chain of weak equivalences, natural both in $M$ and $N$, connecting $\D_n\Sigma^\infty \sEbar(M, N)$ and 
$$\underset{\Lambda\in \calE_n}{\widetilde\nat}\left((M^\Lambda, \Delta^\Lambda M);  \map_*\left(T_\Lambda, \Sigma^\infty \Omega^\Lambda N_+\wedge S^{d\excess(\Lambda)}\right)\right).$$
\end{proposition}

\section{Proof of the main theorem}\label{S: EndOfProof}
In this section we prove Theorem~\ref{T: MainTheorem}. Proposition~\ref{P: MainThmForStandard} essentially gives the theorem in the special case when $M$ is a disjoint union of standard balls. It remains to deduce the theorem for general manifolds $M$, using the ideas of embedding calculus~\cite{WeissEmb}. But first, we need to analyze our main construction a little further. In particular, we want to prove Proposition~\ref{P: EmbeddingOfOrthogonal} from the introduction.

\begin{lemma}\label{L: EssentiallyCofibrant}
The functor $ \left(M^\Lambda,  \Delta^\Lambda M\right)$ is essentially cofibrant (in the sense of Definition~\ref{D: EssentiallyCofibrant}) as a covariant functor from $\calE_n^{\op}$ to pairs of spaces.
\end{lemma}
\begin{proof}
Let $\Lambda$ be an object of $\calE_n^i\setminus \calE_n^{i-1}$ for some $n$ and $i$. Thus $\Lambda$ is an irreducible partition of excess $n$, and support of size $n+i$. Our task is to analyze the pushout of the diagram
$$\underset{(\Lambda\to \Lambda')\in \Lambda\downarrow \calE_n^{i-1}}{\operatorname{colim }} M^{\Lambda'} \longleftarrow \underset{(\Lambda\to \Lambda')\in \Lambda\downarrow \calE_n^{i-1}}{\operatorname{colim }} \Delta^{\Lambda'} M \longrightarrow \Delta^\Lambda M$$
and the natural map from this pushout to $M^\Lambda$. Identifying $M^\Lambda$ with $M^{n+i}$, we claim that the above pushout maps homeomorphically onto $\Delta^{n+i}M$ (the full fat diagonal). We leave the proof of this claim as an easy exercise for the reader. It follows that the map from the pushout to $M^\Lambda$ is a cofibration. Clearly, the group of automorphisms of $\Lambda$, being a subgroup of the symmetric group $\Sigma_{n+i}$, acts freely on the quotient of $M^{n+i}$ by the full fat diagonal. This is what we needed to check.
\end{proof}

Proposition~\ref{P: EmbeddingOfOrthogonal} now follows from Lemma~\ref{L: EssentiallyCofibrant} together with Corollary~\ref{C: InductiveModel} and induction. As pointed out in the introduction (Remark~\ref{R: Consequences}) it follows from Proposition~\ref{P: EmbeddingOfOrthogonal} that 
our model for $\D_n\Sigma^\infty\Ebar(M, N)$ is a homogeneous functor of degree $n$ in the variable $N$, and also is a polynomial functor for degree $2n$ in the variable $M$ (in the sense of Embedding Calculus). 

Now, the proof of the Main Theorem. To begin with, it is enough to prove the theorem in the case when $N$ is a nice submanifold of some Euclidean space $W$. Indeed, if $N$ is a tame stably parallelizable manifold, then for some Euclidean space $W_0$ $N\times W_0$ is a nice submanifold of some $W$, so we would know that the theorem is true for $N\times W_0$. But it follows from the classification of homogeneous functors that the value of a homogeneous functor at ${\mathbb R}\sp{0}$ is determined by its value at $W_0$, for any Euclidean space $W_0$. So, let us assume that $N$ is a nice codimension zero submanifold of $W$. We may also assume that the basepoint map $h\colon M\longrightarrow N$ that is used to define $\Ebar(M, N)$ is a codimension zero embedding as well. This is so because every immersion becomes isotopic to an embedding after crossing the target with a large enough Euclidean space (and so we may assume that $h$ is an embedding), and because $\Ebar(M, N)$ does not change its homotopy type if $M$ is replaced with a tubular neighborhood.

Let $\calO^s(M)$ be the poset/category whose objects are subsets of $M$ that are disjoint unions of standard balls in $M$ (this is well-defined, because $M$ has been identified with an open subset of $W$). Let $\calO^s_k(M)$ be the sub-poset consisting of sets that are the disjoint union of at most $k$ balls. We will often want to only consider unions of balls whose radius is smaller than a fixed $\epsilon>0$, but we will not make $\epsilon$ visible in the notation. For the purposes of this paper, we define $\operatorname{T}_k\Sigma^\infty\Ebar(M, N)$ (the $k$-th Taylor approximation of $\Sigma^\infty\Ebar(M, N)$ in the sense of embedding calculus) by the formula
$$\operatorname{T}_k\Sigma^\infty\Ebar(M, N) :=\underset{U\in \calO^s_k(M)}{\holim} \Sigma^\infty\Ebar(U, N).$$
This definition differs from Weiss' original definition in~\cite{WeissEmb} in that he uses the category $\calO_k(M)$ of all subsets of $M$ that are homeomorphic to a disjoint union of open balls, rather than the smaller category of subsets that are actually the union of standard balls, but it is not hard to show that the restriction map from Weiss' definition to our definition is a homotopy equivalence. Next, we may replace $\Ebar(U, N)$ with $\sEbar(U, N)$, because there is a bi-natural homotopy equivalence $\sEbar(U, N)\stackrel{\simeq}{\longrightarrow} \Ebar(U, N)$. Thus, we may take our model for $\operatorname{T}_k\Sigma^\infty\Ebar(M, N)$ to be the homotopy limit
$$\underset{U\in \calO^s_k(M)}{\holim} \Sigma^\infty\sEbar(U, N).$$ It is not hard to show that the above homotopy limit can be used to provide a good model for $\D_n\Sigma^\infty\Ebar(M, N)$. More precisely, for $k$ large enough (in fact, $k$ needs to be roughly $2n$), there is a natural equivalence (\cite{ALV}, Proposition 10.9)
$$\D_n\Sigma^\infty\Ebar(M, N)\simeq \underset{U\in \calO^s_k(M)}{\holim} \D_n \Sigma^\infty\sEbar(U, N).$$
By Proposition~\ref{P: MainThmForStandard}, the right and side is equivalent to 
$$\underset{U\in \calO^s_k(M)}{\holim}\,\underset{\Lambda\in \calE_n}{\widetilde\nat}\left((U^\Lambda, \Delta^\Lambda U);  \map_*\left(T_\Lambda, \Sigma^\infty \Omega^\Lambda N_+\wedge S^{d\excess(\Lambda)}\right)\right).$$
By definition, this is is the same as 
$$\operatorname{T}_k \underset{\Lambda\in \calE_n}{\widetilde\nat}\left((M^\Lambda, \Delta^\Lambda M);  \map_*\left(T_\Lambda, \Sigma^\infty \Omega^\Lambda N_+\wedge S^{d\excess(\Lambda)}\right)\right).$$
Clearly, this is equivalent to 
$$\underset{\Lambda\in \calE_n}{\widetilde\nat}\left((M^\Lambda, \Delta^\Lambda M);  \map_*\left(T_\Lambda, \Sigma^\infty \Omega^\Lambda N_+\wedge S^{d\excess(\Lambda)}\right)\right)$$
for $k\ge 2n$, because the latter is a functor of degree $2n$ in the variable $M$. This completes the
proof.

\end{document}